\documentclass[final]{dmtcs-episciences}
\pdfoutput=1
\usepackage[round]{natbib}
\usepackage[utf8]{inputenc}
\usepackage[T1]{fontenc}

\usepackage{mathtools}
\usepackage{amsmath,amssymb,amsbsy}
\usepackage{graphicx}
\usepackage{multirow}
\usepackage{amsfonts}
\usepackage{amsthm}
\usepackage{amssymb}
\usepackage{graphicx, enumerate}
\usepackage{amsmath}
\usepackage{latexsym}
\usepackage{longtable}
\usepackage{tabularx}
\usepackage{amsmath}
\usepackage{amsfonts}
\usepackage{amsthm}
\usepackage{setspace}
\usepackage{graphicx}
\usepackage{float}
\usepackage{rotating}
\usepackage{tikz}
\usepackage{verbatim}
\usepackage{appendix}

\PassOptionsToPackage{unicode}{hyperref}
\PassOptionsToPackage{naturalnames}{hyperref}

\newtheorem{thm}{Theorem}[section]
\newtheorem{lem}[thm]{Lemma}

\newtheorem{cor}[thm]{Corollary}
\newtheorem{prob}[thm]{Problem}
\newtheorem{dfn}[thm]{Definition}

\newtheorem{conj}[thm]{Conjecture}

\makeatletter
\def\imod#1{\allowbreak\mkern10mu({\operator@font mod}\,\,#1)}
\makeatother

\def \bsk {\bigskip}

\newcommand{\zet}{\mathbb{Z}}

\author[Sylwia Cichacz, Karol Suchan]{Sylwia Cichacz\affiliationmark{1} \and Karol Suchan\affiliationmark{1,2}}
\title[Zero-sum partitions of Abelian groups of order $2^n$]{Zero-sum partitions of Abelian groups of order $2^n$}\affiliation{
  AGH University of Science and Technology, Krakow, Poland\\
Universidad Diego Portales,  Santiago, Chile}

\begin{document}
\publicationdata{vol. 25:1}{2023}{6}{10.46298/dmtcs.9914}{2022-08-11; 2022-08-11; 2022-12-16; 2023-02-02}{2023-02-02}
\maketitle
\begin{abstract}
The following problem has been known since the 80's. Let $\Gamma$ be an Abelian group of order $m$ (denoted $|\Gamma|=m$), and let $t$ and $m_i$, $1 \leq i \leq t$, be positive integers  such that $\sum_{i=1}^t m_i=m-1$. Determine when $\Gamma^*=\Gamma\setminus\{0\}$, the set of non-zero elements of $\Gamma$, can be partitioned into disjoint subsets $S_i$, $1 \leq i \leq t$, such that $|S_i|=m_i$ and $\sum_{s\in S_i}s=0$ for every $i\in[1,t]$.

It is easy to check that $m_i\geq 2$ (for every $i\in[1,t]$) and $|I(\Gamma)|\neq 1$ are necessary conditions for the existence of such partitions, where $I(\Gamma)$ is the set of involutions of $\Gamma$. 
It was proved that the condition $m_i\geq 2$ is sufficient if and only if $|I(\Gamma)|\in\{0,3\}$ (see Zeng, (2015)). 

For other groups (i.e., for which $|I(\Gamma)|\neq 3$ and $|I(\Gamma)|>1$), only the case of any group $\Gamma$ with $\Gamma\cong(\zet_2)^n$ for some positive integer $n$ has been analyzed completely so far, and it was shown independently by several authors that $m_i\geq 3$ is sufficient in this case. Moreover, recently Cichacz and Tuza (2021) proved that, if $|\Gamma|$ is large enough and $|I(\Gamma)|>1$, then $m_i\geq 4$ is sufficient.

In this paper we generalize this result for every Abelian group of order $2^n$. Namely, we show that the condition $m_i\geq 3$ is sufficient for $\Gamma$ such that $|I(\Gamma)|>1$ and $|\Gamma|=2^n$, for every positive integer $n$.
We also present some applications of this result to graph magic- and anti-magic-type labelings.
\end{abstract}
%05E16, 20K01, 05C25, 05C78
%\vfil

 \section{Introduction}\label{sec:intro}
 
 \subsection{Main problem}\label{subsec:problem}

Let $\Gamma$ be an Abelian group of order $m$ with the operation denoted by $+$\footnote{For standard notions, notations and results in finite algebra the reader is referred to the textbook by \citet{Gallian}. In this article we recall only the ones that are most directly related to the presentation of our work.}.  For convenience, we will denote $\sum_{i=1}^k a$ by $ka$, the inverse of $a$ by $-a$, and $a+(-b)$ by $a - b$. Moreover, we will write $\sum_{a\in S}{a}$ for the sum of all elements in $S$. The identity element of $\Gamma$ will be denoted by $0$, and the set of non-zero elements of  $\Gamma$ by $ \Gamma^*$.
Recall that any element $\iota\in\Gamma$ of order 2 (i.e., $\iota\neq 0$ and $2\iota=0$) is called an \emph{involution}. We will write $I(\Gamma)$ for the set of involutions of $\Gamma$. A non-trivial finite group has an involution if and only if the order of the group is even. 
The fundamental theorem of finite Abelian groups states that every finite Abelian group $\Gamma$ is isomorphic to the direct product of some cyclic subgroups of prime-power order \citep{Gallian}. In other words, there exists a positive integer $k$, (not necessarily distinct) prime numbers $\{p_i\}_{i=1}^{k}$, and positive integers $\{\alpha_i\}_{i=1}^{k}$, such that
$$\Gamma\cong\zet_{p_1^{\alpha_1}}\times\zet_{p_2^{\alpha_2}}\times\ldots\times\zet_{p_k^{\alpha_k}} \mathrm{, where}\;\;\; n = p_1^{\alpha_1}\cdot p_2^{\alpha_2}\cdot\ldots\cdot p_k^{\alpha_k}.$$ Moreover, this factorization is unique (up to the order of terms in the direct product).  Since any cyclic finite group of even order has exactly one involution, if $e$ is the number of cyclic subgroups in the factorization of $\Gamma$ whose order is even, then $|I(\Gamma)|=2^e-1$.

For positive integers $a$ and $b$ such that $a<b$ let $[a,b]=\{a,a+1,\ldots,b\}$.

Because the results presented in this paper are invariant under the isomorphism between groups ($\cong$), we only need to consider one group in every isomorphism class. Our presentation will be focused on groups of the form $\zet_{p_1^{\alpha_1}}\times\zet_{p_2^{\alpha_2}}\times\ldots\times\zet_{p_k^{\alpha_k}}$, with prime numbers $\{p_i\}_{i=1}^{k}$ and positive integers $\{\alpha_i\}_{i=1}^{k}$. For ease of presentation, when speaking of a subgroup of a group in which all elements have $0$ on some positions (corresponding to some terms in the factorization), we will omit these zeros (under group isomorphism). For an Abelian group $\Gamma \cong U \times H$, given a pair $u \in U$ and $v \in H$, we will use the notation $(u,v)$ (concatenation of $u$ and $v$) for the corresponding element of $\Gamma$.

The sum of all elements of a group $\Gamma$ is equal to the sum of its involutions and the identity element \citep{Gallian}. {The following lemma is well known. The readers can consult \citet{ref:ComNelPal} for a proof.}

\begin{lem}[\citet{ref:ComNelPal}]\label{involutions} Let $\Gamma$  be an Abelian group.
\begin{itemize}
 \item[-] If $|I(\Gamma)|=1$, then $\sum_{g\in \Gamma}g= \iota$, where $\iota$ is the involution.
\item[-] If $I|(\Gamma)|\neq 1$, then $\sum_{g\in \Gamma}g=0$.
\end{itemize}
\end{lem}

%The number of involutions in an Abelian group is $2^{r}-1$, where $r$ is the minimal number of generators of the Sylow $2$-subgroup. x It is well known that $\sum_{g\in\Gamma}g=0$ if and only if $|I(\Gamma)| \neq 1$ (see e.g. \citet{ref:ComNelPal}, Lemma 8).

%In \citet{Tannenbaum1} Tannenbaum introduced the notion of zero-sum partitions of in Abelian groups.
%Let $\Gamma$ be an Abelian group of order $n$ let $r_1, r_2,\ldots, r_t$ be positive integers greater than two such that $r_1 + r_2 + \ldots + r_t=n-1$ .% and let
 %$w_1, w_2,\ldots, w_t$ be arbitrary elements of $\Gamma$ (not necessary distinct).\\
 %If the non-zero elements of $\Gamma$ be partitioned into disjoint subsets $S_1,\ldots,S_t$ such that $|S_i|=r_i$ and $\sum\limits_{s\in S_i}s=0$ for every $i=1,\ldots,t$, then we say that  $\Gamma$ has the \emph{zero-sum-partition} (ZSP).

In 1981 Tannenbaum introduced the following problem of partitioning in Abelian groups.
\begin{prob}[\citet{Tannenbaum1}]\label{problemT}
Let $\Gamma$ be an Abelian group of order $m$. Let $t$ and $m_i$, {$i\in[1,t]$}, be positive integers  such that $\sum_{i=1}^t m_i=m-1$. In other words, $\{m_i\}_{i=1}^t$ is an integer partition of $m-1$. Let $w_i$, {$i\in[1,t]$}, be arbitrary elements of $\Gamma$ (not necessarily distinct). Determine when there exists a subset partition $\{S_i\}_{i=1}^t$ of $\Gamma^*$ (i.e., subsets $S_i$ are pairwise disjoint and their union is $ \Gamma^*$) such that $|S_i|=m_i$ and $\sum_{s\in S_i}s=w_i$ for every {$i\in[1,t]$}.\end{prob}

If a respective subset partition of $\Gamma^*$ exists, then we say that $\{m_i\}_{i=1}^t$ is \textit{realizable} in $\Gamma^*$ with $\{w_i\}_{i=1}^{t}$. In such a case, we say that $\{S_i\}_{i=1}^t$ \textit{realizes} $\{m_i\}_{i=1}^t$ in $\Gamma^*$ with $\{w_i\}_{i=1}^{t}$. We will not {specify} the sequence $\{w_i\}_{i=1}^{t}$ whenever it is clear from the context. Realizability of the sequence implies that $\sum_{i=1}^t w_i=\sum_{g\in \Gamma}g$. Analogously, in the context of any subset $Z$ of $\Gamma$, we will say that $\{S_i\}_{i=1}^t$, a partition of $Z$ into zero-sum subsets, realizes $\{m_i\}_{i=1}^t$ in $Z$ if $|S_i|=m_i$ for every {$i\in[1,t]$}.

Note that we use sequences $\{m_i\}_{i=1}^t$ only for ease of presentation, whereas we interpret them as multisets. In other words, only the values present in the sequence and their multiplicities are relevant for a sequence to be realizable, not their order. 

Since we will consider subsets of fixed cardinalities throughout the paper, let us present an abbreviated notation. Given a set, any of its subsets of cardinality $k$ is called a \textit{$k$-subset}.

The case most studied in the literature is the {\em Zero-Sum Partition (ZSP)} problem, i.e., when $w_i=0$ for every  {$i\in[1,t]$}. Note that,
by Lemma \ref{involutions}, for an Abelian group $\Gamma$ with $|I(\Gamma)|=1$, $\Gamma^*$ does not admit zero-sum partitions. Moreover, given a group $\Gamma$, if $\{m_i\}_{i=1}^t$ is realizable in $\Gamma^*$ (with $w_i=0$ for every {$i\in[1,t]$}), then necessarily $m_i\geq 2$ for every {$i\in[1,t]$}. It was proved that this condition is sufficient when $|I(\Gamma)|=0$  \citep{Tannenbaum1,Zeng} or $|I(\Gamma)|=3$ \citep{Zeng}. Let us generalize this situation with the following definition.

\begin{dfn}Let $\Gamma$ be a finite Abelian group of order $m$. We say that $\Gamma$ has {\em $x$-Zero-Sum Partition Property ($x$-ZSPP)} if, for every positive integer $t$ and every integer partition $\{m_i\}_{i=1}^t$ of $m-1$ (i.e., $\sum_{i=1}^t m_i=m-1$) with $m_i \geq x$ for every  {$i\in[1,t]$}, there exists a subset partition $\{S_i\}_{i=1}^t$ of $\Gamma^*$ (i.e., subsets $S_i$ are pairwise disjoint and their union is $ \Gamma^*$) with $|S_i| = m_i$ and $\sum_{s\in S_i}s = 0$ for every {$i\in[1,t]$}. %Moreover there exists a partition $\{m_i\}_{i=1}^t$ of $m-1$ such that $m_1=x$ and there is a partition $\{S_i\}_{i=1}^t$ of $\Gamma-\{0\}$ such that $|S_i| = m_i$ and $\sum_{s\in S_i}s = 0$ for every {$i\in[1,t]$.
\end{dfn}

The following theorem was first conjectured by \citet{KLR} (they also showed the necessity), and later proved by \citet{Zeng}.
\begin{thm}[\citet{Zeng}]\label{Zeng}Let $\Gamma$ be a finite Abelian group of order $m$. $\Gamma$ has $2$-ZSPP if and only if $|I(\Gamma)|\in\{0,3\}$. 
\end{thm}

The above theorem confirms, for the case of $|I(\Gamma)|=3$, the following conjecture stated by Tannenbaum. 
\begin{conj}[\citet{Tannenbaum1}]\label{conjectureT}
Let $\Gamma$ be a finite Abelian group of order $m$ with $|I(\Gamma)|>1$. Let $R =\Gamma \setminus (I(\Gamma)\cup\{0\})$. For every positive integer $t$ and every integer partition $\{m_i\}_{i=1}^t$ of $m-1$, with $m_i \geq 2$ for every {$i\in [1, |R|/2]$}, and {$m_i \geq 3$} for every {$i\in[|R|/2+1, t]$}, there is a subset partition $\{S_i\}_{i=1}^t$ of $\Gamma^*$ such that $|S_i| = m_i$ and $\sum_{s\in S_i}s = 0$ for every {$i\in[ 1, t]$}. 
\end{conj}
It was shown independently by a few authors that Conjecture \ref{conjectureT} is also true for $\Gamma\cong (\zet_2)^n$ for any integer $n$, with $n>1$ (Note that in this case $I(\Gamma)=\Gamma^*$, so $R=\emptyset$).

\begin{thm}[\citet{ref:CaccJia,Egawa,Tannenbaum2}]\label{Sylow}Let $\Gamma\cong (\zet_2)^n$ for some integer $n$, with $n>1$, then $\Gamma$ has $3$-ZSPP. 
\end{thm}

Note that, in general, even the weaker version of Conjecture~\ref{conjectureT} posed by \citet{CicZ} is still open.
\begin{conj}[\citet{CicZ}]\label{conjecture}Let $\Gamma$ be a finite Abelian group with $|I(\Gamma)|>1$. Then $\Gamma$ has $3$-ZSPP. 
\end{conj}

Recently \citet{CicTuz} proved that, if $|\Gamma|$ is large enough and $|I(\Gamma)|>1$, then $\Gamma$ has $4$-ZSPP. 

The main contribution of this paper is to show that Conjecture~\ref{conjecture} holds for every Abelian group $\Gamma$ the cardinality of which is a power of $2$ and $|I(\Gamma)|>1$. This result is a generalization of the result for $\Gamma\cong(\zet_2)^n$ by \citet{ref:CaccJia,Egawa,Tannenbaum2}. The structure of the paper is as follows.
In Section \ref{sec:rt}, we start by establishing the result for $\Gamma\cong \zet_4\times(\zet_2)^{n-2}$, for any integer $n$ with $n \geq 3$. Then, in Section \ref{sec:mr}, we prove the main result that any Abelian group $\Gamma$ of order $2^n$, for a positive integer $n$, such that $|I(\Gamma)|>1$ has $3$-ZSPP. We use some methods developed by \citet{Egawa,Zeng}. Subset partitions of
$\left(\zet_4\times(\zet_2)^{2}\right)^*$, $\left(\zet_4\times(\zet_2)^{3}\right)^*$, $\left(\zet_4\times(\zet_2)^{4}\right)^*$, $\left(\zet_4\times(\zet_2)^{5}\right)^*$,
$\left((\zet_4)^{2}\times\zet_2\right)^*$,
$\left((\zet_4)^{3}\right)^*$,
$\left(\zet_8\times(\zet_2)^{2}\right)^*$, 
and $\left(\zet_4\times(\zet_2)^{2}\right)^*+\left((\zet_2)^{2}\right)^*$
were analyzed by a computer program that we created, and sample zero-sum partitions that certify that the corresponding sequences are realizable are given in the annexes. In Section \ref{sec:sa}, we present some applications of our main result. In Section \ref{sec:final}, we conclude with some final remarks and propose some conjectures. 
%The annexes present zero-sum partitions for the necessary cases in the following order:
%\begin{itemize}
%\item Annex A: $\left(\zet_4\times(\zet_2)^{2}\right)^*$
%\item Annex B: $\left(\zet_4\times(\zet_2)^{3}\right)^*$
%\item Annex C: $\left(\zet_4\times(\zet_2)^{4}\right)^*$
%\item Annex D: $\left(\zet_4\times(\zet_2)^{5}\right)^*$
%\item Annex E: $\left((\zet_4)^{2}\times\zet_2\right)^*$
%\item Annex F: $\left((\zet_4)^{3}\right)^*$
%\item Annex G: $\left(\zet_8\times(\zet_2)^{2}\right)^*$
%\item Annex H: $\left(\zet_4\times(\zet_2)^{2}\right)^*+\left((\zet_2)^{2}\right)^*$
%\end{itemize}

\subsection{Related work and applications}\label{subsec:rw}

$x$-Zero-Sum Partition Property of groups can be applied in magic- and anti-magic-type labelings of graphs\footnote{For standard terms and notations in graph theory, the reader is referred to the textbook by \citet{Diestel}. For an introduction to magic- and anti-magic-type labelings, the reader is referred to the monograph by \citet{BMRS}.}. Generally speaking, such a labeling of a graph $G=(V,E)$ is a mapping from only $V$ or $E$, or their union $V\cup E$, to a set of labels, which most often is a set of integers or elements of a group. Then the weight of a graph element is typically the sum of labels of the adjacent or incident elements of one or both types. When the weight of all elements is required to be equal, then we speak of a magic-type labeling; when the weights should be all different, then we speak of an anti-magic-type labeling. Probably the best known problem in this area is the {\em anti-magic conjecture} by \citet{HR}, which claims that the edges of every graph except $K_2$ can be labeled bijectively with integers $1, 2, \dots, |E|$ so that the weight of every vertex is unique. This conjecture is still open. 

It is easy to check, by the Pigeonhole Principle, that in any simple graph $G$ there exist two vertices of the same degree. The situation changes if we consider an edge labeling $f:E(G)\rightarrow \{1, \ldots, k\}$ (Note that the labels do not have to be distinct) and calculate the so-called \textit{weighted degree} of every vertex $v$ as the sum of labels of all the edges incident to $v$. The labeling $f$ is called \textit{irregular} if the weighted degrees of all the vertices are unique (so it is an anti-magic-type labeling). The smallest value of $k$ that allows some irregular labeling is called the \textit{irregularity strength of $G$} and denoted by $s(G)$. The problem of finding $s(G)$ was introduced by \citet{ref:ChaJacLehOelRuiSab1} and investigated by numerous authors  \citep{ref:AigTri2,ref:AmaTog,ref:KalKarPfe1}. 

\citet{ref:AnhCic1} introduced the \textit{group irregular labeling}. They label the edges with (not necessarily distinct) elements of a group $\Gamma$, and then the weighted degree of a vertex  $v$ is the sum (taken in $\Gamma$) of labels of the edges incident to $v$. A labeling is \textit{$\Gamma$-irregular}, if the resulting weighted degrees are pairwise distinct. The smallest $k$, such that for every Abelian group $\Gamma$ of order $k$ there exists a $\Gamma$-irregular labeling of $G$, is called  the \textit{group irregularity strength} of $G$ and denoted $s_g(G)$. Note that $s(G)\leq s_g(G)$ for every graph $G$.

\citet{ref:AnhCic} used $2$-ZSPP of groups of odd order for bounding the group irregularity strength of disconnected graphs without a star as a connected component. Roughly speaking, the authors divide every connected component into $2$-subsets and $3$-subsets of vertices and partition the set of non-zero elements of the corresponding group into the same number of zero-sum $2$-subsets and $3$-subsets, and later use the method of augmented paths to do the labeling. They obtain that $s_g(G)\leq 2\lfloor|V(G)|/2\rfloor+1$.

An analogous \textit{group irregular labeling for directed graphs}, where the weight of every vertex is the sum of the labels flowing into the vertex minus the sum of labels flowing out of the vertex, was considered by \citet{ref:AigTri}, \citet{ref:CaccJia}, \citet{CicTuz}, \citet{Egawa}, \citet{Fukuchi} and \citet{Tuza}. It turns out that there exists a $\Gamma$-irregular labeling of a directed graph $\overrightarrow{G}$ with weakly connected components $\{\overrightarrow{G_i}\}_{i=1}^t$ if and only if there exist in $\Gamma$ pairwise disjoint subsets $\{S_i\}_{i=1}^t$ such that $|S_i|=|V(\overrightarrow{G_i})|$  and $\sum_{s\in S_i}s=0$ for every  {$i\in[1, t]$} {(see Theorem 2.1 in \citet{CicTuz})}.

\citet{KLR} and \citet{Zeng} used results on zero-sum partitions of Abelian groups for another type of anti-magic labeling. They label the edges bijectively with the non-zero elements of a group $\Gamma$. The weight of every vertex is calculated as the sum (taken in $\Gamma$) of the labels of incident edges. All the weights should be different. They showed that every $2$-tree\footnote{A $2$-tree $T$ is a rooted tree, where every vertex which is not a leaf has at least two children.} $T$ of order $n$ admits a $\Gamma$-anti-magic labeling (in the above sense) for any group $\Gamma$ of order $n$ such that $|I(\Gamma)|\in\{0,3\}$. 

\citet{Fro1} defined the notion of \textit{group distance magic labeling}. In this case, the vertices of the graph are labeled through a bijection with the elements of an Abelian group $\Gamma$. The weight of every vertex is computed as the sum (in $\Gamma$) of the labels assigned to its neighbors. If all the weights are the same, then it is a \textit{$\Gamma$-distance magic labeling}.

Note that a group $\Gamma$ of order $m$ has a constant sum partition (of all of its elements) if and only if there exists a certain complete multipartite $\Gamma$-distance magic graph. Indeed, suppose we have a constant sum partition of $\Gamma$: $\{S_i\}_{i=1}^t$ of $\Gamma$ with $|S_i| = m_i$ and $\sum_{s\in S_i}s = \nu$ for every {$i\in[1,t]$} and some $\nu \in \Gamma$. Let $G$ be a complete $t$-partite graph with the color classes $\{A_i\}_{i=1}^t$, where $|A_i|=m_i$ for every  {$i\in[1 ,t]$}. Let us label the vertices of $A_i$ with distinct elements of $S_i$ for every {$i\in[1 ,t]$}, and compute the weight of every vertex as the sum (in $\Gamma$) of the labels of its neighbors. It is easy to see that all the weights are equal and thus we have a $\Gamma$-distance magic labeling. On the other hand, suppose $G=K_{m_1,\ldots,m_t}$ is a complete $t$-partite graph of order $m$ with the color classes $\{A_i\}_{i=1}^t$ that is $\Gamma$-distance magic  {with a labeling $\ell$}. So $\sum_{i=1,i\neq j}^t\sum_{x\in A_i}\ell(x) = \mu$ for every {$j\in[1 , t]$}, which implies that $\sum_{x\in A_j}\ell(x) = \nu$ for every  {$j\in[1 , t]$}, and some $\nu\in\Gamma$. 
 
\section{\texorpdfstring{Realizable triples in $\left(\zet_4\times(\zet_2)^{n-2}\right)^*$}{Realizable triples}}\label{sec:rt}

Throughout this section, let $n$ be an integer, $n \geq 3$, and let  $Y\cong\zet_4\times(\zet_2)^{n-2}$. For a subset $S$ of $Y$, let $\langle S\rangle$ denote the subgroup of $Y$ generated by $S$. To simplify the notation, we will write just $\langle s_1, \dots, s_k \rangle$ instead of $\langle \{ s_1, \dots, s_k \} \rangle$. For any subsets $S$ and $T$ of $Y$, let $S+T=\{u+v\colon u\in S, v\in T\}$. Recall that the exponent of a group $\Gamma$, denoted by $e(\Gamma)$, is defined as the least common multiple of the orders of all elements of the group. Thus every non-zero subgroup of $Y$ has its exponent equal to either $2$ or $4$. The order of $Y$ is $m=2^n$, and we will consider integer partitions $\{m_i\}_{i=1}^t$ of $m-1$ for some positive integer $t$, with $m_i \geq 3$ for every {$i\in[1,t]$}.

Note that, since $m_i\geq 3$ for every $i\in[1,t]$, we can modify the sequence $\{m_i\}_{i=1}^t$ by subdividing every term larger than $5$ into a combination of terms $3$, $4$, and $5$. It is easy to check that, if the non-zero elements of a group $\Gamma$ can be partitioned into zero-sum subsets of cardinalities corresponding to the elements of the new sequence, then the same holds for the old sequence. So we can focus only on sequences with $m_i \in \{3, 4, 5\}$ for every $i\in[1,t]$.

Let $Z$ be a subset of $Y$. Let $\mathcal{K}$ be a partition of $Z$ into zero-sum subsets with cardinalities in $\{3,4,5\}$, and let $a = |\{S \in\mathcal{K}\colon |S|=3\}|$, $b = |\{S \in\mathcal{K}\colon |S|=4\}|$ and $c = |\{S \in\mathcal{K}\colon |S|=5\}|$. In this situation, we say that $\mathcal{K}$ \textit{realizes} the \textit{triple} $(a, b,c)$ in $Z$. If there exists a partition realizing $(a, b,c)$ in $Z$, we say that $(a, b,c)$ is \textit{realizable} in $Z$. Note that a triple $(a,b,c)$ can be seen as a compact representation of a sequence $\{m_i\}_{i=1}^{t}$ in which $a$ elements are equal to $3$, $b$ elements are equal to $4$, and $c$ elements are equal to $5$.

Let us start by recalling three useful lemmas.

\begin{lem}[\citet{Egawa}]\label{l2.2} Let $r$ be a non-negative integer. Let $a,$ $b,$ $c$ be non-negative integers such that $3a+4b+5c\geq 45r+12$ and
\begin{align}\label{2.9}
[(b-1)/9]\leq(a/3)+c.\end{align}
Then there exist non-negative integers $\{x_i\}_{i=1}^r$, $\{y_i\}_{i=1}^r$, $\{z_i\}_{i=1}^r$, such that $3x_{i}+4y_{i}+5z_{i}=45$ for every  {$i\in [1 , r]$},  $\sum_{i=1}^r x_{i}\leq a$, $\sum_{i=1}^r y_{i}\leq b$, and $\sum_{i=1}^r z_{i}\leq c$.\end{lem}

\begin{lem}[\citet{Egawa}]\label{l3.11} Let $X\cong (\zet_{2})^n$, for some integer $n$, with $n\ge 2$. If $n$ is odd, let $W$ denote a subgroup of $X$ of order $2^3$; if $n$ is even, let
$W=\{0\}$. Then the triple $((|X|-|W|)/3,0,0)$ is realizable in $X\backslash W$.\end{lem} 

\begin{lem}[\citet{CicZ}]\label{cichacz} Let $\Gamma$ be an Abelian group such that $|I(\Gamma)|>1$. Let $k$ be a positive even integer, with $k\geq 4$, such that $k$ divides $|\Gamma|$. Then there exists a partition of $\Gamma$ into sets $\{A_i\}_{i=1}^{|\Gamma|/k}$ such that $|A_i|= k$ and $\sum_{a\in A_i} a=0$ for every {$i\in[1 , |\Gamma|/k]$}.\end{lem}

We will also need the following lemmas:

%\comS{Nowa, bardziej ogólna wersja:}
\begin{lem}\label{l3.12} Let $\Gamma\cong H\times \zet_{2^{n_1}}$ be an Abelian group such that $|H|=4n$, $|I(H)|>1$, $n_1\in\{1,2\}$, and $n$ is a positive integer. Let $a$, $b$, $c$ be non-negative integers with $3a+4b+5c=|\Gamma|-1$ and $b\geq (2^{n_1}-1)n$. If there exists a partition $\mathcal{K}$ realizing $(a, b-(2^{n_1}-1)n, c)$ in $H^{*}$, then $(a, b, c)$ is realizable in $\Gamma^{*}$. \end{lem}

\begin{proof} 
By Lemma~\ref{cichacz}, there exists a partition $\mathcal{A}$ of $H$ into zero-sum $4$-subsets. For every $K \in \mathcal{K}$, let $K_0 = \{(g,0) : g\in K\}$, and let $\mathcal{K}_0 = \{K_0 \colon K \in \mathcal{K}\}$. For every $A \in \mathcal{A}$ and every  {$i\in[1, 2^{n_1}-1]$}, let $A_i = \{(g,i) \colon g \in A\}$, and let $\mathcal{A}_i = \{A_i \colon A \in \mathcal{A}\}$. Observe that $\mathcal{K}_0$ realizes the triple $(a, b-(2^{n_1}-1)n, c)$ in $H^*\times\{0\}$, whereas $\bigcup_{i=1}^{2^{n_1}-1}\mathcal{A}_{i}$ realizes the triple $(0, (2^{n_1}-1)n,0)$ in $H\times (\zet_{2^{n_1}})^*$. Thus $\mathcal{K}_0 \cup \bigcup_{i=1}^{2^{n_1}-1}\mathcal{A}_{i}$ realizes $(a, b, c)$ in $\Gamma^*$.
\end{proof}

%\comS{A to stara wersja:}
%\begin{lem}\label{ls} Let \edtS{$n\geq 4$}, and let $U$ be a subgroup of $Y$ of order $2^{n-1}$. Let $a$, $b$, $c$ be non-negative integers with $3a+4b+5c=2^{n}-1$ and $b\geq 2^{n-3}$. If there exists a family $K$ realizing $(a, b-2^{n-3}, c)$ in $(U\backslash \{0\})^*$, then $(a, b, c)$ is realizable in $Y^*$. \end{lem}

%\begin{proof} Observe that $U\cong\zet_4\times(\zet_2)^{n-3}$. Note that there exists   a subgroup $W\cong (\zet_2)^2$ of $U$. Since $e(Y)\in\{2,4\}$, the family $L$ consisting of those cosets of $W$ which are disjoint from $U$ realizes $(0,2^{n-3},0)$, and hence $K\cup L$ realizes $(a, b, c)$.\end{proof} 

\begin{lem}\label{l3.8} Let $W\cong(\zet_2)^3$ be a subgroup of $Y$, let $P$ be a zero-sum $3$-subset of $Y$, suppose that $W\cap\langle P\rangle=\{0\}$. Then $(8,0,0)$ is realizable in $W+P$.
\end{lem}
\begin{proof}
Choose $v_0, v_{1}, v_{2}$ to be distinct elements of $Y$ that generate $W$, so $W=\langle v_0,v_{1},v_{2}\rangle$. Let $p_0,p_{1},p_{2}$ be the elements of $P$. In the following, the  subscripts should be interpreted modulo 3. For  $k\in\{0,1,2\}$, let
$P_{k}=\{p_{k}, v_{k}+p_{k+1}, v_{k}+p_{k+2}\}$,
$S_{k}=\{v_{k}+p_{k}, v_{k+1}+v_{k+2}+p_{k+2}, v_{k}+v_{k+1}+v_{k+2}+p_{k+1}\}$.
For every $l\in\{0,1\}$, let
$T_{l}=\{v_{i}+v_{i+1}+p_{i+2+l}\colon0\leq i\leq 2\}$.
Then $\{P_{k} : 0\leq k \leq 2\} \cup \{S_{k} : 0\leq k \leq 2\} \cup \{T_{l} : 0 \leq l \leq 1\}$ realizes $(8, 0,0)$ in $W+P$.
\end{proof}

\begin{lem}\label{l3.9} Let $W\cong(\zet_2)^3$ be a subgroup of $Y$, let $R$ be a zero-sum $5$-subset of $Y$, and suppose that $W\cap\langle R\rangle=\{0\}$. Then $(0,0,8)$ is realizable in $W+R$.\end{lem}
\begin{proof}

Let $ W=\langle v_{0}, v_{1}, v_{2}\rangle$, $R=\{p_{0}, p_{1}, p_{2}, q, r\}$, $P_{k}$ and $S_{k}$ for every $0\leq k \leq 2$, and the sets $T_{l}$ for every  $0\leq l\leq 1$ be similar to the definitions in the proof of Lemma~\ref{l3.8}. Here also the subscripts should be interpreted modulo 3.

 {Additionally, let
$R_{k}=P_{k}\cup\{v_{k}+q, v_{k}+r\}$ and
$U_{k}=S_{k}\cup\{v_{k+1}+v_{k+2}+q,
v_{k+1}+v_{k+2}+r\}$ for every  $0\leq k \leq 2$, and
$V_{0}=T_{0}\cup\{q, r\}$, $V_{1}=T_{1}\cup\{v_{0}+v_{1}+v_{2}+q, v_{0}+v_{1}+v_{2}+r\}$.
Then $\{R_{k} : 0\leq k \leq 2\} \cup \{U_{k} : 0\leq k \leq 2\} \cup \{V_{i} : 0 \leq i \leq 1\}$ realizes $(0,0,8)$ in $W+R$.}
\end{proof}

%\comK{Coś tu jest nie tak. Mamy dwa Proof, jeden za drugim.}
%\begin{lem}\label{l4.2} Let $a,$ $b,$ $c$ be nonnegative integers with$3a+4b+5c=63$. Then $(a, b, c)$ is realizable in $\left(\zet_4\times(\zet_2)^4\right)^{*}$.\end{lem}
%\begin{proof}
%The case of $\zet_4\times(\zet_2)^2$ was analyzed by a computerprogram we created, sample zero-sum partition of this group can be found in the annexes.  Thus, by Lemma~\ref{l}, we may assume $b\leq 7$, and hence we have $a>5$ or $c>3$. If $a>5$, let $a_{1}=5$ and $c_{1}=0$; if $a\leq 5$ (so $c>3$), let $a_{1}=0$ and $c_{1}=3$. Let $U,$ $V$ be subgroups of $Y$ such that $U\cong \zet_4\times(\zet_2)^2$, $V\cong \zet_2$ and $U\cap V=\{0\}$.  $(a_{1},0, c_{1})$ is realizable in $(U\backslash \{0\}$ by Lemma 3.13, and $(a-a_{1}, b, c-c_{1})$ is realizable in $\left(U+(V\backslash \{0\}\right)^*$ by Lemma 3.16. Hence $(a, b, c)$ is realizable in $X\backslash \{0\}$.
%\end{proof}

We consider first the cases of $Y$ such that $n\leq 7$
\begin{lem}\label{z224m} Let $n$ be an integer with $3\leq n\leq 7$. Let $a,b,c$ be non-negative integers such that $3a+4b+5c=2^n-1$. Then $(a,b,c)$ is realizable in $Y^*$. \end{lem}
\begin{proof}
By Theorem~\ref{Zeng}, $\zet_4\times\zet_2$ has $2$-ZSPP. For $n\geq 4$, the triples $(a, b, c)$ with $b<2^{n-3}$ were analyzed by a computer program we created, sample realizations can be found in the annexes. Therefore, for $H =\zet_4\times(\zet_2)^{n-3}$ and  $n_1=1$, by Lemma~\ref{l3.12}, we obtain that all triples $(a,b,c)$ are realizable in $Y^*$.
\end{proof}

\begin{thm}\label{z224} Let $n$ be an integer with $n\geq 3$. Let $a,b,c$ be non-negative integers such that $3a+4b+5c=2^n-1$. Then $(a,b,c)$ is realizable in $Y^*$. \end{thm}
%Let $n,$ $a,$ $b,$ $c$ be as in Theorem
%2 and, as in the preceding section, let $Y\cong\zet_4\times(\zet_2)^{n-2}$.
%Let $n\geq 5$. 

\begin{proof}
By Lemma~\ref{z224m}, we can assume that $n\geq 8$.
Taking  $H =\zet_4\times(\zet_2)^{n-3}$ and  $n_1=1$, by Lemma \ref{l3.12}, we may assume
\begin{align} \label{5.1}
b<2^{n-3},
\end{align}
and hence
\begin{align} \label{5.2}
3a+5c>2^{n-1}.
\end{align}

Let $V\cong\zet_4\times (\zet_2)^2$ and $U\cong (\zet_2)^{n-4}$  be subgroups of $Y$ such that $U\cap V=\{0\}$. If $n$ is odd, let $W\cong(\zet_2)^3$
be a subgroup  of $U$; if $n$ is even, let $W=\{0\}$. Since $n\geq 8$, we have
\begin{align} \label{5.3}
|W|\leq 2^{n-6}.
\end{align}

We will show that there is $a=a_{1}+a_{2}+a_{3},$ $b=b_{1}+b_{2}+b_{3}$ and $c=c_{1}+c_{2}+c_{3}$,
so that $(a_{1}, b_{1}, c_{1})$, $(a_{2}, b_{2}, c_{2})$, and $(a_{3}, b_{3}, c_{3})$ are realizable in $W+V^{*}$ ,
$(U\backslash W)+V^{*}$, and $U^{*}$, respectively. 

Let us first consider $W+V^{*}$. By~(\ref{5.2}), we have $3a>2^{n-2}$ or $5c>2^{n-2}$. Assume
first that $3a>2^{n-2}$. In this case, let $(a_{1}, b_{1}, c_{1})=(5|W|, 0,0)$. By~(\ref{5.3}), we have $a_{1}<a$.
By Lemma~\ref{z224m}, we can partition $V^{*}$ into five zero-sum $3$-subsets $P_{0},$ $\cdots,$ $P_{4}$. Then, for every $i\in[0,4]$, there exists a subset partition $\mathcal{K}_{i}$ of $W+P_{i}$ that realizes $(|W|, 0,0)$ (for $n$ odd we apply Lemma~\ref{l3.8}). Hence
$\mathcal{K}=\bigcup_{i=0}^{4}\mathcal{K}_{i}$ realizes $(a_{1}, b_{1}, c_{1})$ in $W+V^{*}$. Assume now that $3a\leq 2^{n-2}$, so
$5c>2^{n-2}$. In this case, we take $(a_{1}, b_{1}, c_{1})=(0,0,3|W|)$. By (\ref{5.3}), we have $c_{1}<c$. By Lemma~\ref{z224m}, we can partition $V^{*}$ into three zero-sum $5$-subsets $R_{0},$ $R_{1},$ $R_{2}$. Then we can see from Lemma~\ref{l3.9} that, for every $i\in\{0,1,2\}$, there exists a subset partition $\mathcal{K}_{i}$ of $W+R_{i}$ realizing $(0, 0, |W|)$. Thus, the partition $\mathcal{K}=\bigcup_{i=0}^{2}\mathcal{K}_{i}$ realizes
$(a_{1}, b_{1}, c_{1})$ in $W+V^{*}$.

Let us now consider $(U\backslash W)+V^*$ and $U^*$. Let $r=(|U|-|W|)/3$. Then
\begin{align} \label{5.4}
3(a-a_{1})+4(b-b_{1})+5(c-c_{1})=45r+|U^*|.\end{align}
Like \citet{Egawa}, we also have
$$[((b-b_{1})-1)/9]<(a-a_{1})/3+(c-c_{1}).$$
Since (\ref{5.4}) implies $3(a-a_{1})+4(b-b_{1})+5(c-c_{1})=45r+(2^{n-4}-1)\geq 45r+15$, it now
follows from Lemma~\ref{l2.2} that there exist non-negative integers $\{x_{i},y_{i},z_{i}\}_{i=1}^{r}$ such that
\begin{align}\label{5.5}
3x_{i}+4y_{i}+5z_{i}=45\end{align}
for every  {$i\in[1,r]$},  $\sum_{i=1}^{r} x_{i}\leq a-a_{1},$ $\sum_{i=1}^{r} y_{i}\leq b-b_{1}$ and $\sum_{i=1}^{r} z_{i}\leq c-c_{1}$. Let $a_{2}=\sum_{i=1}^{r} x_{i},$ $b_{2}=\sum_{i=1}^{r} y_{i},$ $c_{2}=\sum_{i=1}^{r} z_{i}$,
$a_{3}=a-a_{1}-a_{2},$ $b_{3}=b-b_{1}-b_{2}$, and $c_{3}=c-c_{1}-c_{2}$. By Theorem~\ref{Sylow}, it follows from (\ref{5.4}) and (\ref{5.5}) that there exists a subset partition $\mathcal{L}$ of $U^*$ realizing $(a_{3}, b_{3}, c_{3})$.
By Lemma \ref{l3.11}, we can partition $U\backslash W$ into $r$ zero-sum $3$-subsets $\{S_{i}\}_{i=0}^{r-1}$. Since $S_i\cong ((\zet_2)^2)^*$ for every  {$i\in[0,r-1]$}, the solutions for $S_{i}+V^*$,  available in the annexes certify that, for every $i$, there exists a partition $\mathcal{N}_i$ realizing $(x_{i}, y_{i}, z_{i})$ in  $S_{i}+V^*$. %, and we see from Lemma~\ref{Karol} that for every $1\leq i\leq r$, there exists a sy $N_{i}$ of subsets of $S_{i}+(V\backslash \{0\})$ realizing $(x_{i}, y_{i}, z_{i})$. 
Then $\bigcup_{i=0}^{r-1}\mathcal{N}_{i}$ realizes $(a_{2}, b_{2}, c_{2})$ in $(U\backslash W)+V^*$. Hence, the partition $\mathcal{K}\cup\mathcal{ L}\cup(\bigcup_{i=0}^{r-1}\mathcal{N}_{i})$ realizes $(a, b, c)$ in $Y^*$.
\end{proof}

\section{Main result}\label{sec:mr}
 In this section we will use the idea of good subsets from \citet{Zeng}. We call a 6-subset $C$ of an Abelian group $\Gamma$ \textit{good} if $C = \{c, d,-c -d,-c,-d, c + d\}$ for some $c$ and $d$ in $\Gamma$. This idea is strongly connected with Skolem partitions of groups \citep{Tannenbaum1}. Note that a good 6-subset is partitionable into three zero-sum 2-subsets
as well as into two zero-sum 3-subsets. We will need the following lemma.

\begin{lem}[\citet{CicZ,Zeng}]\label{bijection}
Let $\Gamma$ be a finite Abelian group such that $|I(\Gamma)|\neq 1$. Let Bij$(\Gamma)$ denote the set of all bijections from $\Gamma$ to itself. Then there exist $\phi,\varphi\in$Bij$(\Gamma)$ (not necessarily distinct) such that $g+\phi(g)+\varphi(g)=0$  for every $g\in\Gamma$. %In particular, we mayassumethat #(0)='(0)=0.
\end{lem}

%For a group $H$, let $H^{*}=H\setminus \{0\}$.
We start with some small cases.

\begin{lem}\label{male2} Let $\Gamma$ be any of the following groups: $(\zet_4)^2\times \zet_2$, $(\zet_4)^3$, $\zet_8\times (\zet_2)^2$. Let $a,b,c$ be non-negative integers such that $3a+4b+5c=|\Gamma|-1$. Then $(a,b,c)$ is realizable in $\Gamma^*$. \end{lem}

\begin{proof}
$\zet_4\times\zet_4$ has $2$-ZSPP by Theorem~\ref{Zeng}. For $\Gamma\cong(\zet_4)^2\times \zet_{2^{n_1}}$, with $n_1\in\{1,2\}$, the triples $(a, b, c)$ with $b<(2^{n_1}-1)4$ were analyzed by a computer program we created, sample realizations can be found in the annexes. Therefore, for $H =(\zet_4)^2$ and $n_1\in\{1,2\}$, applying Lemma~\ref{l3.12}, we obtain that all triples $(a,b,c)$ are realizable in $\Gamma^*$.

The case for $\zet_8\times\zet_2$ is done by Theorem~\ref{Zeng}. For $\Gamma\cong\zet_8\times (\zet_2)^2$, the triples $(a, b, c)$ with $b<4$ were analyzed by a computer program we created, sample realizations can be found in the annexes. Therefore, for $H =\zet_8\times\zet_2$ and  $n_1=1$, applying Lemma~\ref{l3.12}, we obtain that all triples $(a,b,c)$ are realizable in $\Gamma^*$.
\end{proof}

Recall that the quotient group of $\Gamma$ by $H$, for a subgroup $H$ of $\Gamma$, is denoted by $\Gamma/H$.
Now we state our main result.
\begin{thm}Let $\Gamma$ be such that $|I(\Gamma)|>1$ and $|\Gamma|=2^n$ for some integer $n$, with $n>1$. Then $\Gamma$ has $3$-Zero-Sum Partition Property.\label{mainSK}
\end{thm}

The proof is  by contradiction (using the method
of smallest counterexample). Suppose the theorem is false and let $\Gamma$ be the smallest group the order of which is a power of $2$, with $|I(\Gamma)|>1$, without $3$-ZSPP.

\textit{Case 1.} Suppose $\Gamma \cong \zet_{2^{\alpha}}\times H$, where $|H|$ is even and $\alpha\geq 3$. Note that for $H$ such that $|I(H)|=1$, we have $\Gamma \cong \zet_{2^{\alpha}}\times \zet_{2^{\beta}}$ for $\beta>0$ and $\Gamma$ has $3$-ZSPP by Theorem~\ref{Zeng}. Thus we can assume that $|I(H)|\geq3$. There exists a subgroup $L \cong   \zet_{2^{\alpha-3}}\times H$ of $\Gamma$ such that $|I(L)|>1$ and
$\Gamma/L \cong \zet_8$. Since $\zet_8=\{0,4\}\cup\{1,2,5,-1,-2,-5\}$, we can choose a set of coset representatives for the subgroup $L$ in $\Gamma$, say  $A$, such
that $$A = \{0, e\} \cup \{d, f,-d - f,-d,-f, d + f\},$$
where $2e \in L$. Since  $|I(L)|>1$, by Lemma~\ref{bijection}, there exist $\phi,\varphi\in$Bij$(L)$ such that $g+\phi(g)+\varphi(g)=0$ for every $g\in L$. Thus
$$\Gamma^{*} = L^{*} \cup(e+L) \cup \left(\bigcup \limits_{g\in L}\{d+g, f+\phi(g),-d-f+\varphi(g),-d-g,-f-\phi(g), d+f-\varphi(g)\}\right),$$
where the latter 6-subsets are good. Note that $B=L\cup(e+L)$ is a subgroup of $\Gamma$ such that $|I(B)|\geq 3$. So, $B$ has $3$-ZSPP.
Observe that $|B|=2^{\gamma}$ for some positive integer $\gamma\geq 2$. Let
 $$W=\left(\bigcup\limits_{g\in L}\{d+g, f+\phi(g),-d-f+\varphi(g),-d-g,-f-\phi(g), d+f-\varphi(g)\}\right).$$ Note that $|W|=6|L|\equiv0\pmod 4$. We will prove that any triple $(a,b,c)$ such that $3a+4b+5c=|\Gamma|-1$ is realizable in $\Gamma^*$ and obtain a contradiction. 
 Assume that $r_1,\dots,r_{a}=3$, $r_{a+1},\ldots,r_{a+c}=5$, and $r_{a+c+1},\dots,r_{a+b+c}=4$. Let $l$ be such that $\sum_{i=1}^{l-1} r_i\leq|B^*|$ and $\sum_{i=1}^{l} r_i> |B^*|$.
Let $r_{l}'=|B^{*}|- \sum_{i=1}^{l-1} r_i$ and $r_l''=r_l-r_l'$. 
If ($r_l'=0$ or $r_l'\geq 3$) and ($r_l''=0$ or  $r_l''\geq 2$), then the sequence $r_1,\ldots,r_{l-1},r_{l}'$ is realized by a zero-sum partition $A_1,\ldots,A_{l-1},A_{l}'$ of $B^*$. Moreover, the sequence $r_{l}'', r_{l+1}, r_{l+2}, \ldots, r_{a+b+c}$ is realized by a zero-sum partition $A_{l}'',A_{l+1},A_{l+2},\ldots,A_{a+b+c}$ of $W$, since $W$ is the union of good $6$-subsets and $|W|=\sum_{i=l+1}^{a+b+c}r_i+r_{l}''$, so we are done. Hence we have to settle only the cases where $r''_l=1$ or $r'_l=1$ or $r'_l=2$.

 \textit{Case 1.a.}  $r_l''=1$.
Then $l$ is even and $r_{l+1}$ is odd, since $|B^{*}|$ is odd and $|W|$ is even
 (thus $|W|-r''_l$ requires at least one odd term among $r_{l+1},\dots,r_{a+b+c}$,
 and all odd terms are put at the beginning of the sequence).
{Suppose first $r_l=3$, then $3l-1=|B^{*}|=2^{\gamma}-1$ and so
$3l=2^{\gamma}$, a contradiction. Thus, suppose} that $r_l=5$. If $b>0$, then the triple $(a,1,l-a-1)$ is realizable in $B^*$, and $(0,b-1,c+a-l+1)$ is realizable in $W$. Assume now that $b=0$. If $|B^{*}|>7$, then, for $a\geq 2$, the triple $(a-2,0,l-a+1)$ is realizable in $B^*$ and $(2,0,c+a-l-1)$ is realizable in $W$. For $a\leq1$, there is  $r_{l-1}=5$, and $(a+3,0,l-a-2)$ is realizable in $B^*$, and $(1,2,c+a-l-2)$ is realizable in $W$.
Observe that in this case, by the construction of $W$, the two zero-sum $4$-subset in $W$ can be split into four zero-sum $2$-subsets, and we obtain a realization of $(a,0,c)$ in $\Gamma^*$. If $|B^{*}|=7$, then $a=1$. Now  $|I(B)|=7$, so $B\cong(\zet_2)^3$, which implies that $\Gamma \cong (\zet_2)^2\times\zet_8$, and this case was considered in Lemma~\ref{male2}. 

 \textit{Case 1.b.}  $r_l'=2$.
If $r_l=3$, then we are done as in Case 1.a.. Since $|W|\equiv 0\pmod 4$, we can assume $r_l=5$. For $a\geq1$, we obtain that {$(a-1,0,l-a)$ is realizable in $B^*$, and $(1,b,c+a-l)$ is realizable in $W$. For $b\geq 1$, since $|B^{*}|\geq 7$},  we have a realization of $(1,1,l-2)$ in $B^*$ and a realization of $(1,b,c-l)$ in $W$. By the construction of $W$, one zero-sum $4$-subset in $W$ can be split into two zero-sum $2$-subsets, and we obtain a realization of $(a,b,c)$ in $\Gamma^*$.
The only missing case is $a=b=0$. If $|B^{*}|>7$, then actually $|B^{*}|\geq 15$ and $r_{l-1}=r_{l-2}=5$, thus $(4,0,l-3)$ is realizable in $B^*$. Since $W$ is the union of good 6-subsets, we have a partition of $W$ into four zero-sum $2$-subsets and $(c-l-1)$ zero-sum $5$-subsets.

Assume now that $|B^{*}|=7= |I(B)|$. But then, as in Case 1.a., we have $\Gamma \cong (\zet_2)^2\times\zet_8$, and this case was considered in Lemma~\ref{male2}.

\textit{Case 1.c.} $r_l'=1$.
{Then $l> 1$}, and $r_l''\geq 2$ is even. 
Assume first that $r_l=3$. Since $|W|\equiv 0\pmod 3$, we have $b+c>0$.
If now $b>0$, then there exists a realization of $(l-2,1,0)$  in $B^*$ and a realization of $(a-l+2,b-1,c)$ in $W$.
 If $|B^{*}|>7$, then $|B^{*}|\geq 15$ and $r_{l-1}=r_{l-2}=r_{l-3}=3$. If $b=0$ and $c=1$, then $|\Gamma|\equiv 0\pmod 3$, a contradiction. Therefore we can assume that $b=0$ and $c\geq 2$, but then  $(l-4,0,2)$ is realizable in $B^*$, and {$(a-l+4,0,c-2)$} is realizable in $W$. For the case $|B^{*}|=7= |I(B)|$ (i.e. $\Gamma \cong (\zet_2)^2\times\zet_8$), we apply Lemma~\ref{male2}.
Suppose now that $r_l=5$. If $r_{l-1}=5$, then we set $r_{l-1}'=3$ and
 $r_{l-1}''=2$, and re-define $r_l':=3$ and $r_{l}'':=r_l-3$.
The sequence $r_1,\ldots,r_{l-2},r_{l-1}',r_l'$ is realized by a zero-sum partition $A_1,\ldots,A_{l-2}, A_{l-1}', A_l'$ in $B^{*}$, and the sequence $r_{l-1}'', r_l'', r_{l+1}, \ldots, r_{a+b+c}$ is realized by a zero-sum partition $A_{l-1}'', A_l'', A_{l+1}, \ldots, A_{a+b+c}$ of $W$. Thus $r_{l-1}=3$. If $b>0$, then we have a realization of $(a-1,1,0)$ in $B^*$ and a realization of $(1,b-1,c)$  in $W$. If $|B^{*}|>7$, then $|B^{*}|\geq 15$ and $r_{l-1}=r_{l-2}=r_{l-3}=3$. Since $b=0$ and $c\geq 2$, we obtain that $(a-3,0,2)$ is realizable in $B^*$, and $(3,0,c-2)$ is realizable in $W$. The case $|B^{*}|=7= |I(B)|$ implies that $\Gamma \cong (\zet_2)^2\times\zet_8$ and we apply Lemma~\ref{male2}.

\textit{Case 2.} Suppose $\Gamma \cong \zet_{4}\times \zet_4\times H$. By Theorem~\ref{Zeng} and Lemma~\ref{male2}, we can assume that $|H|\geq 8$ or $H\cong\zet_2\times\zet_2$. Moreover, by Case 1., we can assume that $e(\Gamma)=4$, which implies that we can assume that $|I(H)|>1$.
 Then there is a subgroup $L \cong H$ of $\Gamma$ such that
$\Gamma/L \cong \zet_4\times\zet_4$. Because $\zet_4\times\zet_4=\{(0,0), (0,2), (2,0), (2,2)\}\cup\{(0,1), (1,2), (3,1), (0,3), (3,2), (1,3)\}\cup
\{(1,0), (1,1), (2,3), (3,0),$\\$(3,3), (2,1)\}$, we can choose a set of coset representatives $A$ for the subgroup $L$ in $\Gamma$, such that 
%$$B = \{0, e_1,e_2,e_1+e_2\} \cup \{b_1, c_1,-b_1 - c_1,-b_1,-c_1, b_1 + c_1\}\cup \{b_2, c_2,-b_2 - c_2,-b_2,-c_2, b_2 + c_2\},$$
$$A = \{0, e_1,e_2,e_1+e_2\}\cup\left( \bigcup\limits_{i=1,2}\{d_i, f_i,-d_i - f_i,-d_i,-f_i, d_i + f_i\}\right),$$
where $2e_1,2e_2 \in L$. Since  $|I(L)|>1$, by Lemma~\ref{bijection}, there exist
 $\phi,\varphi\in$Bij$(L)$ such that $g+\phi(g)+\varphi(g)=0$  for every $g\in L$. Thus
$$\begin{array}{ll}
\Gamma^{*} =& L^{*}\cup(e_1+L)\cup(e_2+L)\cup(e_1+e_2+L)\cup\\
&\left(\bigcup\limits_{\substack{i=1,2 \\ g\in L}}\left\{d_i+g, f_i+\phi(g),-d_i-f_i+\varphi(g),-d_i-g,-f_i-\phi(g), d_i+f_i-\varphi(g)\right\}\right),\end{array}$$
where the latter 6-subsets are good. Note that $L\cup(e_1+L)\cup(e_2+L)\cup(e_1+e_2+L)$ is a subgroup of $\Gamma$. So, $L\cup(e_1+L)\cup(e_2+L)\cup(e_1+e_2+L)$ has $3$-ZSPP. Note that for $B=L\cup(e_1+L)\cup(e_2+L)\cup(e_1+e_2+L)$ there is $B\cong \zet_{2}\times \zet_2\times H$ and  $|I(B)|\geq 15$. Now, a similar reasoning to that applied in Case 1. leads to obtaining that $\Gamma$ has $3$-ZSPP, a contradiction.\\

\textit{Case 3.} $\Gamma \cong \zet_{4}\times (\zet_2)^{n-2}$. In this case we are done by~Theorem~\ref{z224}.  \\

Since the case $\Gamma\cong(\zet_2)^n$ has been solved by Theorem~\ref{Sylow} this completes the proof
of the theorem.~\qed

\section{Some applications}\label{sec:sa}
In this section we present some applications of Theorem~\ref{mainSK} to the problems mentioned in Subsection~\ref{subsec:rw}.

\subsection{Group irregular labeling for directed graphs}\label{subsec:gild}
Let $\overrightarrow{G}$ be  a directed graph. 
Recall that, if there exists a mapping $\psi\colon E(\overrightarrow{G})\rightarrow\Gamma$ such that the mapping $\varphi_{\psi}\colon V(\overrightarrow{G})\rightarrow\Gamma$ defined by
$$\varphi_{\psi}(x)=\sum_{y\in N^-(x)}\psi(yx)-\sum_{y\in N^+(x)}\psi(xy),\;\;\;(x\in V(\overrightarrow{G}))$$
is injective, then we say that $\psi$ is a $\Gamma$-irregular labeling of $\overrightarrow{G}$. By Theorem~\ref{mainSK} we obtain the following:
\begin{cor} Let $\overrightarrow{G}$ be a directed graph of order $n$ with no component of order less than $3$, and let $\Gamma$ be a finite Abelian group with more than one involution such that $|\Gamma|=2^m$ for some $m$, $|\Gamma|> n$  and $|\Gamma|\not\in\{n+2,n+3\}$.  Then there exists a $\Gamma$-irregular labeling of $\overrightarrow{G}$.\end{cor}
\begin{proof}
Let $\{\overrightarrow{G_i}\}_{i=1}^t$ be the weakly connected components of $\overrightarrow{G}$. Recall that there exists a $\Gamma$-irregular labeling of a directed graph $\overrightarrow{G}$ with weakly connected components $\{\overrightarrow{G_i}\}_{i=1}^t$ if and only if there exist in $\Gamma$ pairwise disjoint subsets $\{S_i\}_{i=1}^t$ such that $|S_i|=|V(\overrightarrow{G_i})|$  and $\sum_{s\in S_i}s=0$ for every {$i\in[1 ,t]$} {(see Theorem 2.1 in \citet{CicTuz})}.

If $|\Gamma|=n+1$, then we can take the integer partition $\{m_i\}_{i=1}^t$ of $n$ with $m_i = |V(\overrightarrow{G}_i)|$ for every $i\in[1,t]$, and we get the result by Theorem~\ref{mainSK}.

Assume now that $|\Gamma|\geq n+4$. Let $m_i=|V(\overrightarrow{G}_i)|$ for every $i\in[1,t]$, and let $m_{t+1}=|\Gamma^*|-n$. Since $m_{t+1}\geq 3$, using the sequence $\{m_i\}_{i=1}^{t+1}$, by Theorem~\ref{mainSK}, again, we get the result.
\end{proof}

\subsection{\texorpdfstring{Realizable triples in $\Gamma$-anti-magic labeling}{Realizable triples in Gamma-anti-magic labeling}}\label{subsec:gaml}
Recall that, given an Abelian group $\Gamma$, a $\Gamma$-anti-magic labeling of a graph $G$ is a bijection $f\colon E(G)\to \Gamma^*$ such that the weight of every vertex (i.e. the sum of labels of incident edges) is unique. \citet{KLR} showed that, if $\Gamma$ has a unique involution, then every tree on $n$ vertices is not $\Gamma$-anti-magic. They conjectured that a tree with $|\Gamma|$ vertices admits a $\Gamma$-anti-magic labeling if and only if $\Gamma$ is not a group with a unique involution. Using the same method as \citet{KLR} for $2$-trees, by Theorem~\ref{mainSK}, we obtain the following:
\begin{cor}Every $3$-tree\footnote{A $3$-tree $T$ is a rooted tree, where every vertex which is not a leaf has at least three children.} $T$ of order $2^n$ admits a $\Gamma$-anti-magic labeling if and only if $\Gamma \not \cong \zet_{2^n}$.\end{cor}
\begin{proof} 
The necessity of the condition follows directly from the work of \citet{KLR}.
For the sufficiency, note that the only group  of order $2^n$ with a unique involution is the cyclic group $\zet_{2^n}$. Thus we can assume that $\Gamma \not \cong \zet_{2^n}$. Let $\{v_i\}_{i=1}^{t}$ be the vertices of $T$ which are not leaves. Let us denote their corresponding numbers of children by $\{m_i\}_{i=1}^{t}$. Thus $\sum_{i=1}^{t}m_i=2^n-1$. Since $T$ is a
$3$-tree, we have that $m_i \geq 3$ for every $i\in[1,t]$. By Theorem~\ref{mainSK}, $\Gamma$ has $3$-ZSPP, and there exists a zero-sum partition $\{S_i\}_{i=1}^{t}$ of $\Gamma^*$.
For all $i\in[1,t]$, we label the edges of the set $E_i=\{v_iw\colon w$ is a child of $v_i\}$ by the elements of $S_i$. The edges of $T$ are labeled bijectively with the non-zero elements of $\Gamma$, and the sum of
the labels in every $E_i$ is $0$. Since every vertex of $T$, except the root (whose weight is $0$), has a unique parent, it easily follows that the weights are pairwise distinct.
\end{proof}

\subsection{Group distance magic labeling}\label{subsec:gdml}
Recall that in group distance magic labeling, the vertices of the graph are labeled through a bijection with the elements of an Abelian group $\Gamma$. The weight of every vertex is computed as the sum (in $\Gamma$) of the labels assigned to its neighbors. If all the weights are the same, then it is a \textit{$\Gamma$-distance magic labeling}. 
% A graph that admits a $\Gamma$-distance magic labeling is a {$\Gamma$-distance magic graph} and the weigh is called \textit{magic constant}. 

Let $G=K_{m_1,\ldots,m_t}$ be a complete $t$-partite graph of order $m$. Let now $1\leq m_1\leq m_2 \leq\ldots\leq m_t$. Using some constant-sum partition properties of every finite Abelian group $\Gamma$ of order $m$, it was shown that:
\begin{enumerate}
\item  For $t=2$, if  $m_1 + m_2\not\equiv 2\pmod 4$ then the graph $G$ is  $\Gamma$-distance magic \citep{CichaczOM}.
	\item  For $t=3$, if ($m_2>1$ and  $m_1 + m_2+m_3\neq 2^p$ for any positive integer $p$) or  ($m_1\neq 2$ and $m_2>2$), then the graph $G$ admits a $\Gamma$-distance magic labeling \citep{Cic3}.
	\item If $m_1=m_2=\ldots=m_t>2$ and $|I(\Gamma)|\neq 1$, then the graph $G$ admits a $\Gamma$-distance magic labeling \citep{CicZ}.
	\item If $m_1\geq 3$ and $m_t\geq \frac{1}{2}(m +\sqrt{2m+1}) -1$, and $|I(\Gamma)|\neq 1$, then\/ $G$ admits a $\Gamma$-distance magic labeling \citep{CicTuz}.	
	\item  If $m_1\geq 4$, $m$ is large enough, and $|I(\Gamma)|\neq 1$, then\/ $G$ admits a $\Gamma$-distance magic labeling \citep{CicTuz},
	\item If $m_2\geq 2$, $t$ is odd, and $\Gamma\cong\zet_n$, then\/ $G$ admits a $\Gamma$-distance magic labeling \citep{FREYBERG2019}.
\end{enumerate}

Our main result allows us to prove the following corollary.

\begin{cor}\label{even} Let\/ $G = K_{m_1,m_2,\ldots,m_t}$ be a complete\/
 $t$-partite graph such that\/ $2^n = m_1 +m_2 +\ldots+m_{t}$ and\/
 $m_i\geq 3$ for every\/ $1\leq i \leq t$.
 Let\/ $\Gamma\not\cong \zet_{2^n}$ be an Abelian group of order\/ $2^n$, then the graph\/ $G$ admits a $\Gamma$-distance magic labeling.
  \end{cor}
\begin{proof}
Since $\Gamma\not\cong \zet_{2^n}$, and the required cartinalities are greater than or equal $3$, by Theorem~\ref{mainSK}, there exists a zero-sum partition $\{A'_i\}_{i=1}^{t}$ of $\Gamma^*$ such that $|A_t'|=m_t-1$ and $|A_i'| =m_i$ for every  {$i\in[1, t-1]$}.
Let $A_t=A'_t\cup \{0\}$ and $A_i=A_i'$ for every  {$i\in[1,t-1]$}. For every $i\in[1,t]$, we can label the vertices of $V_i$, where $V_i$ is the color class of cardinality $m_i$, using the elements of the set $A_i$, thus obtaining the required labeling.
\end{proof}\bsk

\section{Final remarks}\label{sec:final}
Our main result confirms Conjecture \ref{conjecture}, that any finite Abelian group with more than one involution has $3$-Zero-Sum Partition Property for the case of groups of order $2^n$. The other cases of groups of even order with more than one involution are still open (recall that the case of groups of odd order is completely solved \citep{Tannenbaum1,Zeng}). We believe that, using similar techniques (with some more effort) as in the proof of Theorem~\ref{mainSK}, Conjecture \ref{conjecture} could be confirmed for other groups of even order with more than one involution. The problem was partially solved in a very recent work by \citet{MP}, where the authors showed that, if $|\Gamma|$ is large enough and $|I(\Gamma)|>1$, then $\Gamma$ has $3$-ZSPP.

On the side of applications, Corollary \ref{even} covers all groups of order $2^n$, except for $\zet_{2^n}$. It would be interesting to see what happens in this case. We know the following sufficient and necessary condition for $\zet_m$-distance magic labeling of a complete $t$-partite graph only for $t$ or $m$ odd.

\begin{thm}[\citet{Cic3,FREYBERG2019}]
Let $G = K_{m_1,m_2,\ldots,m_t}$ be a complete $t$-partite graph of order $m$ such that\/ $ 1\leq m_1 \leq m_2 \leq\ldots\leq m_{t}$ and $m$ or $t$ is odd. The graph $G$ admits a $\zet_m$-distance magic labeling if and only if $m_2\geq 2$.
\end{thm}

For $m$ and $t$ even, we only know the following.
\begin{thm}[\citet{CichaczOM,CicZ}]
Let $G = K_{m_1,m_2,\ldots,m_t}$ be a complete $t$-partite graph of order $m$ such that $ 1\leq m_1 \leq m_2 \leq\ldots\leq m_{t}$ and $m$ and $t$ are even. If $t=2$, the graph $G$ admits a $\zet_m$-distance magic labeling if and only if $m\not\equiv2  \pmod 4$. If $ 2\leq m_1 = m_2 =\ldots=m_{t}$ then the graph $G$ admits a $\zet_m$-distance magic labeling if and only if $m_1$ is even. 
\end{thm}

Let $m=2^{n}(2k+1)$, $t=2^{n'}(2k'+1)$  for $n,n',k,k'\in \mathbb{N}$, $n\geq1$ and $n'\geq n$.
Note that, on the one hand, for every $a\in \zet_m$, we have $ta\not\equiv m/2\pmod m$. On the other hand, if the graph $G = K_{m_1,m_2,\ldots,m_t}$ is $\zet_m$-distance magic, with the magic constant $\mu\in\zet_m$, then $t\mu = (t-1)\sum_{g\in \Gamma}g \equiv m/2\pmod m$ - a contradiction. Observe that for $m=2^n$ and $t=2^{n'}(2k'+1)$ such that $m=\sum_{i=1}^{t}m_i$, with  $m_i\geq 2$ for every $i\in[1,t]$, there is $n'<n$. Thus we state the following conjecture:

\begin{conj}
Let $G = K_{m_1,m_2,\ldots,m_t}$ be a complete $t$-partite graph such that $2^n = \sum_{i=1}^{t} m_i$ and $m_i\geq 2$ for every $i\in[1,t]$. Then $G$ admits a $\zet_{2^n}$-distance magic labeling.
\end{conj}

We also state a similar conjecture in the language of Problem~\ref{problemT} of partitions of the group $\zet_{2^n}$.

\begin{conj}
Let $\zet_m$ be a group of order $m=2^n$. Let $t$ and $m_i$, $i\in[1,t]$, be positive integers  such that $\sum_{i=1}^t m_i=m-1$, $m_1\geq1$, $m_i\geq 2$ for every $i\in[2, t]$. Then there exists $\mu\in\zet_m$ such that the elements of $\zet_m^*$ can be partitioned into disjoint subsets $S_i$, $i\in[1,t]$, such that $|S_i|=m_i$ and $\sum_{s\in S_i}s=\mu$ for every $i\in[1,t]$.
\end{conj}

\textbf{Acknowledgement}\\
The work of the first author was partially supported by the Faculty of Applied Mathematics AGH UST statutory tasks within subsidy of Ministry of Education and Science. The authors would like to thank anonymous referees for theirs valuable comments.
\bibliographystyle{abbrvnat}
\bibliography{grupy}

\newpage

\begin{appendices}

Throughout the annexes, the notation ``\verb|a*3  b*4  c*5|'' refers to $a$ sets of size $3$, $b$ sets of size $4$, and $c$ sets of size $5$. In each partition, we present one set per line. Each tuple represents one element of the partitioned set, with the elements of the tuple following the order of the direct product.
\scriptsize

\section{\texorpdfstring{Zero sum partitions of $\left(\zet_4\times(\zet_2)^{2}\right)^*$ (with $b \leq 1$)}{Zero sum partitions of (Z4 x (Z2)**2)* (with b <= 1)}}\label{a:422}

\begin{verbatim}
[4, 2, 2]
[5, 0, 0]
[
[[0, 0, 1],[1, 0, 0],[3, 0, 1]],
[[0, 1, 0],[1, 0, 1],[3, 1, 1]],
[[2, 1, 0],[3, 0, 0],[3, 1, 0]],
[[0, 1, 1],[2, 0, 0],[2, 1, 1]],
[[1, 1, 0],[1, 1, 1],[2, 0, 1]]
]
A partition for sets of sizes:  5*3  0*4  0*5

[4, 2, 2]
[0, 0, 3]
[
[[0, 0, 1],[0, 1, 0],[0, 1, 1],[1, 0, 0],[3, 0, 0]],
[[1, 1, 0],[1, 1, 1],[2, 0, 0],[2, 1, 0],[2, 1, 1]],
[[1, 0, 1],[2, 0, 1],[3, 0, 1],[3, 1, 0],[3, 1, 1]]
]
A partition for sets of sizes:  0*3  0*4  3*5
\end{verbatim}

\section{\texorpdfstring{Zero sum partitions of $\left(\zet_4\times(\zet_2)^{3}\right)^*$ (with $b\leq 3$)}{Zero sum partitions of (Z4 x (Z2)**3)* (with b <= 3)}}\label{a:4222}
\begin{verbatim}
[4, 2, 2, 2]
[7, 0, 2]
[
[[0, 0, 0, 1],[0, 0, 1, 0],[0, 0, 1, 1]],
[[0, 1, 0, 0],[1, 0, 0, 0],[3, 1, 0, 0]],
[[0, 1, 1, 1],[1, 0, 0, 1],[3, 1, 1, 0]],
[[1, 0, 1, 0],[1, 0, 1, 1],[2, 0, 0, 1]],
[[1, 1, 0, 1],[1, 1, 1, 0],[2, 0, 1, 1]],
[[2, 0, 1, 0],[3, 0, 0, 0],[3, 0, 1, 0]],
[[2, 1, 0, 0],[3, 0, 0, 1],[3, 1, 0, 1]],
[[1, 1, 0, 0],[1, 1, 1, 1],[2, 0, 0, 0],[2, 1, 0, 1],[2, 1, 1, 0]],
[[0, 1, 0, 1],[0, 1, 1, 0],[2, 1, 1, 1],[3, 0, 1, 1],[3, 1, 1, 1]]
]
A partition for sets of sizes:  7*3  0*4  2*5

[4, 2, 2, 2]
[2, 0, 5]
[
[[0, 0, 0, 1],[0, 0, 1, 0],[0, 0, 1, 1]],
[[0, 1, 0, 0],[1, 0, 0, 0],[3, 1, 0, 0]],
[[0, 1, 1, 1],[1, 0, 0, 1],[1, 0, 1, 0],[1, 0, 1, 1],[1, 1, 1, 1]],
[[1, 1, 0, 0],[1, 1, 0, 1],[1, 1, 1, 0],[2, 0, 0, 0],[3, 1, 1, 1]],
[[2, 0, 0, 1],[2, 0, 1, 0],[2, 1, 0, 0],[3, 0, 0, 1],[3, 1, 1, 0]],
[[2, 0, 1, 1],[2, 1, 1, 0],[2, 1, 1, 1],[3, 0, 0, 0],[3, 0, 1, 0]],
[[0, 1, 0, 1],[0, 1, 1, 0],[2, 1, 0, 1],[3, 0, 1, 1],[3, 1, 0, 1]]
]
A partition for sets of sizes:  2*3  0*4  5*5
\end{verbatim}

\newpage
\section{\texorpdfstring{Zero sum partitions of $\left(\zet_4\times(\zet_2)^{4}\right)^*$ (with $b\leq 7$)}{Zero sum partitions of (Z4 x (Z2)**4)* (with b <= 7)}}\label{a:42222}

\begin{verbatim}
[4, 2, 2, 2, 2]
[21, 0, 0]
[
[[0, 0, 0, 0, 1],[0, 0, 0, 1, 0],[0, 0, 0, 1, 1]],
[[0, 0, 1, 0, 0],[0, 1, 0, 0, 0],[0, 1, 1, 0, 0]],
[[0, 0, 1, 1, 1],[0, 1, 0, 0, 1],[0, 1, 1, 1, 0]],
[[0, 1, 0, 1, 0],[1, 0, 0, 0, 0],[3, 1, 0, 1, 0]],
[[0, 1, 1, 0, 1],[1, 0, 0, 0, 1],[3, 1, 1, 0, 0]],
[[0, 1, 0, 1, 1],[1, 0, 0, 1, 0],[3, 1, 0, 0, 1]],
[[1, 0, 0, 1, 1],[1, 0, 1, 0, 0],[2, 0, 1, 1, 1]],
[[1, 0, 1, 1, 0],[1, 0, 1, 1, 1],[2, 0, 0, 0, 1]],
[[1, 1, 0, 0, 1],[1, 1, 0, 1, 0],[2, 0, 0, 1, 1]],
[[1, 1, 1, 0, 0],[1, 1, 1, 1, 0],[2, 0, 0, 1, 0]],
[[1, 1, 0, 0, 0],[1, 1, 1, 0, 1],[2, 0, 1, 0, 1]],
[[1, 1, 0, 1, 1],[1, 1, 1, 1, 1],[2, 0, 1, 0, 0]],
[[2, 0, 1, 1, 0],[3, 0, 0, 0, 0],[3, 0, 1, 1, 0]],
[[2, 1, 0, 0, 0],[3, 0, 1, 0, 1],[3, 1, 1, 0, 1]],
[[2, 1, 1, 0, 0],[3, 0, 1, 0, 0],[3, 1, 0, 0, 0]],
[[2, 1, 0, 0, 1],[3, 0, 1, 1, 1],[3, 1, 1, 1, 0]],
[[2, 1, 1, 0, 1],[3, 0, 0, 1, 0],[3, 1, 1, 1, 1]],
[[0, 0, 1, 0, 1],[2, 1, 0, 1, 1],[2, 1, 1, 1, 0]],
[[0, 0, 1, 1, 0],[1, 0, 1, 0, 1],[3, 0, 0, 1, 1]],
[[0, 1, 1, 1, 1],[2, 0, 0, 0, 0],[2, 1, 1, 1, 1]],
[[2, 1, 0, 1, 0],[3, 0, 0, 0, 1],[3, 1, 0, 1, 1]]
]
A partition for sets of sizes: 21*3  0*4  0*5

[4, 2, 2, 2, 2]
[18, 1, 1]
[
[[0, 0, 0, 0, 1],[0, 0, 0, 1, 0],[0, 0, 0, 1, 1]],
[[0, 0, 1, 0, 0],[0, 1, 0, 0, 0],[0, 1, 1, 0, 0]],
[[0, 0, 1, 1, 1],[0, 1, 0, 0, 1],[0, 1, 1, 1, 0]],
[[0, 1, 0, 1, 0],[1, 0, 0, 0, 0],[3, 1, 0, 1, 0]],
[[0, 1, 1, 0, 1],[1, 0, 0, 0, 1],[3, 1, 1, 0, 0]],
[[0, 1, 0, 1, 1],[1, 0, 0, 1, 0],[3, 1, 0, 0, 1]],
[[1, 0, 0, 1, 1],[1, 0, 1, 0, 0],[2, 0, 1, 1, 1]],
[[1, 0, 1, 1, 0],[1, 0, 1, 1, 1],[2, 0, 0, 0, 1]],
[[1, 1, 0, 0, 1],[1, 1, 0, 1, 0],[2, 0, 0, 1, 1]],
[[1, 1, 1, 0, 0],[1, 1, 1, 1, 0],[2, 0, 0, 1, 0]],
[[1, 1, 0, 0, 0],[1, 1, 1, 0, 1],[2, 0, 1, 0, 1]],
[[1, 1, 0, 1, 1],[1, 1, 1, 1, 1],[2, 0, 1, 0, 0]],
[[2, 0, 1, 1, 0],[3, 0, 0, 0, 0],[3, 0, 1, 1, 0]],
[[2, 1, 0, 0, 0],[3, 0, 0, 1, 1],[3, 1, 0, 1, 1]],
[[2, 1, 0, 1, 1],[3, 0, 1, 0, 0],[3, 1, 1, 1, 1]],
[[2, 1, 1, 1, 1],[3, 0, 0, 0, 1],[3, 1, 1, 1, 0]],
[[0, 0, 1, 0, 1],[2, 1, 0, 0, 1],[2, 1, 1, 0, 0]],
[[2, 1, 0, 1, 0],[3, 0, 1, 1, 1],[3, 1, 1, 0, 1]],
[[1, 0, 1, 0, 1],[2, 0, 0, 0, 0],[2, 1, 1, 0, 1],[3, 1, 0, 0, 0]],
[[0, 0, 1, 1, 0],[0, 1, 1, 1, 1],[2, 1, 1, 1, 0],[3, 0, 0, 1, 0],[3, 0, 1, 0, 1]]
]
A partition for sets of sizes: 18*3  1*4  1*5
\end{verbatim}

\newpage

\begin{verbatim}
[4, 2, 2, 2, 2]
[16, 0, 3]
[
[[0, 0, 0, 0, 1],[0, 0, 0, 1, 0],[0, 0, 0, 1, 1]],
[[0, 0, 1, 0, 0],[0, 1, 0, 0, 0],[0, 1, 1, 0, 0]],
[[0, 0, 1, 1, 1],[0, 1, 0, 0, 1],[0, 1, 1, 1, 0]],
[[0, 1, 0, 1, 0],[1, 0, 0, 0, 0],[3, 1, 0, 1, 0]],
[[0, 1, 1, 0, 1],[1, 0, 0, 0, 1],[3, 1, 1, 0, 0]],
[[0, 1, 0, 1, 1],[1, 0, 0, 1, 0],[3, 1, 0, 0, 1]],
[[1, 0, 0, 1, 1],[1, 0, 1, 0, 0],[2, 0, 1, 1, 1]],
[[1, 0, 1, 1, 0],[1, 0, 1, 1, 1],[2, 0, 0, 0, 1]],
[[1, 1, 0, 0, 1],[1, 1, 0, 1, 0],[2, 0, 0, 1, 1]],
[[1, 1, 1, 0, 0],[1, 1, 1, 1, 0],[2, 0, 0, 1, 0]],
[[1, 1, 0, 0, 0],[1, 1, 1, 0, 1],[2, 0, 1, 0, 1]],
[[1, 1, 0, 1, 1],[1, 1, 1, 1, 1],[2, 0, 1, 0, 0]],
[[2, 0, 1, 1, 0],[3, 0, 0, 0, 0],[3, 0, 1, 1, 0]],
[[2, 1, 0, 0, 0],[3, 0, 0, 1, 1],[3, 1, 0, 1, 1]],
[[2, 1, 0, 1, 1],[3, 0, 1, 0, 0],[3, 1, 1, 1, 1]],
[[2, 1, 1, 1, 1],[3, 0, 0, 0, 1],[3, 1, 1, 1, 0]],
[[2, 1, 0, 0, 1],[2, 1, 1, 0, 0],[2, 1, 1, 0, 1],[3, 0, 1, 0, 1],[3, 1, 1, 0, 1]],
[[0, 0, 1, 0, 1],[0, 1, 1, 1, 1],[2, 0, 0, 0, 0],[3, 0, 0, 1, 0],[3, 1, 0, 0, 0]],
[[0, 0, 1, 1, 0],[1, 0, 1, 0, 1],[2, 1, 0, 1, 0],[2, 1, 1, 1, 0],[3, 0, 1, 1, 1]]
]
A partition for sets of sizes: 16*3  0*4  3*5

[4, 2, 2, 2, 2]
[13, 1, 4]
[
[[0, 0, 0, 0, 1],[0, 0, 0, 1, 0],[0, 0, 0, 1, 1]],
[[0, 0, 1, 0, 0],[0, 1, 0, 0, 0],[0, 1, 1, 0, 0]],
[[0, 0, 1, 1, 1],[0, 1, 0, 0, 1],[0, 1, 1, 1, 0]],
[[0, 1, 0, 1, 0],[1, 0, 0, 0, 0],[3, 1, 0, 1, 0]],
[[0, 1, 1, 0, 1],[1, 0, 0, 0, 1],[3, 1, 1, 0, 0]],
[[0, 1, 0, 1, 1],[1, 0, 0, 1, 0],[3, 1, 0, 0, 1]],
[[1, 0, 0, 1, 1],[1, 0, 1, 0, 0],[2, 0, 1, 1, 1]],
[[1, 0, 1, 1, 0],[1, 0, 1, 1, 1],[2, 0, 0, 0, 1]],
[[1, 1, 0, 0, 1],[1, 1, 0, 1, 0],[2, 0, 0, 1, 1]],
[[1, 1, 1, 0, 0],[1, 1, 1, 1, 0],[2, 0, 0, 1, 0]],
[[1, 1, 0, 0, 0],[1, 1, 1, 0, 1],[2, 0, 1, 0, 1]],
[[1, 1, 0, 1, 1],[1, 1, 1, 1, 1],[2, 0, 1, 0, 0]],
[[2, 0, 1, 1, 0],[3, 0, 0, 0, 0],[3, 0, 1, 1, 0]],
[[2, 1, 0, 0, 0],[2, 1, 0, 0, 1],[2, 1, 0, 1, 0],[2, 1, 0, 1, 1]],
[[2, 1, 1, 0, 0],[2, 1, 1, 0, 1],[2, 1, 1, 1, 0],[3, 0, 0, 0, 1],[3, 1, 1, 1, 0]],
[[1, 0, 1, 0, 1],[2, 1, 1, 1, 1],[3, 0, 0, 1, 1],[3, 0, 1, 0, 0],[3, 1, 1, 0, 1]],
[[0, 0, 1, 0, 1],[0, 0, 1, 1, 0],[2, 0, 0, 0, 0],[3, 1, 0, 0, 0],[3, 1, 0, 1, 1]],
[[0, 1, 1, 1, 1],[3, 0, 0, 1, 0],[3, 0, 1, 0, 1],[3, 0, 1, 1, 1],[3, 1, 1, 1, 1]]
]
A partition for sets of sizes: 13*3  1*4  4*5
\end{verbatim}

\newpage

\begin{verbatim}
[4, 2, 2, 2, 2]
[11, 0, 6]
[
[[0, 0, 0, 0, 1],[0, 0, 0, 1, 0],[0, 0, 0, 1, 1]],
[[0, 0, 1, 0, 0],[0, 1, 0, 0, 0],[0, 1, 1, 0, 0]],
[[0, 0, 1, 1, 1],[0, 1, 0, 0, 1],[0, 1, 1, 1, 0]],
[[0, 1, 0, 1, 0],[1, 0, 0, 0, 0],[3, 1, 0, 1, 0]],
[[0, 1, 1, 0, 1],[1, 0, 0, 0, 1],[3, 1, 1, 0, 0]],
[[0, 1, 0, 1, 1],[1, 0, 0, 1, 0],[3, 1, 0, 0, 1]],
[[1, 0, 0, 1, 1],[1, 0, 1, 0, 0],[2, 0, 1, 1, 1]],
[[1, 0, 1, 1, 0],[1, 0, 1, 1, 1],[2, 0, 0, 0, 1]],
[[1, 1, 0, 0, 1],[1, 1, 0, 1, 0],[2, 0, 0, 1, 1]],
[[1, 1, 1, 0, 0],[1, 1, 1, 1, 0],[2, 0, 0, 1, 0]],
[[1, 1, 0, 0, 0],[1, 1, 1, 0, 1],[2, 0, 1, 0, 1]],
[[2, 0, 0, 0, 0],[2, 0, 1, 0, 0],[2, 0, 1, 1, 0],[3, 0, 0, 0, 0],[3, 0, 0, 1, 0]],
[[1, 0, 1, 0, 1],[2, 1, 0, 0, 0],[3, 0, 0, 0, 1],[3, 0, 0, 1, 1],[3, 1, 1, 1, 1]],
[[2, 1, 1, 0, 0],[2, 1, 1, 0, 1],[2, 1, 1, 1, 0],[3, 0, 1, 0, 0],[3, 1, 0, 1, 1]],
[[2, 1, 0, 0, 1],[2, 1, 0, 1, 0],[2, 1, 0, 1, 1],[3, 0, 1, 0, 1],[3, 1, 1, 0, 1]],
[[0, 0, 1, 0, 1],[0, 0, 1, 1, 0],[0, 1, 1, 1, 1],[1, 1, 0, 1, 1],[3, 0, 1, 1, 1]],
[[1, 1, 1, 1, 1],[2, 1, 1, 1, 1],[3, 0, 1, 1, 0],[3, 1, 0, 0, 0],[3, 1, 1, 1, 0]]
]
A partition for sets of sizes: 11*3  0*4  6*5

[4, 2, 2, 2, 2]
[8, 1, 7]
[
[[0, 0, 0, 0, 1],[0, 0, 0, 1, 0],[0, 0, 0, 1, 1]],
[[0, 0, 1, 0, 0],[0, 1, 0, 0, 0],[0, 1, 1, 0, 0]],
[[0, 0, 1, 1, 1],[0, 1, 0, 0, 1],[0, 1, 1, 1, 0]],
[[0, 1, 0, 1, 0],[1, 0, 0, 0, 0],[3, 1, 0, 1, 0]],
[[0, 1, 1, 0, 1],[1, 0, 0, 0, 1],[3, 1, 1, 0, 0]],
[[0, 1, 0, 1, 1],[1, 0, 0, 1, 0],[3, 1, 0, 0, 1]],
[[1, 0, 0, 1, 1],[1, 0, 1, 0, 0],[2, 0, 1, 1, 1]],
[[1, 0, 1, 1, 0],[1, 0, 1, 1, 1],[2, 0, 0, 0, 1]],
[[1, 1, 0, 0, 1],[1, 1, 0, 1, 0],[1, 1, 1, 0, 0],[1, 1, 1, 1, 1]],
[[1, 1, 1, 0, 1],[1, 1, 1, 1, 0],[2, 0, 0, 0, 0],[2, 0, 1, 0, 1],[2, 0, 1, 1, 0]],
[[2, 0, 0, 1, 0],[2, 0, 0, 1, 1],[2, 0, 1, 0, 0],[3, 0, 0, 0, 0],[3, 0, 1, 0, 1]],
[[1, 0, 1, 0, 1],[2, 1, 0, 0, 0],[3, 0, 0, 0, 1],[3, 0, 0, 1, 0],[3, 1, 1, 1, 0]],
[[2, 1, 1, 0, 0],[2, 1, 1, 0, 1],[2, 1, 1, 1, 0],[3, 0, 1, 0, 0],[3, 1, 0, 1, 1]],
[[2, 1, 0, 0, 1],[2, 1, 0, 1, 0],[2, 1, 1, 1, 1],[3, 0, 0, 1, 1],[3, 1, 1, 1, 1]],
[[0, 0, 1, 0, 1],[0, 0, 1, 1, 0],[0, 1, 1, 1, 1],[1, 1, 0, 1, 1],[3, 0, 1, 1, 1]],
[[1, 1, 0, 0, 0],[2, 1, 0, 1, 1],[3, 0, 1, 1, 0],[3, 1, 0, 0, 0],[3, 1, 1, 0, 1]]
]
A partition for sets of sizes:  8*3  1*4  7*5

[4, 2, 2, 2, 2]
[6, 0, 9]
[
[[0, 0, 0, 0, 1],[0, 0, 0, 1, 0],[0, 0, 0, 1, 1]],
[[0, 0, 1, 0, 0],[0, 1, 0, 0, 0],[0, 1, 1, 0, 0]],
[[0, 0, 1, 1, 1],[0, 1, 0, 0, 1],[0, 1, 1, 1, 0]],
[[0, 1, 0, 1, 0],[1, 0, 0, 0, 0],[3, 1, 0, 1, 0]],
[[0, 1, 1, 0, 1],[1, 0, 0, 0, 1],[3, 1, 1, 0, 0]],
[[0, 1, 0, 1, 1],[1, 0, 0, 1, 0],[3, 1, 0, 0, 1]],
[[1, 0, 0, 1, 1],[1, 0, 1, 0, 0],[1, 0, 1, 0, 1],[2, 0, 0, 0, 0],[3, 0, 0, 1, 0]],
[[1, 1, 0, 0, 0],[1, 1, 0, 0, 1],[1, 1, 0, 1, 0],[2, 0, 0, 1, 1],[3, 1, 0, 0, 0]],
[[1, 1, 1, 0, 1],[1, 1, 1, 1, 0],[1, 1, 1, 1, 1],[2, 0, 0, 0, 1],[3, 1, 1, 0, 1]],
[[2, 0, 0, 1, 0],[2, 0, 1, 0, 0],[2, 0, 1, 0, 1],[3, 0, 0, 0, 0],[3, 0, 0, 1, 1]],
[[2, 0, 1, 1, 1],[2, 1, 0, 0, 0],[2, 1, 0, 0, 1],[3, 0, 0, 0, 1],[3, 0, 1, 1, 1]],
[[2, 1, 1, 0, 0],[2, 1, 1, 0, 1],[2, 1, 1, 1, 0],[3, 0, 1, 0, 0],[3, 1, 0, 1, 1]],
[[1, 0, 1, 1, 1],[2, 1, 0, 1, 1],[3, 0, 1, 0, 1],[3, 0, 1, 1, 0],[3, 1, 1, 1, 1]],
[[0, 0, 1, 0, 1],[0, 0, 1, 1, 0],[0, 1, 1, 1, 1],[2, 0, 1, 1, 0],[2, 1, 0, 1, 0]],
[[1, 0, 1, 1, 0],[1, 1, 0, 1, 1],[1, 1, 1, 0, 0],[2, 1, 1, 1, 1],[3, 1, 1, 1, 0]]
]
A partition for sets of sizes:  6*3  0*4  9*5

[4, 2, 2, 2, 2]
[3, 1, 10]
[
[[0, 0, 0, 0, 1],[0, 0, 0, 1, 0],[0, 0, 0, 1, 1]],
[[0, 0, 1, 0, 0],[0, 1, 0, 0, 0],[0, 1, 1, 0, 0]],
[[0, 0, 1, 1, 1],[0, 1, 0, 0, 1],[0, 1, 1, 1, 0]],
[[0, 1, 0, 1, 0],[0, 1, 0, 1, 1],[1, 0, 0, 0, 0],[3, 0, 0, 0, 1]],
[[0, 0, 1, 0, 1],[0, 1, 1, 1, 1],[1, 0, 0, 0, 1],[1, 0, 0, 1, 0],[2, 1, 0, 0, 1]],
[[1, 0, 0, 1, 1],[1, 0, 1, 0, 0],[1, 0, 1, 0, 1],[2, 0, 0, 0, 0],[3, 0, 0, 1, 0]],
[[1, 1, 0, 0, 0],[1, 1, 0, 0, 1],[1, 1, 0, 1, 0],[2, 0, 0, 0, 1],[3, 1, 0, 1, 0]],
[[1, 1, 1, 0, 1],[1, 1, 1, 1, 0],[1, 1, 1, 1, 1],[2, 0, 0, 1, 0],[3, 1, 1, 1, 0]],
[[1, 0, 1, 1, 0],[2, 0, 0, 1, 1],[3, 0, 0, 0, 0],[3, 0, 0, 1, 1],[3, 0, 1, 1, 0]],
[[2, 0, 1, 1, 1],[2, 1, 0, 0, 0],[2, 1, 0, 1, 0],[3, 1, 0, 0, 0],[3, 1, 1, 0, 1]],
[[2, 1, 1, 0, 0],[2, 1, 1, 0, 1],[2, 1, 1, 1, 0],[3, 0, 1, 0, 0],[3, 1, 0, 1, 1]],
[[0, 1, 1, 0, 1],[1, 0, 1, 1, 1],[2, 0, 1, 0, 1],[2, 0, 1, 1, 0],[3, 1, 0, 0, 1]],
[[0, 0, 1, 1, 0],[1, 1, 1, 0, 0],[2, 0, 1, 0, 0],[2, 1, 0, 1, 1],[3, 0, 1, 0, 1]],
[[1, 1, 0, 1, 1],[2, 1, 1, 1, 1],[3, 0, 1, 1, 1],[3, 1, 1, 0, 0],[3, 1, 1, 1, 1]]
]
A partition for sets of sizes:  3*3  1*4 10*5


[4, 2, 2, 2, 2]
[1, 0, 12]
[
[[0, 0, 0, 0, 1],[0, 0, 0, 1, 0],[0, 0, 0, 1, 1]],
[[0, 0, 1, 0, 0],[0, 0, 1, 0, 1],[0, 0, 1, 1, 0],[0, 1, 0, 0, 0],[0, 1, 1, 1, 1]],
[[0, 1, 0, 0, 1],[0, 1, 0, 1, 0],[0, 1, 0, 1, 1],[1, 0, 0, 0, 0],[3, 1, 0, 0, 0]],
[[0, 1, 1, 1, 0],[1, 0, 0, 0, 1],[1, 0, 0, 1, 0],[1, 0, 0, 1, 1],[1, 1, 1, 1, 0]],
[[0, 0, 1, 1, 1],[1, 0, 1, 0, 0],[1, 0, 1, 0, 1],[1, 1, 0, 0, 1],[1, 1, 1, 1, 1]],
[[1, 1, 0, 0, 0],[1, 1, 0, 1, 0],[1, 1, 0, 1, 1],[2, 0, 0, 0, 0],[3, 1, 0, 0, 1]],
[[1, 1, 1, 0, 1],[2, 0, 0, 0, 1],[3, 0, 0, 0, 0],[3, 0, 0, 0, 1],[3, 1, 1, 0, 1]],
[[2, 0, 0, 1, 0],[2, 0, 0, 1, 1],[2, 0, 1, 0, 0],[3, 0, 0, 1, 0],[3, 0, 1, 1, 1]],
[[2, 0, 1, 1, 1],[2, 1, 0, 0, 0],[2, 1, 0, 0, 1],[3, 0, 0, 1, 1],[3, 0, 1, 0, 1]],
[[2, 1, 1, 0, 0],[2, 1, 1, 0, 1],[2, 1, 1, 1, 0],[3, 0, 1, 0, 0],[3, 1, 0, 1, 1]],
[[1, 0, 1, 1, 0],[2, 0, 1, 0, 1],[3, 0, 1, 1, 0],[3, 1, 0, 1, 0],[3, 1, 1, 1, 1]],
[[2, 0, 1, 1, 0],[2, 1, 0, 1, 1],[2, 1, 1, 1, 1],[3, 1, 1, 0, 0],[3, 1, 1, 1, 0]],
[[0, 1, 1, 0, 0],[0, 1, 1, 0, 1],[1, 0, 1, 1, 1],[1, 1, 1, 0, 0],[2, 1, 0, 1, 0]]
]
A partition for sets of sizes:  1*3  0*4 12*5
\end{verbatim}

\newpage
\section{\texorpdfstring{Zero sum partitions of $\left(\zet_4\times(\zet_2)^{5}\right)^*$ (with with $b\leq 15$)}{Zero sum partitions of (Z4 x (Z2)**5)* (with with b <= 15)}}\label{a:422222}

\begin{verbatim} 
[4, 2, 2, 2, 2, 2]
[41, 1, 0]
[
[[0, 0, 0, 0, 0, 1],[0, 0, 0, 0, 1, 0],[0, 0, 0, 0, 1, 1]],
[[0, 0, 0, 1, 0, 0],[0, 0, 1, 0, 0, 0],[0, 0, 1, 1, 0, 0]],
[[0, 0, 0, 1, 1, 1],[0, 0, 1, 0, 0, 1],[0, 0, 1, 1, 1, 0]],
[[0, 0, 1, 0, 1, 0],[0, 1, 0, 0, 0, 0],[0, 1, 1, 0, 1, 0]],
[[0, 0, 1, 1, 0, 1],[0, 1, 0, 0, 0, 1],[0, 1, 1, 1, 0, 0]],
[[0, 0, 1, 0, 1, 1],[0, 1, 0, 0, 1, 0],[0, 1, 1, 0, 0, 1]],
[[0, 1, 0, 0, 1, 1],[1, 0, 0, 0, 0, 0],[3, 1, 0, 0, 1, 1]],
[[0, 1, 0, 1, 1, 0],[1, 0, 0, 0, 0, 1],[3, 1, 0, 1, 1, 1]],
[[0, 0, 0, 1, 0, 1],[0, 1, 1, 0, 1, 1],[0, 1, 1, 1, 1, 0]],
[[0, 0, 1, 1, 1, 1],[0, 1, 0, 1, 1, 1],[0, 1, 1, 0, 0, 0]],
[[0, 1, 1, 1, 1, 1],[1, 0, 0, 0, 1, 0],[3, 1, 1, 1, 0, 1]],
[[0, 1, 0, 1, 0, 0],[1, 0, 0, 1, 0, 0],[3, 1, 0, 0, 0, 0]],
[[1, 0, 0, 1, 0, 1],[1, 0, 0, 1, 1, 0],[2, 0, 0, 0, 1, 1]],
[[1, 0, 1, 0, 0, 0],[1, 0, 1, 0, 0, 1],[2, 0, 0, 0, 0, 1]],
[[1, 0, 1, 0, 1, 1],[1, 0, 1, 1, 0, 0],[2, 0, 0, 1, 1, 1]],
[[1, 0, 1, 1, 1, 0],[1, 1, 0, 0, 0, 0],[2, 1, 1, 1, 1, 0]],
[[1, 1, 0, 0, 0, 1],[1, 1, 0, 0, 1, 1],[2, 0, 0, 0, 1, 0]],
[[1, 1, 0, 1, 0, 0],[1, 1, 1, 0, 0, 0],[2, 0, 1, 1, 0, 0]],
[[1, 1, 0, 1, 1, 1],[1, 1, 1, 0, 0, 1],[2, 0, 1, 1, 1, 0]],
[[1, 1, 1, 0, 1, 0],[1, 1, 1, 1, 0, 0],[2, 0, 0, 1, 1, 0]],
[[1, 0, 1, 0, 1, 0],[1, 1, 1, 1, 1, 0],[2, 1, 0, 1, 0, 0]],
[[1, 1, 0, 0, 1, 0],[1, 1, 1, 1, 0, 1],[2, 0, 1, 1, 1, 1]],
[[1, 0, 0, 1, 1, 1],[1, 1, 1, 0, 1, 1],[2, 1, 1, 1, 0, 0]],
[[2, 0, 0, 1, 0, 0],[3, 0, 0, 0, 0, 0],[3, 0, 0, 1, 0, 0]],
[[2, 0, 1, 0, 0, 1],[3, 0, 0, 0, 0, 1],[3, 0, 1, 0, 0, 0]],
[[2, 0, 1, 1, 0, 1],[3, 0, 0, 0, 1, 0],[3, 0, 1, 1, 1, 1]],
[[2, 1, 0, 0, 0, 0],[3, 0, 0, 1, 0, 1],[3, 1, 0, 1, 0, 1]],
[[2, 1, 0, 0, 1, 0],[3, 0, 0, 0, 1, 1],[3, 1, 0, 0, 0, 1]],
[[2, 1, 0, 1, 0, 1],[3, 0, 0, 1, 1, 1],[3, 1, 0, 0, 1, 0]],
[[2, 1, 1, 0, 0, 0],[3, 0, 0, 1, 1, 0],[3, 1, 1, 1, 1, 0]],
[[2, 1, 1, 0, 1, 1],[3, 0, 1, 1, 0, 1],[3, 1, 0, 1, 1, 0]],
[[1, 0, 1, 1, 1, 1],[1, 1, 0, 1, 1, 0],[2, 1, 1, 0, 0, 1]],
[[2, 1, 0, 0, 1, 1],[3, 0, 1, 0, 1, 0],[3, 1, 1, 0, 0, 1]],
[[2, 0, 1, 0, 0, 0],[3, 1, 0, 1, 0, 0],[3, 1, 1, 1, 0, 0]],
[[2, 0, 0, 1, 0, 1],[3, 1, 1, 0, 1, 0],[3, 1, 1, 1, 1, 1]],
[[1, 0, 0, 0, 1, 1],[1, 1, 0, 1, 0, 1],[2, 1, 0, 1, 1, 0]],
[[2, 1, 0, 1, 1, 1],[3, 0, 1, 1, 0, 0],[3, 1, 1, 0, 1, 1]],
[[0, 1, 0, 1, 0, 1],[2, 0, 1, 0, 1, 0],[2, 1, 1, 1, 1, 1]],
[[2, 1, 0, 0, 0, 1],[3, 0, 1, 0, 0, 1],[3, 1, 1, 0, 0, 0]],
[[0, 0, 0, 1, 1, 0],[1, 0, 1, 1, 0, 1],[3, 0, 1, 0, 1, 1]],
[[0, 1, 1, 1, 0, 1],[2, 0, 0, 0, 0, 0],[2, 1, 1, 1, 0, 1]],
[[1, 1, 1, 1, 1, 1],[2, 0, 1, 0, 1, 1],[2, 1, 1, 0, 1, 0],[3, 0, 1, 1, 1, 0]]
]
A partition for sets of sizes: 41*3  1*4  0*5
\end{verbatim}

\newpage

\begin{verbatim}
[4, 2, 2, 2, 2, 2]
[38, 2, 1]
[
[[0, 0, 0, 0, 0, 1],[0, 0, 0, 0, 1, 0],[0, 0, 0, 0, 1, 1]],
[[0, 0, 0, 1, 0, 0],[0, 0, 1, 0, 0, 0],[0, 0, 1, 1, 0, 0]],
[[0, 0, 0, 1, 1, 1],[0, 0, 1, 0, 0, 1],[0, 0, 1, 1, 1, 0]],
[[0, 0, 1, 0, 1, 0],[0, 1, 0, 0, 0, 0],[0, 1, 1, 0, 1, 0]],
[[0, 0, 1, 1, 0, 1],[0, 1, 0, 0, 0, 1],[0, 1, 1, 1, 0, 0]],
[[0, 0, 1, 0, 1, 1],[0, 1, 0, 0, 1, 0],[0, 1, 1, 0, 0, 1]],
[[0, 1, 0, 0, 1, 1],[1, 0, 0, 0, 0, 0],[3, 1, 0, 0, 1, 1]],
[[0, 1, 0, 1, 1, 0],[1, 0, 0, 0, 0, 1],[3, 1, 0, 1, 1, 1]],
[[0, 0, 0, 1, 0, 1],[0, 1, 1, 0, 1, 1],[0, 1, 1, 1, 1, 0]],
[[0, 0, 1, 1, 1, 1],[0, 1, 0, 1, 1, 1],[0, 1, 1, 0, 0, 0]],
[[0, 1, 1, 1, 1, 1],[1, 0, 0, 0, 1, 0],[3, 1, 1, 1, 0, 1]],
[[0, 1, 0, 1, 0, 0],[1, 0, 0, 1, 0, 0],[3, 1, 0, 0, 0, 0]],
[[1, 0, 0, 1, 0, 1],[1, 0, 0, 1, 1, 0],[2, 0, 0, 0, 1, 1]],
[[1, 0, 1, 0, 0, 0],[1, 0, 1, 0, 0, 1],[2, 0, 0, 0, 0, 1]],
[[1, 0, 1, 0, 1, 1],[1, 0, 1, 1, 0, 0],[2, 0, 0, 1, 1, 1]],
[[1, 0, 1, 1, 1, 0],[1, 1, 0, 0, 0, 0],[2, 1, 1, 1, 1, 0]],
[[1, 1, 0, 0, 0, 1],[1, 1, 0, 0, 1, 1],[2, 0, 0, 0, 1, 0]],
[[1, 1, 0, 1, 0, 0],[1, 1, 1, 0, 0, 0],[2, 0, 1, 1, 0, 0]],
[[1, 1, 0, 1, 1, 1],[1, 1, 1, 0, 0, 1],[2, 0, 1, 1, 1, 0]],
[[1, 1, 1, 0, 1, 0],[1, 1, 1, 1, 0, 0],[2, 0, 0, 1, 1, 0]],
[[1, 0, 1, 0, 1, 0],[1, 1, 1, 1, 1, 0],[2, 1, 0, 1, 0, 0]],
[[1, 1, 0, 0, 1, 0],[1, 1, 1, 1, 0, 1],[2, 0, 1, 1, 1, 1]],
[[1, 0, 0, 1, 1, 1],[1, 1, 1, 0, 1, 1],[2, 1, 1, 1, 0, 0]],
[[2, 0, 0, 1, 0, 0],[3, 0, 0, 0, 0, 0],[3, 0, 0, 1, 0, 0]],
[[2, 0, 1, 0, 0, 1],[3, 0, 0, 0, 0, 1],[3, 0, 1, 0, 0, 0]],
[[2, 0, 1, 1, 0, 1],[3, 0, 0, 0, 1, 0],[3, 0, 1, 1, 1, 1]],
[[2, 1, 0, 0, 0, 0],[3, 0, 0, 1, 0, 1],[3, 1, 0, 1, 0, 1]],
[[2, 1, 0, 0, 1, 0],[3, 0, 0, 0, 1, 1],[3, 1, 0, 0, 0, 1]],
[[2, 1, 0, 1, 0, 1],[3, 0, 0, 1, 1, 1],[3, 1, 0, 0, 1, 0]],
[[2, 1, 1, 0, 0, 0],[3, 0, 0, 1, 1, 0],[3, 1, 1, 1, 1, 0]],
[[2, 1, 1, 0, 1, 1],[3, 0, 1, 1, 0, 1],[3, 1, 0, 1, 1, 0]],
[[1, 0, 1, 1, 1, 1],[1, 1, 0, 1, 1, 0],[2, 1, 1, 0, 0, 1]],
[[2, 1, 1, 1, 1, 1],[3, 0, 1, 0, 1, 1],[3, 1, 0, 1, 0, 0]],
[[2, 1, 0, 1, 1, 0],[3, 0, 1, 0, 0, 1],[3, 1, 1, 1, 1, 1]],
[[2, 0, 0, 1, 0, 1],[3, 1, 1, 0, 0, 1],[3, 1, 1, 1, 0, 0]],
[[1, 1, 0, 1, 0, 1],[1, 1, 1, 1, 1, 1],[2, 0, 1, 0, 1, 0]],
[[2, 1, 0, 1, 1, 1],[3, 0, 1, 1, 0, 0],[3, 1, 1, 0, 1, 1]],
[[0, 1, 0, 1, 0, 1],[2, 0, 1, 0, 0, 0],[2, 1, 1, 1, 0, 1]],
[[1, 0, 0, 0, 1, 1],[2, 1, 0, 0, 1, 1],[2, 1, 1, 0, 1, 0],[3, 0, 1, 0, 1, 0]],
[[0, 1, 1, 1, 0, 1],[2, 0, 1, 0, 1, 1],[3, 0, 1, 1, 1, 0],[3, 1, 1, 0, 0, 0]],
[[0, 0, 0, 1, 1, 0],[1, 0, 1, 1, 0, 1],[2, 0, 0, 0, 0, 0],[2, 1, 0, 0, 0, 1],[3, 1, 1, 0, 1, 0]]
]
A partition for sets of sizes: 38*3  2*4  1*5
\end{verbatim}

\newpage

\begin{verbatim}
[4, 2, 2, 2, 2, 2]
[39, 0, 2]
[
[[0, 0, 0, 0, 0, 1],[0, 0, 0, 0, 1, 0],[0, 0, 0, 0, 1, 1]],
[[0, 0, 0, 1, 0, 0],[0, 0, 1, 0, 0, 0],[0, 0, 1, 1, 0, 0]],
[[0, 0, 0, 1, 1, 1],[0, 0, 1, 0, 0, 1],[0, 0, 1, 1, 1, 0]],
[[0, 0, 1, 0, 1, 0],[0, 1, 0, 0, 0, 0],[0, 1, 1, 0, 1, 0]],
[[0, 0, 1, 1, 0, 1],[0, 1, 0, 0, 0, 1],[0, 1, 1, 1, 0, 0]],
[[0, 0, 1, 0, 1, 1],[0, 1, 0, 0, 1, 0],[0, 1, 1, 0, 0, 1]],
[[0, 1, 0, 0, 1, 1],[1, 0, 0, 0, 0, 0],[3, 1, 0, 0, 1, 1]],
[[0, 1, 0, 1, 1, 0],[1, 0, 0, 0, 0, 1],[3, 1, 0, 1, 1, 1]],
[[0, 0, 0, 1, 0, 1],[0, 1, 1, 0, 1, 1],[0, 1, 1, 1, 1, 0]],
[[0, 0, 1, 1, 1, 1],[0, 1, 0, 1, 1, 1],[0, 1, 1, 0, 0, 0]],
[[0, 1, 1, 1, 1, 1],[1, 0, 0, 0, 1, 0],[3, 1, 1, 1, 0, 1]],
[[0, 1, 0, 1, 0, 0],[1, 0, 0, 1, 0, 0],[3, 1, 0, 0, 0, 0]],
[[1, 0, 0, 1, 0, 1],[1, 0, 0, 1, 1, 0],[2, 0, 0, 0, 1, 1]],
[[1, 0, 1, 0, 0, 0],[1, 0, 1, 0, 0, 1],[2, 0, 0, 0, 0, 1]],
[[1, 0, 1, 0, 1, 1],[1, 0, 1, 1, 0, 0],[2, 0, 0, 1, 1, 1]],
[[1, 0, 1, 1, 1, 0],[1, 1, 0, 0, 0, 0],[2, 1, 1, 1, 1, 0]],
[[1, 1, 0, 0, 0, 1],[1, 1, 0, 0, 1, 1],[2, 0, 0, 0, 1, 0]],
[[1, 1, 0, 1, 0, 0],[1, 1, 1, 0, 0, 0],[2, 0, 1, 1, 0, 0]],
[[1, 1, 0, 1, 1, 1],[1, 1, 1, 0, 0, 1],[2, 0, 1, 1, 1, 0]],
[[1, 1, 1, 0, 1, 0],[1, 1, 1, 1, 0, 0],[2, 0, 0, 1, 1, 0]],
[[1, 0, 1, 0, 1, 0],[1, 1, 1, 1, 1, 0],[2, 1, 0, 1, 0, 0]],
[[1, 1, 0, 0, 1, 0],[1, 1, 1, 1, 0, 1],[2, 0, 1, 1, 1, 1]],
[[1, 0, 0, 1, 1, 1],[1, 1, 1, 0, 1, 1],[2, 1, 1, 1, 0, 0]],
[[2, 0, 0, 1, 0, 0],[3, 0, 0, 0, 0, 0],[3, 0, 0, 1, 0, 0]],
[[2, 0, 1, 0, 0, 1],[3, 0, 0, 0, 0, 1],[3, 0, 1, 0, 0, 0]],
[[2, 0, 1, 1, 0, 1],[3, 0, 0, 0, 1, 0],[3, 0, 1, 1, 1, 1]],
[[2, 1, 0, 0, 0, 0],[3, 0, 0, 1, 0, 1],[3, 1, 0, 1, 0, 1]],
[[2, 1, 0, 0, 1, 0],[3, 0, 0, 0, 1, 1],[3, 1, 0, 0, 0, 1]],
[[2, 1, 0, 1, 0, 1],[3, 0, 0, 1, 1, 1],[3, 1, 0, 0, 1, 0]],
[[2, 1, 1, 0, 0, 0],[3, 0, 0, 1, 1, 0],[3, 1, 1, 1, 1, 0]],
[[2, 1, 1, 0, 1, 1],[3, 0, 1, 1, 0, 1],[3, 1, 0, 1, 1, 0]],
[[1, 0, 1, 1, 1, 1],[1, 1, 0, 1, 1, 0],[2, 1, 1, 0, 0, 1]],
[[2, 1, 1, 1, 1, 1],[3, 0, 1, 0, 1, 1],[3, 1, 0, 1, 0, 0]],
[[2, 1, 0, 1, 1, 0],[3, 0, 1, 0, 0, 1],[3, 1, 1, 1, 1, 1]],
[[2, 0, 0, 1, 0, 1],[3, 1, 1, 0, 0, 1],[3, 1, 1, 1, 0, 0]],
[[1, 1, 0, 1, 0, 1],[1, 1, 1, 1, 1, 1],[2, 0, 1, 0, 1, 0]],
[[2, 1, 0, 1, 1, 1],[3, 0, 1, 1, 0, 0],[3, 1, 1, 0, 1, 1]],
[[0, 1, 0, 1, 0, 1],[1, 0, 1, 1, 0, 1],[3, 1, 1, 0, 0, 0]],
[[0, 1, 1, 1, 0, 1],[2, 0, 0, 0, 0, 0],[2, 1, 1, 1, 0, 1]],
[[0, 0, 0, 1, 1, 0],[1, 0, 0, 0, 1, 1],[2, 1, 0, 0, 0, 1],[2, 1, 1, 0, 1, 0],[3, 0, 1, 1, 1, 0]],
[[2, 0, 1, 0, 0, 0],[2, 0, 1, 0, 1, 1],[2, 1, 0, 0, 1, 1],[3, 0, 1, 0, 1, 0],[3, 1, 1, 0, 1, 0]]
]
A partition for sets of sizes: 39*3  0*4  2*5
\end{verbatim}
\newpage

\begin{verbatim} 
[4, 2, 2, 2, 2, 2]
[35, 3, 2]
[
[[0, 0, 0, 0, 0, 1],[0, 0, 0, 0, 1, 0],[0, 0, 0, 0, 1, 1]],
[[0, 0, 0, 1, 0, 0],[0, 0, 1, 0, 0, 0],[0, 0, 1, 1, 0, 0]],
[[0, 0, 0, 1, 1, 1],[0, 0, 1, 0, 0, 1],[0, 0, 1, 1, 1, 0]],
[[0, 0, 1, 0, 1, 0],[0, 1, 0, 0, 0, 0],[0, 1, 1, 0, 1, 0]],
[[0, 0, 1, 1, 0, 1],[0, 1, 0, 0, 0, 1],[0, 1, 1, 1, 0, 0]],
[[0, 0, 1, 0, 1, 1],[0, 1, 0, 0, 1, 0],[0, 1, 1, 0, 0, 1]],
[[0, 1, 0, 0, 1, 1],[1, 0, 0, 0, 0, 0],[3, 1, 0, 0, 1, 1]],
[[0, 1, 0, 1, 1, 0],[1, 0, 0, 0, 0, 1],[3, 1, 0, 1, 1, 1]],
[[0, 0, 0, 1, 0, 1],[0, 1, 1, 0, 1, 1],[0, 1, 1, 1, 1, 0]],
[[0, 0, 1, 1, 1, 1],[0, 1, 0, 1, 1, 1],[0, 1, 1, 0, 0, 0]],
[[0, 1, 1, 1, 1, 1],[1, 0, 0, 0, 1, 0],[3, 1, 1, 1, 0, 1]],
[[0, 1, 0, 1, 0, 0],[1, 0, 0, 1, 0, 0],[3, 1, 0, 0, 0, 0]],
[[1, 0, 0, 1, 0, 1],[1, 0, 0, 1, 1, 0],[2, 0, 0, 0, 1, 1]],
[[1, 0, 1, 0, 0, 0],[1, 0, 1, 0, 0, 1],[2, 0, 0, 0, 0, 1]],
[[1, 0, 1, 0, 1, 1],[1, 0, 1, 1, 0, 0],[2, 0, 0, 1, 1, 1]],
[[1, 0, 1, 1, 1, 0],[1, 1, 0, 0, 0, 0],[2, 1, 1, 1, 1, 0]],
[[1, 1, 0, 0, 0, 1],[1, 1, 0, 0, 1, 1],[2, 0, 0, 0, 1, 0]],
[[1, 1, 0, 1, 0, 0],[1, 1, 1, 0, 0, 0],[2, 0, 1, 1, 0, 0]],
[[1, 1, 0, 1, 1, 1],[1, 1, 1, 0, 0, 1],[2, 0, 1, 1, 1, 0]],
[[1, 1, 1, 0, 1, 0],[1, 1, 1, 1, 0, 0],[2, 0, 0, 1, 1, 0]],
[[1, 0, 1, 0, 1, 0],[1, 1, 1, 1, 1, 0],[2, 1, 0, 1, 0, 0]],
[[1, 1, 0, 0, 1, 0],[1, 1, 1, 1, 0, 1],[2, 0, 1, 1, 1, 1]],
[[1, 0, 0, 1, 1, 1],[1, 1, 1, 0, 1, 1],[2, 1, 1, 1, 0, 0]],
[[2, 0, 0, 1, 0, 0],[3, 0, 0, 0, 0, 0],[3, 0, 0, 1, 0, 0]],
[[2, 0, 1, 0, 0, 1],[3, 0, 0, 0, 0, 1],[3, 0, 1, 0, 0, 0]],
[[2, 0, 1, 1, 0, 1],[3, 0, 0, 0, 1, 0],[3, 0, 1, 1, 1, 1]],
[[2, 1, 0, 0, 0, 0],[3, 0, 0, 1, 0, 1],[3, 1, 0, 1, 0, 1]],
[[2, 1, 0, 0, 1, 0],[3, 0, 0, 0, 1, 1],[3, 1, 0, 0, 0, 1]],
[[2, 1, 0, 1, 0, 1],[3, 0, 0, 1, 1, 1],[3, 1, 0, 0, 1, 0]],
[[2, 1, 1, 0, 0, 0],[3, 0, 0, 1, 1, 0],[3, 1, 1, 1, 1, 0]],
[[2, 1, 1, 0, 1, 1],[3, 0, 1, 1, 0, 1],[3, 1, 0, 1, 1, 0]],
[[1, 0, 1, 1, 1, 1],[1, 1, 0, 1, 1, 0],[2, 1, 1, 0, 0, 1]],
[[2, 1, 1, 1, 1, 1],[3, 0, 1, 0, 1, 1],[3, 1, 0, 1, 0, 0]],
[[2, 1, 0, 1, 1, 0],[3, 0, 1, 0, 0, 1],[3, 1, 1, 1, 1, 1]],
[[2, 0, 0, 1, 0, 1],[3, 1, 1, 0, 0, 1],[3, 1, 1, 1, 0, 0]],
[[3, 0, 1, 1, 0, 0],[3, 0, 1, 1, 1, 0],[3, 1, 1, 0, 0, 0],[3, 1, 1, 0, 1, 0]],
[[1, 1, 0, 1, 0, 1],[2, 1, 0, 0, 1, 1],[2, 1, 1, 1, 0, 1],[3, 1, 1, 0, 1, 1]],
[[0, 1, 0, 1, 0, 1],[0, 1, 1, 1, 0, 1],[2, 0, 0, 0, 0, 0],[2, 0, 1, 0, 0, 0]],
[[1, 0, 0, 0, 1, 1],[1, 1, 1, 1, 1, 1],[2, 1, 0, 0, 0, 1],[2, 1, 0, 1, 1, 1],[2, 1, 1, 0, 1, 0]],
[[0, 0, 0, 1, 1, 0],[1, 0, 1, 1, 0, 1],[2, 0, 1, 0, 1, 0],[2, 0, 1, 0, 1, 1],[3, 0, 1, 0, 1, 0]]
]
A partition for sets of sizes: 35*3  3*4  2*5
\end{verbatim}
\newpage

\begin{verbatim} 
[4, 2, 2, 2, 2, 2]
[36, 1, 3]
[
[[0, 0, 0, 0, 0, 1],[0, 0, 0, 0, 1, 0],[0, 0, 0, 0, 1, 1]],
[[0, 0, 0, 1, 0, 0],[0, 0, 1, 0, 0, 0],[0, 0, 1, 1, 0, 0]],
[[0, 0, 0, 1, 1, 1],[0, 0, 1, 0, 0, 1],[0, 0, 1, 1, 1, 0]],
[[0, 0, 1, 0, 1, 0],[0, 1, 0, 0, 0, 0],[0, 1, 1, 0, 1, 0]],
[[0, 0, 1, 1, 0, 1],[0, 1, 0, 0, 0, 1],[0, 1, 1, 1, 0, 0]],
[[0, 0, 1, 0, 1, 1],[0, 1, 0, 0, 1, 0],[0, 1, 1, 0, 0, 1]],
[[0, 1, 0, 0, 1, 1],[1, 0, 0, 0, 0, 0],[3, 1, 0, 0, 1, 1]],
[[0, 1, 0, 1, 1, 0],[1, 0, 0, 0, 0, 1],[3, 1, 0, 1, 1, 1]],
[[0, 0, 0, 1, 0, 1],[0, 1, 1, 0, 1, 1],[0, 1, 1, 1, 1, 0]],
[[0, 0, 1, 1, 1, 1],[0, 1, 0, 1, 1, 1],[0, 1, 1, 0, 0, 0]],
[[0, 1, 1, 1, 1, 1],[1, 0, 0, 0, 1, 0],[3, 1, 1, 1, 0, 1]],
[[0, 1, 0, 1, 0, 0],[1, 0, 0, 1, 0, 0],[3, 1, 0, 0, 0, 0]],
[[1, 0, 0, 1, 0, 1],[1, 0, 0, 1, 1, 0],[2, 0, 0, 0, 1, 1]],
[[1, 0, 1, 0, 0, 0],[1, 0, 1, 0, 0, 1],[2, 0, 0, 0, 0, 1]],
[[1, 0, 1, 0, 1, 1],[1, 0, 1, 1, 0, 0],[2, 0, 0, 1, 1, 1]],
[[1, 0, 1, 1, 1, 0],[1, 1, 0, 0, 0, 0],[2, 1, 1, 1, 1, 0]],
[[1, 1, 0, 0, 0, 1],[1, 1, 0, 0, 1, 1],[2, 0, 0, 0, 1, 0]],
[[1, 1, 0, 1, 0, 0],[1, 1, 1, 0, 0, 0],[2, 0, 1, 1, 0, 0]],
[[1, 1, 0, 1, 1, 1],[1, 1, 1, 0, 0, 1],[2, 0, 1, 1, 1, 0]],
[[1, 1, 1, 0, 1, 0],[1, 1, 1, 1, 0, 0],[2, 0, 0, 1, 1, 0]],
[[1, 0, 1, 0, 1, 0],[1, 1, 1, 1, 1, 0],[2, 1, 0, 1, 0, 0]],
[[1, 1, 0, 0, 1, 0],[1, 1, 1, 1, 0, 1],[2, 0, 1, 1, 1, 1]],
[[1, 0, 0, 1, 1, 1],[1, 1, 1, 0, 1, 1],[2, 1, 1, 1, 0, 0]],
[[2, 0, 0, 1, 0, 0],[3, 0, 0, 0, 0, 0],[3, 0, 0, 1, 0, 0]],
[[2, 0, 1, 0, 0, 1],[3, 0, 0, 0, 0, 1],[3, 0, 1, 0, 0, 0]],
[[2, 0, 1, 1, 0, 1],[3, 0, 0, 0, 1, 0],[3, 0, 1, 1, 1, 1]],
[[2, 1, 0, 0, 0, 0],[3, 0, 0, 1, 0, 1],[3, 1, 0, 1, 0, 1]],
[[2, 1, 0, 0, 1, 0],[3, 0, 0, 0, 1, 1],[3, 1, 0, 0, 0, 1]],
[[2, 1, 0, 1, 0, 1],[3, 0, 0, 1, 1, 1],[3, 1, 0, 0, 1, 0]],
[[2, 1, 1, 0, 0, 0],[3, 0, 0, 1, 1, 0],[3, 1, 1, 1, 1, 0]],
[[2, 1, 1, 0, 1, 1],[3, 0, 1, 1, 0, 1],[3, 1, 0, 1, 1, 0]],
[[1, 0, 1, 1, 1, 1],[1, 1, 0, 1, 1, 0],[2, 1, 1, 0, 0, 1]],
[[2, 1, 1, 1, 1, 1],[3, 0, 1, 0, 1, 1],[3, 1, 0, 1, 0, 0]],
[[2, 1, 0, 1, 1, 0],[3, 0, 1, 0, 0, 1],[3, 1, 1, 1, 1, 1]],
[[2, 0, 0, 1, 0, 1],[3, 1, 1, 0, 0, 1],[3, 1, 1, 1, 0, 0]],
[[1, 1, 0, 1, 0, 1],[1, 1, 1, 1, 1, 1],[2, 0, 1, 0, 1, 0]],
[[1, 0, 0, 0, 1, 1],[2, 0, 1, 0, 0, 0],[2, 1, 0, 0, 1, 1],[3, 1, 1, 0, 0, 0]],
[[2, 1, 0, 1, 1, 1],[2, 1, 1, 0, 1, 0],[2, 1, 1, 1, 0, 1],[3, 0, 1, 0, 1, 0],[3, 1, 1, 0, 1, 0]],
[[0, 0, 0, 1, 1, 0],[0, 1, 0, 1, 0, 1],[2, 1, 0, 0, 0, 1],[3, 0, 1, 1, 0, 0],[3, 0, 1, 1, 1, 0]],
[[0, 1, 1, 1, 0, 1],[1, 0, 1, 1, 0, 1],[2, 0, 0, 0, 0, 0],[2, 0, 1, 0, 1, 1],[3, 1, 1, 0, 1, 1]]
]
A partition for sets of sizes: 36*3  1*4  3*5
\end{verbatim}
\newpage

\begin{verbatim} 
[4, 2, 2, 2, 2, 2]
[33, 2, 4]
[
[[0, 0, 0, 0, 0, 1],[0, 0, 0, 0, 1, 0],[0, 0, 0, 0, 1, 1]],
[[0, 0, 0, 1, 0, 0],[0, 0, 1, 0, 0, 0],[0, 0, 1, 1, 0, 0]],
[[0, 0, 0, 1, 1, 1],[0, 0, 1, 0, 0, 1],[0, 0, 1, 1, 1, 0]],
[[0, 0, 1, 0, 1, 0],[0, 1, 0, 0, 0, 0],[0, 1, 1, 0, 1, 0]],
[[0, 0, 1, 1, 0, 1],[0, 1, 0, 0, 0, 1],[0, 1, 1, 1, 0, 0]],
[[0, 0, 1, 0, 1, 1],[0, 1, 0, 0, 1, 0],[0, 1, 1, 0, 0, 1]],
[[0, 1, 0, 0, 1, 1],[1, 0, 0, 0, 0, 0],[3, 1, 0, 0, 1, 1]],
[[0, 1, 0, 1, 1, 0],[1, 0, 0, 0, 0, 1],[3, 1, 0, 1, 1, 1]],
[[0, 0, 0, 1, 0, 1],[0, 1, 1, 0, 1, 1],[0, 1, 1, 1, 1, 0]],
[[0, 0, 1, 1, 1, 1],[0, 1, 0, 1, 1, 1],[0, 1, 1, 0, 0, 0]],
[[0, 1, 1, 1, 1, 1],[1, 0, 0, 0, 1, 0],[3, 1, 1, 1, 0, 1]],
[[0, 1, 0, 1, 0, 0],[1, 0, 0, 1, 0, 0],[3, 1, 0, 0, 0, 0]],
[[1, 0, 0, 1, 0, 1],[1, 0, 0, 1, 1, 0],[2, 0, 0, 0, 1, 1]],
[[1, 0, 1, 0, 0, 0],[1, 0, 1, 0, 0, 1],[2, 0, 0, 0, 0, 1]],
[[1, 0, 1, 0, 1, 1],[1, 0, 1, 1, 0, 0],[2, 0, 0, 1, 1, 1]],
[[1, 0, 1, 1, 1, 0],[1, 1, 0, 0, 0, 0],[2, 1, 1, 1, 1, 0]],
[[1, 1, 0, 0, 0, 1],[1, 1, 0, 0, 1, 1],[2, 0, 0, 0, 1, 0]],
[[1, 1, 0, 1, 0, 0],[1, 1, 1, 0, 0, 0],[2, 0, 1, 1, 0, 0]],
[[1, 1, 0, 1, 1, 1],[1, 1, 1, 0, 0, 1],[2, 0, 1, 1, 1, 0]],
[[1, 1, 1, 0, 1, 0],[1, 1, 1, 1, 0, 0],[2, 0, 0, 1, 1, 0]],
[[1, 0, 1, 0, 1, 0],[1, 1, 1, 1, 1, 0],[2, 1, 0, 1, 0, 0]],
[[1, 1, 0, 0, 1, 0],[1, 1, 1, 1, 0, 1],[2, 0, 1, 1, 1, 1]],
[[1, 0, 0, 1, 1, 1],[1, 1, 1, 0, 1, 1],[2, 1, 1, 1, 0, 0]],
[[2, 0, 0, 1, 0, 0],[3, 0, 0, 0, 0, 0],[3, 0, 0, 1, 0, 0]],
[[2, 0, 1, 0, 0, 1],[3, 0, 0, 0, 0, 1],[3, 0, 1, 0, 0, 0]],
[[2, 0, 1, 1, 0, 1],[3, 0, 0, 0, 1, 0],[3, 0, 1, 1, 1, 1]],
[[2, 1, 0, 0, 0, 0],[3, 0, 0, 1, 0, 1],[3, 1, 0, 1, 0, 1]],
[[2, 1, 0, 0, 1, 0],[3, 0, 0, 0, 1, 1],[3, 1, 0, 0, 0, 1]],
[[2, 1, 0, 1, 0, 1],[3, 0, 0, 1, 1, 1],[3, 1, 0, 0, 1, 0]],
[[2, 1, 1, 0, 0, 0],[3, 0, 0, 1, 1, 0],[3, 1, 1, 1, 1, 0]],
[[2, 1, 1, 0, 1, 1],[3, 0, 1, 1, 0, 1],[3, 1, 0, 1, 1, 0]],
[[1, 0, 1, 1, 1, 1],[1, 1, 0, 1, 1, 0],[2, 1, 1, 0, 0, 1]],
[[2, 1, 1, 1, 1, 1],[3, 0, 1, 0, 1, 1],[3, 1, 0, 1, 0, 0]],
[[2, 0, 1, 0, 0, 0],[2, 0, 1, 0, 1, 0],[2, 1, 0, 0, 0, 1],[2, 1, 0, 0, 1, 1]],
[[3, 0, 1, 0, 0, 1],[3, 0, 1, 0, 1, 0],[3, 1, 1, 0, 0, 0],[3, 1, 1, 0, 1, 1]],
[[2, 0, 0, 1, 0, 1],[2, 0, 1, 0, 1, 1],[2, 1, 1, 1, 0, 1],[3, 0, 1, 1, 0, 0],[3, 1, 1, 1, 1, 1]],
[[0, 1, 0, 1, 0, 1],[1, 0, 0, 0, 1, 1],[2, 1, 0, 1, 1, 0],[2, 1, 1, 0, 1, 0],[3, 1, 1, 0, 1, 0]],
[[0, 0, 0, 1, 1, 0],[1, 0, 1, 1, 0, 1],[2, 0, 0, 0, 0, 0],[2, 1, 0, 1, 1, 1],[3, 1, 1, 1, 0, 0]],
[[0, 1, 1, 1, 0, 1],[1, 1, 0, 1, 0, 1],[1, 1, 1, 1, 1, 1],[3, 0, 1, 1, 1, 0],[3, 1, 1, 0, 0, 1]]
]
A partition for sets of sizes: 33*3  2*4  4*5
\end{verbatim}
\newpage

\begin{verbatim} 
[4, 2, 2, 2, 2, 2]
[34, 0, 5]
[
[[0, 0, 0, 0, 0, 1],[0, 0, 0, 0, 1, 0],[0, 0, 0, 0, 1, 1]],
[[0, 0, 0, 1, 0, 0],[0, 0, 1, 0, 0, 0],[0, 0, 1, 1, 0, 0]],
[[0, 0, 0, 1, 1, 1],[0, 0, 1, 0, 0, 1],[0, 0, 1, 1, 1, 0]],
[[0, 0, 1, 0, 1, 0],[0, 1, 0, 0, 0, 0],[0, 1, 1, 0, 1, 0]],
[[0, 0, 1, 1, 0, 1],[0, 1, 0, 0, 0, 1],[0, 1, 1, 1, 0, 0]],
[[0, 0, 1, 0, 1, 1],[0, 1, 0, 0, 1, 0],[0, 1, 1, 0, 0, 1]],
[[0, 1, 0, 0, 1, 1],[1, 0, 0, 0, 0, 0],[3, 1, 0, 0, 1, 1]],
[[0, 1, 0, 1, 1, 0],[1, 0, 0, 0, 0, 1],[3, 1, 0, 1, 1, 1]],
[[0, 0, 0, 1, 0, 1],[0, 1, 1, 0, 1, 1],[0, 1, 1, 1, 1, 0]],
[[0, 0, 1, 1, 1, 1],[0, 1, 0, 1, 1, 1],[0, 1, 1, 0, 0, 0]],
[[0, 1, 1, 1, 1, 1],[1, 0, 0, 0, 1, 0],[3, 1, 1, 1, 0, 1]],
[[0, 1, 0, 1, 0, 0],[1, 0, 0, 1, 0, 0],[3, 1, 0, 0, 0, 0]],
[[1, 0, 0, 1, 0, 1],[1, 0, 0, 1, 1, 0],[2, 0, 0, 0, 1, 1]],
[[1, 0, 1, 0, 0, 0],[1, 0, 1, 0, 0, 1],[2, 0, 0, 0, 0, 1]],
[[1, 0, 1, 0, 1, 1],[1, 0, 1, 1, 0, 0],[2, 0, 0, 1, 1, 1]],
[[1, 0, 1, 1, 1, 0],[1, 1, 0, 0, 0, 0],[2, 1, 1, 1, 1, 0]],
[[1, 1, 0, 0, 0, 1],[1, 1, 0, 0, 1, 1],[2, 0, 0, 0, 1, 0]],
[[1, 1, 0, 1, 0, 0],[1, 1, 1, 0, 0, 0],[2, 0, 1, 1, 0, 0]],
[[1, 1, 0, 1, 1, 1],[1, 1, 1, 0, 0, 1],[2, 0, 1, 1, 1, 0]],
[[1, 1, 1, 0, 1, 0],[1, 1, 1, 1, 0, 0],[2, 0, 0, 1, 1, 0]],
[[1, 0, 1, 0, 1, 0],[1, 1, 1, 1, 1, 0],[2, 1, 0, 1, 0, 0]],
[[1, 1, 0, 0, 1, 0],[1, 1, 1, 1, 0, 1],[2, 0, 1, 1, 1, 1]],
[[1, 0, 0, 1, 1, 1],[1, 1, 1, 0, 1, 1],[2, 1, 1, 1, 0, 0]],
[[2, 0, 0, 1, 0, 0],[3, 0, 0, 0, 0, 0],[3, 0, 0, 1, 0, 0]],
[[2, 0, 1, 0, 0, 1],[3, 0, 0, 0, 0, 1],[3, 0, 1, 0, 0, 0]],
[[2, 0, 1, 1, 0, 1],[3, 0, 0, 0, 1, 0],[3, 0, 1, 1, 1, 1]],
[[2, 1, 0, 0, 0, 0],[3, 0, 0, 1, 0, 1],[3, 1, 0, 1, 0, 1]],
[[2, 1, 0, 0, 1, 0],[3, 0, 0, 0, 1, 1],[3, 1, 0, 0, 0, 1]],
[[2, 1, 0, 1, 0, 1],[3, 0, 0, 1, 1, 1],[3, 1, 0, 0, 1, 0]],
[[2, 1, 1, 0, 0, 0],[3, 0, 0, 1, 1, 0],[3, 1, 1, 1, 1, 0]],
[[2, 1, 1, 0, 1, 1],[3, 0, 1, 1, 0, 1],[3, 1, 0, 1, 1, 0]],
[[1, 0, 1, 1, 1, 1],[1, 1, 0, 1, 1, 0],[2, 1, 1, 0, 0, 1]],
[[2, 1, 1, 1, 1, 1],[3, 0, 1, 0, 1, 1],[3, 1, 0, 1, 0, 0]],
[[2, 1, 0, 1, 1, 0],[3, 0, 1, 0, 0, 1],[3, 1, 1, 1, 1, 1]],
[[2, 0, 1, 0, 0, 0],[2, 0, 1, 0, 1, 1],[2, 1, 0, 1, 1, 1],[3, 0, 1, 1, 0, 0],[3, 1, 1, 0, 0, 0]],
[[2, 0, 0, 1, 0, 1],[2, 1, 0, 0, 0, 1],[2, 1, 0, 0, 1, 1],[3, 1, 1, 0, 1, 1],[3, 1, 1, 1, 0, 0]],
[[2, 0, 0, 0, 0, 0],[2, 0, 1, 0, 1, 0],[2, 1, 1, 1, 0, 1],[3, 0, 1, 1, 1, 0],[3, 1, 1, 0, 0, 1]],
[[0, 1, 0, 1, 0, 1],[0, 1, 1, 1, 0, 1],[1, 0, 1, 1, 0, 1],[1, 1, 1, 1, 1, 1],[2, 1, 1, 0, 1, 0]],
[[0, 0, 0, 1, 1, 0],[1, 0, 0, 0, 1, 1],[1, 1, 0, 1, 0, 1],[3, 0, 1, 0, 1, 0],[3, 1, 1, 0, 1, 0]]
]
A partition for sets of sizes: 34*3  0*4  5*5
\end{verbatim}
\newpage

\begin{verbatim} 
[4, 2, 2, 2, 2, 2]
[30, 3, 5]
[
[[0, 0, 0, 0, 0, 1],[0, 0, 0, 0, 1, 0],[0, 0, 0, 0, 1, 1]],
[[0, 0, 0, 1, 0, 0],[0, 0, 1, 0, 0, 0],[0, 0, 1, 1, 0, 0]],
[[0, 0, 0, 1, 1, 1],[0, 0, 1, 0, 0, 1],[0, 0, 1, 1, 1, 0]],
[[0, 0, 1, 0, 1, 0],[0, 1, 0, 0, 0, 0],[0, 1, 1, 0, 1, 0]],
[[0, 0, 1, 1, 0, 1],[0, 1, 0, 0, 0, 1],[0, 1, 1, 1, 0, 0]],
[[0, 0, 1, 0, 1, 1],[0, 1, 0, 0, 1, 0],[0, 1, 1, 0, 0, 1]],
[[0, 1, 0, 0, 1, 1],[1, 0, 0, 0, 0, 0],[3, 1, 0, 0, 1, 1]],
[[0, 1, 0, 1, 1, 0],[1, 0, 0, 0, 0, 1],[3, 1, 0, 1, 1, 1]],
[[0, 0, 0, 1, 0, 1],[0, 1, 1, 0, 1, 1],[0, 1, 1, 1, 1, 0]],
[[0, 0, 1, 1, 1, 1],[0, 1, 0, 1, 1, 1],[0, 1, 1, 0, 0, 0]],
[[0, 1, 1, 1, 1, 1],[1, 0, 0, 0, 1, 0],[3, 1, 1, 1, 0, 1]],
[[0, 1, 0, 1, 0, 0],[1, 0, 0, 1, 0, 0],[3, 1, 0, 0, 0, 0]],
[[1, 0, 0, 1, 0, 1],[1, 0, 0, 1, 1, 0],[2, 0, 0, 0, 1, 1]],
[[1, 0, 1, 0, 0, 0],[1, 0, 1, 0, 0, 1],[2, 0, 0, 0, 0, 1]],
[[1, 0, 1, 0, 1, 1],[1, 0, 1, 1, 0, 0],[2, 0, 0, 1, 1, 1]],
[[1, 0, 1, 1, 1, 0],[1, 1, 0, 0, 0, 0],[2, 1, 1, 1, 1, 0]],
[[1, 1, 0, 0, 0, 1],[1, 1, 0, 0, 1, 1],[2, 0, 0, 0, 1, 0]],
[[1, 1, 0, 1, 0, 0],[1, 1, 1, 0, 0, 0],[2, 0, 1, 1, 0, 0]],
[[1, 1, 0, 1, 1, 1],[1, 1, 1, 0, 0, 1],[2, 0, 1, 1, 1, 0]],
[[1, 1, 1, 0, 1, 0],[1, 1, 1, 1, 0, 0],[2, 0, 0, 1, 1, 0]],
[[1, 0, 1, 0, 1, 0],[1, 1, 1, 1, 1, 0],[2, 1, 0, 1, 0, 0]],
[[1, 1, 0, 0, 1, 0],[1, 1, 1, 1, 0, 1],[2, 0, 1, 1, 1, 1]],
[[1, 0, 0, 1, 1, 1],[1, 1, 1, 0, 1, 1],[2, 1, 1, 1, 0, 0]],
[[2, 0, 0, 1, 0, 0],[3, 0, 0, 0, 0, 0],[3, 0, 0, 1, 0, 0]],
[[2, 0, 1, 0, 0, 1],[3, 0, 0, 0, 0, 1],[3, 0, 1, 0, 0, 0]],
[[2, 0, 1, 1, 0, 1],[3, 0, 0, 0, 1, 0],[3, 0, 1, 1, 1, 1]],
[[2, 1, 0, 0, 0, 0],[3, 0, 0, 1, 0, 1],[3, 1, 0, 1, 0, 1]],
[[2, 1, 0, 0, 1, 0],[3, 0, 0, 0, 1, 1],[3, 1, 0, 0, 0, 1]],
[[2, 1, 0, 1, 0, 1],[3, 0, 0, 1, 1, 1],[3, 1, 0, 0, 1, 0]],
[[2, 1, 1, 0, 0, 0],[3, 0, 0, 1, 1, 0],[3, 1, 1, 1, 1, 0]],
[[2, 0, 0, 1, 0, 1],[2, 0, 1, 0, 1, 0],[2, 1, 0, 1, 1, 0],[2, 1, 1, 0, 0, 1]],
[[1, 0, 1, 1, 0, 1],[2, 1, 1, 0, 1, 1],[2, 1, 1, 1, 1, 1],[3, 0, 1, 0, 0, 1]],
[[2, 0, 0, 0, 0, 0],[2, 0, 1, 0, 1, 1],[2, 1, 0, 0, 0, 1],[2, 1, 1, 0, 1, 0]],
[[2, 0, 1, 0, 0, 0],[2, 1, 0, 1, 1, 1],[2, 1, 1, 1, 0, 1],[3, 0, 1, 1, 0, 0],[3, 0, 1, 1, 1, 0]],
[[1, 1, 0, 1, 0, 1],[2, 1, 0, 0, 1, 1],[3, 0, 1, 0, 1, 1],[3, 1, 0, 1, 0, 0],[3, 1, 1, 0, 0, 1]],
[[0, 1, 0, 1, 0, 1],[1, 1, 0, 1, 1, 0],[1, 1, 1, 1, 1, 1],[3, 0, 1, 0, 1, 0],[3, 1, 0, 1, 1, 0]],
[[0, 0, 0, 1, 1, 0],[3, 1, 1, 0, 0, 0],[3, 1, 1, 0, 1, 0],[3, 1, 1, 0, 1, 1],[3, 1, 1, 1, 1, 1]],
[[0, 1, 1, 1, 0, 1],[1, 0, 0, 0, 1, 1],[1, 0, 1, 1, 1, 1],[3, 0, 1, 1, 0, 1],[3, 1, 1, 1, 0, 0]]
]
A partition for sets of sizes: 30*3  3*4  5*5
\end{verbatim}
\newpage

\begin{verbatim} 
[4, 2, 2, 2, 2, 2]
[31, 1, 6]
[
[[0, 0, 0, 0, 0, 1],[0, 0, 0, 0, 1, 0],[0, 0, 0, 0, 1, 1]],
[[0, 0, 0, 1, 0, 0],[0, 0, 1, 0, 0, 0],[0, 0, 1, 1, 0, 0]],
[[0, 0, 0, 1, 1, 1],[0, 0, 1, 0, 0, 1],[0, 0, 1, 1, 1, 0]],
[[0, 0, 1, 0, 1, 0],[0, 1, 0, 0, 0, 0],[0, 1, 1, 0, 1, 0]],
[[0, 0, 1, 1, 0, 1],[0, 1, 0, 0, 0, 1],[0, 1, 1, 1, 0, 0]],
[[0, 0, 1, 0, 1, 1],[0, 1, 0, 0, 1, 0],[0, 1, 1, 0, 0, 1]],
[[0, 1, 0, 0, 1, 1],[1, 0, 0, 0, 0, 0],[3, 1, 0, 0, 1, 1]],
[[0, 1, 0, 1, 1, 0],[1, 0, 0, 0, 0, 1],[3, 1, 0, 1, 1, 1]],
[[0, 0, 0, 1, 0, 1],[0, 1, 1, 0, 1, 1],[0, 1, 1, 1, 1, 0]],
[[0, 0, 1, 1, 1, 1],[0, 1, 0, 1, 1, 1],[0, 1, 1, 0, 0, 0]],
[[0, 1, 1, 1, 1, 1],[1, 0, 0, 0, 1, 0],[3, 1, 1, 1, 0, 1]],
[[0, 1, 0, 1, 0, 0],[1, 0, 0, 1, 0, 0],[3, 1, 0, 0, 0, 0]],
[[1, 0, 0, 1, 0, 1],[1, 0, 0, 1, 1, 0],[2, 0, 0, 0, 1, 1]],
[[1, 0, 1, 0, 0, 0],[1, 0, 1, 0, 0, 1],[2, 0, 0, 0, 0, 1]],
[[1, 0, 1, 0, 1, 1],[1, 0, 1, 1, 0, 0],[2, 0, 0, 1, 1, 1]],
[[1, 0, 1, 1, 1, 0],[1, 1, 0, 0, 0, 0],[2, 1, 1, 1, 1, 0]],
[[1, 1, 0, 0, 0, 1],[1, 1, 0, 0, 1, 1],[2, 0, 0, 0, 1, 0]],
[[1, 1, 0, 1, 0, 0],[1, 1, 1, 0, 0, 0],[2, 0, 1, 1, 0, 0]],
[[1, 1, 0, 1, 1, 1],[1, 1, 1, 0, 0, 1],[2, 0, 1, 1, 1, 0]],
[[1, 1, 1, 0, 1, 0],[1, 1, 1, 1, 0, 0],[2, 0, 0, 1, 1, 0]],
[[1, 0, 1, 0, 1, 0],[1, 1, 1, 1, 1, 0],[2, 1, 0, 1, 0, 0]],
[[1, 1, 0, 0, 1, 0],[1, 1, 1, 1, 0, 1],[2, 0, 1, 1, 1, 1]],
[[1, 0, 0, 1, 1, 1],[1, 1, 1, 0, 1, 1],[2, 1, 1, 1, 0, 0]],
[[2, 0, 0, 1, 0, 0],[3, 0, 0, 0, 0, 0],[3, 0, 0, 1, 0, 0]],
[[2, 0, 1, 0, 0, 1],[3, 0, 0, 0, 0, 1],[3, 0, 1, 0, 0, 0]],
[[2, 0, 1, 1, 0, 1],[3, 0, 0, 0, 1, 0],[3, 0, 1, 1, 1, 1]],
[[2, 1, 0, 0, 0, 0],[3, 0, 0, 1, 0, 1],[3, 1, 0, 1, 0, 1]],
[[2, 1, 0, 0, 1, 0],[3, 0, 0, 0, 1, 1],[3, 1, 0, 0, 0, 1]],
[[2, 1, 0, 1, 0, 1],[3, 0, 0, 1, 1, 1],[3, 1, 0, 0, 1, 0]],
[[2, 1, 1, 0, 0, 0],[3, 0, 0, 1, 1, 0],[3, 1, 1, 1, 1, 0]],
[[2, 1, 1, 0, 1, 1],[3, 0, 1, 1, 0, 1],[3, 1, 0, 1, 1, 0]],
[[1, 0, 1, 1, 1, 1],[2, 0, 1, 0, 1, 0],[2, 1, 1, 1, 1, 1],[3, 1, 1, 0, 1, 0]],
[[1, 1, 0, 1, 1, 0],[2, 1, 0, 0, 1, 1],[3, 0, 1, 0, 0, 1],[3, 1, 0, 1, 0, 0],[3, 1, 1, 0, 0, 0]],
[[2, 1, 0, 1, 1, 0],[2, 1, 1, 0, 0, 1],[2, 1, 1, 1, 0, 1],[3, 0, 1, 0, 1, 1],[3, 1, 1, 0, 0, 1]],
[[2, 0, 0, 1, 0, 1],[2, 1, 0, 0, 0, 1],[2, 1, 0, 1, 1, 1],[3, 1, 1, 1, 0, 0],[3, 1, 1, 1, 1, 1]],
[[1, 1, 1, 1, 1, 1],[2, 0, 0, 0, 0, 0],[3, 0, 1, 0, 1, 0],[3, 0, 1, 1, 1, 0],[3, 1, 1, 0, 1, 1]],
[[0, 0, 0, 1, 1, 0],[0, 1, 0, 1, 0, 1],[1, 0, 1, 1, 0, 1],[1, 1, 0, 1, 0, 1],[2, 0, 1, 0, 1, 1]],
[[0, 1, 1, 1, 0, 1],[1, 0, 0, 0, 1, 1],[2, 0, 1, 0, 0, 0],[2, 1, 1, 0, 1, 0],[3, 0, 1, 1, 0, 0]]
]
A partition for sets of sizes: 31*3  1*4  6*5
\end{verbatim}
\newpage

\begin{verbatim} 
[4, 2, 2, 2, 2, 2]
[28, 2, 7]
[
[[0, 0, 0, 0, 0, 1],[0, 0, 0, 0, 1, 0],[0, 0, 0, 0, 1, 1]],
[[0, 0, 0, 1, 0, 0],[0, 0, 1, 0, 0, 0],[0, 0, 1, 1, 0, 0]],
[[0, 0, 0, 1, 1, 1],[0, 0, 1, 0, 0, 1],[0, 0, 1, 1, 1, 0]],
[[0, 0, 1, 0, 1, 0],[0, 1, 0, 0, 0, 0],[0, 1, 1, 0, 1, 0]],
[[0, 0, 1, 1, 0, 1],[0, 1, 0, 0, 0, 1],[0, 1, 1, 1, 0, 0]],
[[0, 0, 1, 0, 1, 1],[0, 1, 0, 0, 1, 0],[0, 1, 1, 0, 0, 1]],
[[0, 1, 0, 0, 1, 1],[1, 0, 0, 0, 0, 0],[3, 1, 0, 0, 1, 1]],
[[0, 1, 0, 1, 1, 0],[1, 0, 0, 0, 0, 1],[3, 1, 0, 1, 1, 1]],
[[0, 0, 0, 1, 0, 1],[0, 1, 1, 0, 1, 1],[0, 1, 1, 1, 1, 0]],
[[0, 0, 1, 1, 1, 1],[0, 1, 0, 1, 1, 1],[0, 1, 1, 0, 0, 0]],
[[0, 1, 1, 1, 1, 1],[1, 0, 0, 0, 1, 0],[3, 1, 1, 1, 0, 1]],
[[0, 1, 0, 1, 0, 0],[1, 0, 0, 1, 0, 0],[3, 1, 0, 0, 0, 0]],
[[1, 0, 0, 1, 0, 1],[1, 0, 0, 1, 1, 0],[2, 0, 0, 0, 1, 1]],
[[1, 0, 1, 0, 0, 0],[1, 0, 1, 0, 0, 1],[2, 0, 0, 0, 0, 1]],
[[1, 0, 1, 0, 1, 1],[1, 0, 1, 1, 0, 0],[2, 0, 0, 1, 1, 1]],
[[1, 0, 1, 1, 1, 0],[1, 1, 0, 0, 0, 0],[2, 1, 1, 1, 1, 0]],
[[1, 1, 0, 0, 0, 1],[1, 1, 0, 0, 1, 1],[2, 0, 0, 0, 1, 0]],
[[1, 1, 0, 1, 0, 0],[1, 1, 1, 0, 0, 0],[2, 0, 1, 1, 0, 0]],
[[1, 1, 0, 1, 1, 1],[1, 1, 1, 0, 0, 1],[2, 0, 1, 1, 1, 0]],
[[1, 1, 1, 0, 1, 0],[1, 1, 1, 1, 0, 0],[2, 0, 0, 1, 1, 0]],
[[1, 0, 1, 0, 1, 0],[1, 1, 1, 1, 1, 0],[2, 1, 0, 1, 0, 0]],
[[1, 1, 0, 0, 1, 0],[1, 1, 1, 1, 0, 1],[2, 0, 1, 1, 1, 1]],
[[1, 0, 0, 1, 1, 1],[1, 1, 1, 0, 1, 1],[2, 1, 1, 1, 0, 0]],
[[2, 0, 0, 1, 0, 0],[3, 0, 0, 0, 0, 0],[3, 0, 0, 1, 0, 0]],
[[2, 0, 1, 0, 0, 1],[3, 0, 0, 0, 0, 1],[3, 0, 1, 0, 0, 0]],
[[2, 0, 1, 1, 0, 1],[3, 0, 0, 0, 1, 0],[3, 0, 1, 1, 1, 1]],
[[2, 1, 0, 0, 0, 0],[3, 0, 0, 1, 0, 1],[3, 1, 0, 1, 0, 1]],
[[2, 1, 0, 0, 1, 0],[3, 0, 0, 0, 1, 1],[3, 1, 0, 0, 0, 1]],
[[2, 1, 0, 1, 0, 1],[2, 1, 0, 1, 1, 0],[2, 1, 1, 0, 0, 0],[2, 1, 1, 0, 1, 1]],
[[2, 0, 0, 1, 0, 1],[2, 0, 1, 0, 0, 0],[2, 1, 0, 1, 1, 1],[2, 1, 1, 0, 1, 0]],
[[1, 0, 1, 1, 1, 1],[2, 1, 1, 1, 0, 1],[3, 0, 0, 1, 1, 0],[3, 0, 1, 0, 1, 0],[3, 1, 1, 1, 1, 0]],
[[1, 1, 0, 1, 1, 0],[2, 1, 0, 0, 1, 1],[3, 0, 0, 1, 1, 1],[3, 0, 1, 0, 0, 1],[3, 0, 1, 0, 1, 1]],
[[2, 0, 1, 0, 1, 1],[2, 1, 1, 0, 0, 1],[2, 1, 1, 1, 1, 1],[3, 1, 0, 0, 1, 0],[3, 1, 1, 1, 1, 1]],
[[1, 1, 0, 1, 0, 1],[2, 1, 0, 0, 0, 1],[3, 0, 1, 1, 1, 0],[3, 1, 0, 1, 1, 0],[3, 1, 1, 1, 0, 0]],
[[1, 1, 1, 1, 1, 1],[2, 0, 1, 0, 1, 0],[3, 1, 0, 1, 0, 0],[3, 1, 1, 0, 0, 0],[3, 1, 1, 0, 0, 1]],
[[0, 0, 0, 1, 1, 0],[0, 1, 0, 1, 0, 1],[0, 1, 1, 1, 0, 1],[1, 0, 0, 0, 1, 1],[3, 0, 1, 1, 0, 1]],
[[1, 0, 1, 1, 0, 1],[2, 0, 0, 0, 0, 0],[3, 0, 1, 1, 0, 0],[3, 1, 1, 0, 1, 0],[3, 1, 1, 0, 1, 1]]
]
A partition for sets of sizes: 28*3  2*4  7*5
\end{verbatim}
\newpage

\begin{verbatim} 
[4, 2, 2, 2, 2, 2]
[29, 0, 8]
[
[[0, 0, 0, 0, 0, 1],[0, 0, 0, 0, 1, 0],[0, 0, 0, 0, 1, 1]],
[[0, 0, 0, 1, 0, 0],[0, 0, 1, 0, 0, 0],[0, 0, 1, 1, 0, 0]],
[[0, 0, 0, 1, 1, 1],[0, 0, 1, 0, 0, 1],[0, 0, 1, 1, 1, 0]],
[[0, 0, 1, 0, 1, 0],[0, 1, 0, 0, 0, 0],[0, 1, 1, 0, 1, 0]],
[[0, 0, 1, 1, 0, 1],[0, 1, 0, 0, 0, 1],[0, 1, 1, 1, 0, 0]],
[[0, 0, 1, 0, 1, 1],[0, 1, 0, 0, 1, 0],[0, 1, 1, 0, 0, 1]],
[[0, 1, 0, 0, 1, 1],[1, 0, 0, 0, 0, 0],[3, 1, 0, 0, 1, 1]],
[[0, 1, 0, 1, 1, 0],[1, 0, 0, 0, 0, 1],[3, 1, 0, 1, 1, 1]],
[[0, 0, 0, 1, 0, 1],[0, 1, 1, 0, 1, 1],[0, 1, 1, 1, 1, 0]],
[[0, 0, 1, 1, 1, 1],[0, 1, 0, 1, 1, 1],[0, 1, 1, 0, 0, 0]],
[[0, 1, 1, 1, 1, 1],[1, 0, 0, 0, 1, 0],[3, 1, 1, 1, 0, 1]],
[[0, 1, 0, 1, 0, 0],[1, 0, 0, 1, 0, 0],[3, 1, 0, 0, 0, 0]],
[[1, 0, 0, 1, 0, 1],[1, 0, 0, 1, 1, 0],[2, 0, 0, 0, 1, 1]],
[[1, 0, 1, 0, 0, 0],[1, 0, 1, 0, 0, 1],[2, 0, 0, 0, 0, 1]],
[[1, 0, 1, 0, 1, 1],[1, 0, 1, 1, 0, 0],[2, 0, 0, 1, 1, 1]],
[[1, 0, 1, 1, 1, 0],[1, 1, 0, 0, 0, 0],[2, 1, 1, 1, 1, 0]],
[[1, 1, 0, 0, 0, 1],[1, 1, 0, 0, 1, 1],[2, 0, 0, 0, 1, 0]],
[[1, 1, 0, 1, 0, 0],[1, 1, 1, 0, 0, 0],[2, 0, 1, 1, 0, 0]],
[[1, 1, 0, 1, 1, 1],[1, 1, 1, 0, 0, 1],[2, 0, 1, 1, 1, 0]],
[[1, 1, 1, 0, 1, 0],[1, 1, 1, 1, 0, 0],[2, 0, 0, 1, 1, 0]],
[[1, 0, 1, 0, 1, 0],[1, 1, 1, 1, 1, 0],[2, 1, 0, 1, 0, 0]],
[[1, 1, 0, 0, 1, 0],[1, 1, 1, 1, 0, 1],[2, 0, 1, 1, 1, 1]],
[[1, 0, 0, 1, 1, 1],[1, 1, 1, 0, 1, 1],[2, 1, 1, 1, 0, 0]],
[[2, 0, 0, 1, 0, 0],[3, 0, 0, 0, 0, 0],[3, 0, 0, 1, 0, 0]],
[[2, 0, 1, 0, 0, 1],[3, 0, 0, 0, 0, 1],[3, 0, 1, 0, 0, 0]],
[[2, 0, 1, 1, 0, 1],[3, 0, 0, 0, 1, 0],[3, 0, 1, 1, 1, 1]],
[[2, 1, 0, 0, 0, 0],[3, 0, 0, 1, 0, 1],[3, 1, 0, 1, 0, 1]],
[[2, 1, 0, 0, 1, 0],[3, 0, 0, 0, 1, 1],[3, 1, 0, 0, 0, 1]],
[[2, 1, 0, 1, 0, 1],[3, 0, 0, 1, 1, 1],[3, 1, 0, 0, 1, 0]],
[[2, 1, 1, 0, 0, 0],[2, 1, 1, 0, 0, 1],[2, 1, 1, 0, 1, 0],[3, 0, 1, 1, 0, 1],[3, 1, 0, 1, 1, 0]],
[[1, 0, 1, 1, 1, 1],[2, 1, 1, 1, 0, 1],[3, 0, 0, 1, 1, 0],[3, 0, 1, 0, 1, 0],[3, 1, 1, 1, 1, 0]],
[[1, 1, 0, 1, 1, 0],[2, 1, 0, 0, 1, 1],[3, 0, 1, 0, 0, 1],[3, 1, 0, 1, 0, 0],[3, 1, 1, 0, 0, 0]],
[[2, 1, 0, 1, 1, 0],[2, 1, 1, 0, 1, 1],[2, 1, 1, 1, 1, 1],[3, 0, 1, 0, 1, 1],[3, 1, 1, 0, 0, 1]],
[[1, 1, 0, 1, 0, 1],[2, 0, 1, 0, 1, 1],[3, 0, 1, 1, 0, 0],[3, 0, 1, 1, 1, 0],[3, 1, 1, 1, 0, 0]],
[[0, 1, 0, 1, 0, 1],[1, 1, 1, 1, 1, 1],[2, 0, 0, 0, 0, 0],[2, 1, 0, 0, 0, 1],[3, 1, 1, 0, 1, 1]],
[[0, 0, 0, 1, 1, 0],[1, 0, 0, 0, 1, 1],[2, 0, 1, 0, 0, 0],[2, 1, 0, 1, 1, 1],[3, 1, 1, 0, 1, 0]],
[[0, 1, 1, 1, 0, 1],[1, 0, 1, 1, 0, 1],[2, 0, 0, 1, 0, 1],[2, 0, 1, 0, 1, 0],[3, 1, 1, 1, 1, 1]]
]
A partition for sets of sizes: 29*3  0*4  8*5
\end{verbatim}
\newpage

\begin{verbatim} 
[4, 2, 2, 2, 2, 2]
[25, 3, 8]
[
[[0, 0, 0, 0, 0, 1],[0, 0, 0, 0, 1, 0],[0, 0, 0, 0, 1, 1]],
[[0, 0, 0, 1, 0, 0],[0, 0, 1, 0, 0, 0],[0, 0, 1, 1, 0, 0]],
[[0, 0, 0, 1, 1, 1],[0, 0, 1, 0, 0, 1],[0, 0, 1, 1, 1, 0]],
[[0, 0, 1, 0, 1, 0],[0, 1, 0, 0, 0, 0],[0, 1, 1, 0, 1, 0]],
[[0, 0, 1, 1, 0, 1],[0, 1, 0, 0, 0, 1],[0, 1, 1, 1, 0, 0]],
[[0, 0, 1, 0, 1, 1],[0, 1, 0, 0, 1, 0],[0, 1, 1, 0, 0, 1]],
[[0, 1, 0, 0, 1, 1],[1, 0, 0, 0, 0, 0],[3, 1, 0, 0, 1, 1]],
[[0, 1, 0, 1, 1, 0],[1, 0, 0, 0, 0, 1],[3, 1, 0, 1, 1, 1]],
[[0, 0, 0, 1, 0, 1],[0, 1, 1, 0, 1, 1],[0, 1, 1, 1, 1, 0]],
[[0, 0, 1, 1, 1, 1],[0, 1, 0, 1, 1, 1],[0, 1, 1, 0, 0, 0]],
[[0, 1, 1, 1, 1, 1],[1, 0, 0, 0, 1, 0],[3, 1, 1, 1, 0, 1]],
[[0, 1, 0, 1, 0, 0],[1, 0, 0, 1, 0, 0],[3, 1, 0, 0, 0, 0]],
[[1, 0, 0, 1, 0, 1],[1, 0, 0, 1, 1, 0],[2, 0, 0, 0, 1, 1]],
[[1, 0, 1, 0, 0, 0],[1, 0, 1, 0, 0, 1],[2, 0, 0, 0, 0, 1]],
[[1, 0, 1, 0, 1, 1],[1, 0, 1, 1, 0, 0],[2, 0, 0, 1, 1, 1]],
[[1, 0, 1, 1, 1, 0],[1, 1, 0, 0, 0, 0],[2, 1, 1, 1, 1, 0]],
[[1, 1, 0, 0, 0, 1],[1, 1, 0, 0, 1, 1],[2, 0, 0, 0, 1, 0]],
[[1, 1, 0, 1, 0, 0],[1, 1, 1, 0, 0, 0],[2, 0, 1, 1, 0, 0]],
[[1, 1, 0, 1, 1, 1],[1, 1, 1, 0, 0, 1],[2, 0, 1, 1, 1, 0]],
[[1, 1, 1, 0, 1, 0],[1, 1, 1, 1, 0, 0],[2, 0, 0, 1, 1, 0]],
[[1, 0, 1, 0, 1, 0],[1, 1, 1, 1, 1, 0],[2, 1, 0, 1, 0, 0]],
[[1, 1, 0, 0, 1, 0],[1, 1, 1, 1, 0, 1],[2, 0, 1, 1, 1, 1]],
[[1, 0, 0, 1, 1, 1],[1, 1, 1, 0, 1, 1],[2, 1, 1, 1, 0, 0]],
[[2, 0, 0, 1, 0, 0],[3, 0, 0, 0, 0, 0],[3, 0, 0, 1, 0, 0]],
[[2, 0, 1, 0, 0, 1],[3, 0, 0, 0, 0, 1],[3, 0, 1, 0, 0, 0]],
[[1, 1, 0, 1, 1, 0],[2, 0, 1, 1, 0, 1],[2, 1, 0, 0, 0, 0],[3, 0, 1, 0, 1, 1]],
[[1, 1, 0, 1, 0, 1],[2, 1, 0, 0, 0, 1],[2, 1, 0, 0, 1, 0],[3, 1, 0, 1, 1, 0]],
[[1, 1, 1, 1, 1, 1],[2, 1, 0, 1, 0, 1],[2, 1, 0, 1, 1, 0],[3, 1, 1, 1, 0, 0]],
[[2, 1, 1, 0, 0, 0],[2, 1, 1, 0, 0, 1],[2, 1, 1, 0, 1, 0],[3, 0, 0, 0, 1, 0],[3, 1, 1, 0, 0, 1]],
[[1, 0, 1, 1, 1, 1],[2, 1, 1, 1, 0, 1],[3, 0, 0, 0, 1, 1],[3, 0, 0, 1, 0, 1],[3, 1, 0, 1, 0, 0]],
[[2, 0, 1, 0, 0, 0],[2, 1, 1, 0, 1, 1],[2, 1, 1, 1, 1, 1],[3, 0, 0, 1, 1, 0],[3, 0, 1, 0, 1, 0]],
[[2, 0, 1, 0, 1, 1],[2, 1, 0, 0, 1, 1],[2, 1, 0, 1, 1, 1],[3, 1, 0, 0, 0, 1],[3, 1, 1, 1, 1, 0]],
[[1, 0, 0, 0, 1, 1],[2, 0, 1, 0, 1, 0],[3, 0, 1, 1, 1, 0],[3, 1, 0, 0, 1, 0],[3, 1, 0, 1, 0, 1]],
[[1, 0, 1, 1, 0, 1],[2, 0, 0, 1, 0, 1],[3, 0, 1, 1, 0, 1],[3, 1, 1, 0, 1, 0],[3, 1, 1, 1, 1, 1]],
[[0, 0, 0, 1, 1, 0],[3, 0, 1, 0, 0, 1],[3, 0, 1, 1, 0, 0],[3, 1, 1, 0, 0, 0],[3, 1, 1, 0, 1, 1]],
[[0, 1, 0, 1, 0, 1],[0, 1, 1, 1, 0, 1],[2, 0, 0, 0, 0, 0],[3, 0, 0, 1, 1, 1],[3, 0, 1, 1, 1, 1]]
]
A partition for sets of sizes: 25*3  3*4  8*5
\end{verbatim}
\newpage

\begin{verbatim} 
[4, 2, 2, 2, 2, 2]
[26, 1, 9]
[
[[0, 0, 0, 0, 0, 1],[0, 0, 0, 0, 1, 0],[0, 0, 0, 0, 1, 1]],
[[0, 0, 0, 1, 0, 0],[0, 0, 1, 0, 0, 0],[0, 0, 1, 1, 0, 0]],
[[0, 0, 0, 1, 1, 1],[0, 0, 1, 0, 0, 1],[0, 0, 1, 1, 1, 0]],
[[0, 0, 1, 0, 1, 0],[0, 1, 0, 0, 0, 0],[0, 1, 1, 0, 1, 0]],
[[0, 0, 1, 1, 0, 1],[0, 1, 0, 0, 0, 1],[0, 1, 1, 1, 0, 0]],
[[0, 0, 1, 0, 1, 1],[0, 1, 0, 0, 1, 0],[0, 1, 1, 0, 0, 1]],
[[0, 1, 0, 0, 1, 1],[1, 0, 0, 0, 0, 0],[3, 1, 0, 0, 1, 1]],
[[0, 1, 0, 1, 1, 0],[1, 0, 0, 0, 0, 1],[3, 1, 0, 1, 1, 1]],
[[0, 0, 0, 1, 0, 1],[0, 1, 1, 0, 1, 1],[0, 1, 1, 1, 1, 0]],
[[0, 0, 1, 1, 1, 1],[0, 1, 0, 1, 1, 1],[0, 1, 1, 0, 0, 0]],
[[0, 1, 1, 1, 1, 1],[1, 0, 0, 0, 1, 0],[3, 1, 1, 1, 0, 1]],
[[0, 1, 0, 1, 0, 0],[1, 0, 0, 1, 0, 0],[3, 1, 0, 0, 0, 0]],
[[1, 0, 0, 1, 0, 1],[1, 0, 0, 1, 1, 0],[2, 0, 0, 0, 1, 1]],
[[1, 0, 1, 0, 0, 0],[1, 0, 1, 0, 0, 1],[2, 0, 0, 0, 0, 1]],
[[1, 0, 1, 0, 1, 1],[1, 0, 1, 1, 0, 0],[2, 0, 0, 1, 1, 1]],
[[1, 0, 1, 1, 1, 0],[1, 1, 0, 0, 0, 0],[2, 1, 1, 1, 1, 0]],
[[1, 1, 0, 0, 0, 1],[1, 1, 0, 0, 1, 1],[2, 0, 0, 0, 1, 0]],
[[1, 1, 0, 1, 0, 0],[1, 1, 1, 0, 0, 0],[2, 0, 1, 1, 0, 0]],
[[1, 1, 0, 1, 1, 1],[1, 1, 1, 0, 0, 1],[2, 0, 1, 1, 1, 0]],
[[1, 1, 1, 0, 1, 0],[1, 1, 1, 1, 0, 0],[2, 0, 0, 1, 1, 0]],
[[1, 0, 1, 0, 1, 0],[1, 1, 1, 1, 1, 0],[2, 1, 0, 1, 0, 0]],
[[1, 1, 0, 0, 1, 0],[1, 1, 1, 1, 0, 1],[2, 0, 1, 1, 1, 1]],
[[1, 0, 0, 1, 1, 1],[1, 1, 1, 0, 1, 1],[2, 1, 1, 1, 0, 0]],
[[2, 0, 0, 1, 0, 0],[3, 0, 0, 0, 0, 0],[3, 0, 0, 1, 0, 0]],
[[2, 0, 1, 0, 0, 1],[3, 0, 0, 0, 0, 1],[3, 0, 1, 0, 0, 0]],
[[2, 0, 1, 1, 0, 1],[3, 0, 0, 0, 1, 0],[3, 0, 1, 1, 1, 1]],
[[2, 0, 0, 0, 0, 0],[2, 0, 0, 1, 0, 1],[2, 1, 0, 0, 1, 0],[2, 1, 0, 1, 1, 1]],
[[2, 1, 0, 0, 1, 1],[2, 1, 0, 1, 0, 1],[2, 1, 0, 1, 1, 0],[3, 0, 0, 1, 0, 1],[3, 1, 0, 1, 0, 1]],
[[2, 1, 1, 0, 0, 0],[2, 1, 1, 0, 0, 1],[2, 1, 1, 0, 1, 0],[3, 0, 0, 0, 1, 1],[3, 1, 1, 0, 0, 0]],
[[1, 0, 1, 1, 1, 1],[2, 1, 1, 1, 0, 1],[3, 0, 0, 1, 1, 0],[3, 0, 1, 0, 1, 0],[3, 1, 1, 1, 1, 0]],
[[1, 1, 0, 1, 1, 0],[2, 1, 1, 0, 1, 1],[3, 0, 0, 1, 1, 1],[3, 1, 0, 0, 0, 1],[3, 1, 1, 0, 1, 1]],
[[2, 0, 1, 0, 0, 0],[2, 0, 1, 0, 1, 1],[2, 1, 0, 0, 0, 1],[3, 0, 1, 0, 1, 1],[3, 1, 1, 0, 0, 1]],
[[1, 1, 0, 1, 0, 1],[2, 1, 0, 0, 0, 0],[3, 0, 1, 1, 1, 0],[3, 1, 0, 1, 0, 0],[3, 1, 1, 1, 1, 1]],
[[1, 1, 1, 1, 1, 1],[2, 1, 1, 1, 1, 1],[3, 0, 1, 1, 0, 0],[3, 1, 0, 1, 1, 0],[3, 1, 1, 0, 1, 0]],
[[0, 1, 0, 1, 0, 1],[1, 0, 0, 0, 1, 1],[1, 0, 1, 1, 0, 1],[3, 0, 1, 0, 0, 1],[3, 1, 0, 0, 1, 0]],
[[0, 0, 0, 1, 1, 0],[0, 1, 1, 1, 0, 1],[2, 0, 1, 0, 1, 0],[3, 0, 1, 1, 0, 1],[3, 1, 1, 1, 0, 0]]
]
A partition for sets of sizes: 26*3  1*4  9*5
\end{verbatim}
\newpage

\begin{verbatim} 
[4, 2, 2, 2, 2, 2]
[23, 2, 10]
[
[[0, 0, 0, 0, 0, 1],[0, 0, 0, 0, 1, 0],[0, 0, 0, 0, 1, 1]],
[[0, 0, 0, 1, 0, 0],[0, 0, 1, 0, 0, 0],[0, 0, 1, 1, 0, 0]],
[[0, 0, 0, 1, 1, 1],[0, 0, 1, 0, 0, 1],[0, 0, 1, 1, 1, 0]],
[[0, 0, 1, 0, 1, 0],[0, 1, 0, 0, 0, 0],[0, 1, 1, 0, 1, 0]],
[[0, 0, 1, 1, 0, 1],[0, 1, 0, 0, 0, 1],[0, 1, 1, 1, 0, 0]],
[[0, 0, 1, 0, 1, 1],[0, 1, 0, 0, 1, 0],[0, 1, 1, 0, 0, 1]],
[[0, 1, 0, 0, 1, 1],[1, 0, 0, 0, 0, 0],[3, 1, 0, 0, 1, 1]],
[[0, 1, 0, 1, 1, 0],[1, 0, 0, 0, 0, 1],[3, 1, 0, 1, 1, 1]],
[[0, 0, 0, 1, 0, 1],[0, 1, 1, 0, 1, 1],[0, 1, 1, 1, 1, 0]],
[[0, 0, 1, 1, 1, 1],[0, 1, 0, 1, 1, 1],[0, 1, 1, 0, 0, 0]],
[[0, 1, 1, 1, 1, 1],[1, 0, 0, 0, 1, 0],[3, 1, 1, 1, 0, 1]],
[[0, 1, 0, 1, 0, 0],[1, 0, 0, 1, 0, 0],[3, 1, 0, 0, 0, 0]],
[[1, 0, 0, 1, 0, 1],[1, 0, 0, 1, 1, 0],[2, 0, 0, 0, 1, 1]],
[[1, 0, 1, 0, 0, 0],[1, 0, 1, 0, 0, 1],[2, 0, 0, 0, 0, 1]],
[[1, 0, 1, 0, 1, 1],[1, 0, 1, 1, 0, 0],[2, 0, 0, 1, 1, 1]],
[[1, 0, 1, 1, 1, 0],[1, 1, 0, 0, 0, 0],[2, 1, 1, 1, 1, 0]],
[[1, 1, 0, 0, 0, 1],[1, 1, 0, 0, 1, 1],[2, 0, 0, 0, 1, 0]],
[[1, 1, 0, 1, 0, 0],[1, 1, 1, 0, 0, 0],[2, 0, 1, 1, 0, 0]],
[[1, 1, 0, 1, 1, 1],[1, 1, 1, 0, 0, 1],[2, 0, 1, 1, 1, 0]],
[[1, 1, 1, 0, 1, 0],[1, 1, 1, 1, 0, 0],[2, 0, 0, 1, 1, 0]],
[[1, 0, 1, 0, 1, 0],[1, 1, 1, 1, 1, 0],[2, 1, 0, 1, 0, 0]],
[[1, 1, 0, 0, 1, 0],[1, 1, 1, 1, 0, 1],[2, 0, 1, 1, 1, 1]],
[[1, 0, 0, 1, 1, 1],[1, 1, 1, 0, 1, 1],[2, 1, 1, 1, 0, 0]],
[[2, 0, 0, 1, 0, 0],[2, 0, 1, 0, 0, 0],[2, 1, 0, 0, 0, 1],[2, 1, 1, 1, 0, 1]],
[[2, 0, 1, 0, 1, 0],[2, 0, 1, 0, 1, 1],[2, 1, 0, 0, 1, 0],[2, 1, 0, 0, 1, 1]],
[[1, 1, 0, 1, 0, 1],[2, 0, 0, 0, 0, 0],[3, 0, 0, 0, 0, 0],[3, 0, 0, 0, 0, 1],[3, 1, 0, 1, 0, 0]],
[[2, 0, 1, 1, 0, 1],[2, 1, 0, 1, 0, 1],[2, 1, 0, 1, 1, 0],[3, 0, 0, 0, 1, 0],[3, 0, 1, 1, 0, 0]],
[[2, 1, 1, 0, 0, 0],[2, 1, 1, 0, 0, 1],[2, 1, 1, 0, 1, 0],[3, 0, 0, 0, 1, 1],[3, 1, 1, 0, 0, 0]],
[[2, 0, 1, 0, 0, 1],[2, 1, 0, 0, 0, 0],[2, 1, 1, 1, 1, 1],[3, 0, 1, 0, 0, 0],[3, 0, 1, 1, 1, 0]],
[[1, 1, 1, 1, 1, 1],[2, 1, 1, 0, 1, 1],[3, 0, 0, 1, 0, 1],[3, 0, 0, 1, 1, 0],[3, 0, 0, 1, 1, 1]],
[[1, 0, 1, 1, 1, 1],[2, 1, 0, 1, 1, 1],[3, 0, 1, 0, 0, 1],[3, 0, 1, 0, 1, 0],[3, 1, 1, 0, 1, 1]],
[[1, 0, 1, 1, 0, 1],[2, 0, 0, 1, 0, 1],[3, 0, 0, 1, 0, 0],[3, 1, 0, 0, 1, 0],[3, 1, 1, 1, 1, 0]],
[[0, 1, 0, 1, 0, 1],[3, 0, 1, 1, 0, 1],[3, 1, 0, 0, 0, 1],[3, 1, 0, 1, 1, 0],[3, 1, 1, 1, 1, 1]],
[[0, 0, 0, 1, 1, 0],[1, 0, 0, 0, 1, 1],[1, 1, 0, 1, 1, 0],[3, 0, 1, 1, 1, 1],[3, 1, 1, 1, 0, 0]],
[[0, 1, 1, 1, 0, 1],[3, 0, 1, 0, 1, 1],[3, 1, 0, 1, 0, 1],[3, 1, 1, 0, 0, 1],[3, 1, 1, 0, 1, 0]]
]
A partition for sets of sizes: 23*3  2*4 10*5
\end{verbatim}
\newpage

\begin{verbatim} 
[4, 2, 2, 2, 2, 2]
[24, 0, 11]
[
[[0, 0, 0, 0, 0, 1],[0, 0, 0, 0, 1, 0],[0, 0, 0, 0, 1, 1]],
[[0, 0, 0, 1, 0, 0],[0, 0, 1, 0, 0, 0],[0, 0, 1, 1, 0, 0]],
[[0, 0, 0, 1, 1, 1],[0, 0, 1, 0, 0, 1],[0, 0, 1, 1, 1, 0]],
[[0, 0, 1, 0, 1, 0],[0, 1, 0, 0, 0, 0],[0, 1, 1, 0, 1, 0]],
[[0, 0, 1, 1, 0, 1],[0, 1, 0, 0, 0, 1],[0, 1, 1, 1, 0, 0]],
[[0, 0, 1, 0, 1, 1],[0, 1, 0, 0, 1, 0],[0, 1, 1, 0, 0, 1]],
[[0, 1, 0, 0, 1, 1],[1, 0, 0, 0, 0, 0],[3, 1, 0, 0, 1, 1]],
[[0, 1, 0, 1, 1, 0],[1, 0, 0, 0, 0, 1],[3, 1, 0, 1, 1, 1]],
[[0, 0, 0, 1, 0, 1],[0, 1, 1, 0, 1, 1],[0, 1, 1, 1, 1, 0]],
[[0, 0, 1, 1, 1, 1],[0, 1, 0, 1, 1, 1],[0, 1, 1, 0, 0, 0]],
[[0, 1, 1, 1, 1, 1],[1, 0, 0, 0, 1, 0],[3, 1, 1, 1, 0, 1]],
[[0, 1, 0, 1, 0, 0],[1, 0, 0, 1, 0, 0],[3, 1, 0, 0, 0, 0]],
[[1, 0, 0, 1, 0, 1],[1, 0, 0, 1, 1, 0],[2, 0, 0, 0, 1, 1]],
[[1, 0, 1, 0, 0, 0],[1, 0, 1, 0, 0, 1],[2, 0, 0, 0, 0, 1]],
[[1, 0, 1, 0, 1, 1],[1, 0, 1, 1, 0, 0],[2, 0, 0, 1, 1, 1]],
[[1, 0, 1, 1, 1, 0],[1, 1, 0, 0, 0, 0],[2, 1, 1, 1, 1, 0]],
[[1, 1, 0, 0, 0, 1],[1, 1, 0, 0, 1, 1],[2, 0, 0, 0, 1, 0]],
[[1, 1, 0, 1, 0, 0],[1, 1, 1, 0, 0, 0],[2, 0, 1, 1, 0, 0]],
[[1, 1, 0, 1, 1, 1],[1, 1, 1, 0, 0, 1],[2, 0, 1, 1, 1, 0]],
[[1, 1, 1, 0, 1, 0],[1, 1, 1, 1, 0, 0],[2, 0, 0, 1, 1, 0]],
[[1, 0, 1, 0, 1, 0],[1, 1, 1, 1, 1, 0],[2, 1, 0, 1, 0, 0]],
[[1, 1, 0, 0, 1, 0],[1, 1, 1, 1, 0, 1],[2, 0, 1, 1, 1, 1]],
[[1, 0, 0, 1, 1, 1],[1, 1, 1, 0, 1, 1],[2, 1, 1, 1, 0, 0]],
[[2, 0, 0, 1, 0, 0],[3, 0, 0, 0, 0, 0],[3, 0, 0, 1, 0, 0]],
[[2, 0, 1, 0, 0, 1],[2, 0, 1, 0, 1, 0],[2, 0, 1, 0, 1, 1],[3, 0, 0, 0, 0, 1],[3, 0, 1, 0, 0, 1]],
[[1, 1, 0, 1, 0, 1],[2, 0, 0, 0, 0, 0],[3, 0, 0, 0, 1, 0],[3, 0, 0, 0, 1, 1],[3, 1, 0, 1, 0, 0]],
[[2, 1, 0, 0, 1, 1],[2, 1, 0, 1, 0, 1],[2, 1, 0, 1, 1, 0],[3, 0, 0, 1, 0, 1],[3, 1, 0, 1, 0, 1]],
[[2, 1, 1, 0, 0, 0],[2, 1, 1, 0, 0, 1],[2, 1, 1, 0, 1, 0],[3, 0, 0, 1, 1, 1],[3, 1, 1, 1, 0, 0]],
[[1, 0, 1, 1, 1, 1],[2, 1, 1, 1, 0, 1],[3, 0, 0, 1, 1, 0],[3, 0, 1, 0, 1, 0],[3, 1, 1, 1, 1, 0]],
[[2, 1, 0, 0, 0, 0],[2, 1, 0, 0, 0, 1],[2, 1, 1, 1, 1, 1],[3, 0, 1, 0, 0, 0],[3, 1, 0, 1, 1, 0]],
[[1, 1, 1, 1, 1, 1],[2, 1, 1, 0, 1, 1],[3, 0, 1, 1, 0, 0],[3, 1, 0, 0, 0, 1],[3, 1, 1, 0, 0, 1]],
[[2, 0, 1, 1, 0, 1],[2, 1, 0, 0, 1, 0],[2, 1, 0, 1, 1, 1],[3, 1, 0, 0, 1, 0],[3, 1, 1, 0, 1, 0]],
[[1, 0, 1, 1, 0, 1],[2, 0, 0, 1, 0, 1],[3, 0, 1, 0, 1, 1],[3, 0, 1, 1, 0, 1],[3, 0, 1, 1, 1, 0]],
[[0, 0, 0, 1, 1, 0],[0, 1, 0, 1, 0, 1],[0, 1, 1, 1, 0, 1],[1, 1, 0, 1, 1, 0],[3, 1, 1, 0, 0, 0]],
[[1, 0, 0, 0, 1, 1],[2, 0, 1, 0, 0, 0],[3, 0, 1, 1, 1, 1],[3, 1, 1, 0, 1, 1],[3, 1, 1, 1, 1, 1]]
]
A partition for sets of sizes: 24*3  0*4 11*5
\end{verbatim}
\newpage

\begin{verbatim} 
[4, 2, 2, 2, 2, 2]
[20, 3, 11]
[
[[0, 0, 0, 0, 0, 1],[0, 0, 0, 0, 1, 0],[0, 0, 0, 0, 1, 1]],
[[0, 0, 0, 1, 0, 0],[0, 0, 1, 0, 0, 0],[0, 0, 1, 1, 0, 0]],
[[0, 0, 0, 1, 1, 1],[0, 0, 1, 0, 0, 1],[0, 0, 1, 1, 1, 0]],
[[0, 0, 1, 0, 1, 0],[0, 1, 0, 0, 0, 0],[0, 1, 1, 0, 1, 0]],
[[0, 0, 1, 1, 0, 1],[0, 1, 0, 0, 0, 1],[0, 1, 1, 1, 0, 0]],
[[0, 0, 1, 0, 1, 1],[0, 1, 0, 0, 1, 0],[0, 1, 1, 0, 0, 1]],
[[0, 1, 0, 0, 1, 1],[1, 0, 0, 0, 0, 0],[3, 1, 0, 0, 1, 1]],
[[0, 1, 0, 1, 1, 0],[1, 0, 0, 0, 0, 1],[3, 1, 0, 1, 1, 1]],
[[0, 0, 0, 1, 0, 1],[0, 1, 1, 0, 1, 1],[0, 1, 1, 1, 1, 0]],
[[0, 0, 1, 1, 1, 1],[0, 1, 0, 1, 1, 1],[0, 1, 1, 0, 0, 0]],
[[0, 1, 1, 1, 1, 1],[1, 0, 0, 0, 1, 0],[3, 1, 1, 1, 0, 1]],
[[0, 1, 0, 1, 0, 0],[1, 0, 0, 1, 0, 0],[3, 1, 0, 0, 0, 0]],
[[1, 0, 0, 1, 0, 1],[1, 0, 0, 1, 1, 0],[2, 0, 0, 0, 1, 1]],
[[1, 0, 1, 0, 0, 0],[1, 0, 1, 0, 0, 1],[2, 0, 0, 0, 0, 1]],
[[1, 0, 1, 0, 1, 1],[1, 0, 1, 1, 0, 0],[2, 0, 0, 1, 1, 1]],
[[1, 0, 1, 1, 1, 0],[1, 1, 0, 0, 0, 0],[2, 1, 1, 1, 1, 0]],
[[1, 1, 0, 0, 0, 1],[1, 1, 0, 0, 1, 1],[2, 0, 0, 0, 1, 0]],
[[1, 1, 0, 1, 0, 0],[1, 1, 1, 0, 0, 0],[2, 0, 1, 1, 0, 0]],
[[1, 1, 0, 1, 1, 1],[1, 1, 1, 0, 0, 1],[2, 0, 1, 1, 1, 0]],
[[1, 1, 1, 0, 1, 0],[1, 1, 1, 1, 0, 0],[2, 0, 0, 1, 1, 0]],
[[1, 1, 1, 1, 0, 1],[1, 1, 1, 1, 1, 0],[3, 0, 0, 0, 0, 0],[3, 0, 0, 0, 1, 1]],
[[1, 0, 1, 0, 1, 0],[1, 1, 0, 0, 1, 0],[3, 0, 0, 0, 0, 1],[3, 1, 1, 0, 0, 1]],
[[2, 0, 0, 1, 0, 1],[2, 0, 1, 0, 0, 0],[2, 1, 0, 0, 0, 0],[2, 1, 1, 1, 0, 1]],
[[2, 0, 1, 0, 0, 1],[2, 0, 1, 0, 1, 0],[2, 0, 1, 0, 1, 1],[3, 0, 0, 0, 1, 0],[3, 0, 1, 0, 1, 0]],
[[1, 1, 0, 1, 0, 1],[2, 0, 1, 1, 1, 1],[3, 0, 0, 1, 0, 0],[3, 0, 0, 1, 0, 1],[3, 1, 1, 0, 1, 1]],
[[2, 1, 0, 0, 1, 1],[2, 1, 0, 1, 0, 0],[2, 1, 0, 1, 0, 1],[3, 0, 0, 1, 1, 0],[3, 1, 0, 1, 0, 0]],
[[2, 1, 1, 0, 0, 0],[2, 1, 1, 0, 0, 1],[2, 1, 1, 0, 1, 0],[3, 0, 0, 1, 1, 1],[3, 1, 1, 1, 0, 0]],
[[1, 1, 1, 0, 1, 1],[2, 1, 1, 1, 1, 1],[3, 0, 1, 0, 0, 0],[3, 1, 0, 0, 1, 0],[3, 1, 1, 1, 1, 0]],
[[1, 1, 0, 1, 1, 0],[2, 0, 0, 0, 0, 0],[3, 0, 1, 0, 0, 1],[3, 0, 1, 1, 1, 0],[3, 1, 0, 0, 0, 1]],
[[1, 0, 1, 1, 1, 1],[2, 1, 1, 0, 1, 1],[3, 0, 1, 1, 0, 0],[3, 0, 1, 1, 0, 1],[3, 1, 0, 1, 0, 1]],
[[1, 0, 1, 1, 0, 1],[2, 1, 0, 0, 0, 1],[3, 0, 1, 0, 1, 1],[3, 0, 1, 1, 1, 1],[3, 1, 1, 0, 0, 0]],
[[0, 1, 0, 1, 0, 1],[2, 0, 1, 1, 0, 1],[2, 1, 0, 0, 1, 0],[2, 1, 0, 1, 1, 0],[2, 1, 1, 1, 0, 0]],
[[0, 0, 0, 1, 1, 0],[1, 0, 0, 0, 1, 1],[2, 0, 0, 1, 0, 0],[2, 1, 0, 1, 1, 1],[3, 1, 0, 1, 1, 0]],
[[0, 1, 1, 1, 0, 1],[1, 0, 0, 1, 1, 1],[1, 1, 1, 1, 1, 1],[3, 1, 1, 0, 1, 0],[3, 1, 1, 1, 1, 1]]
]
A partition for sets of sizes: 20*3  3*4 11*5
\end{verbatim}
\newpage

\begin{verbatim} 
[4, 2, 2, 2, 2, 2]
[21, 1, 12]
[
[[0, 0, 0, 0, 0, 1],[0, 0, 0, 0, 1, 0],[0, 0, 0, 0, 1, 1]],
[[0, 0, 0, 1, 0, 0],[0, 0, 1, 0, 0, 0],[0, 0, 1, 1, 0, 0]],
[[0, 0, 0, 1, 1, 1],[0, 0, 1, 0, 0, 1],[0, 0, 1, 1, 1, 0]],
[[0, 0, 1, 0, 1, 0],[0, 1, 0, 0, 0, 0],[0, 1, 1, 0, 1, 0]],
[[0, 0, 1, 1, 0, 1],[0, 1, 0, 0, 0, 1],[0, 1, 1, 1, 0, 0]],
[[0, 0, 1, 0, 1, 1],[0, 1, 0, 0, 1, 0],[0, 1, 1, 0, 0, 1]],
[[0, 1, 0, 0, 1, 1],[1, 0, 0, 0, 0, 0],[3, 1, 0, 0, 1, 1]],
[[0, 1, 0, 1, 1, 0],[1, 0, 0, 0, 0, 1],[3, 1, 0, 1, 1, 1]],
[[0, 0, 0, 1, 0, 1],[0, 1, 1, 0, 1, 1],[0, 1, 1, 1, 1, 0]],
[[0, 0, 1, 1, 1, 1],[0, 1, 0, 1, 1, 1],[0, 1, 1, 0, 0, 0]],
[[0, 1, 1, 1, 1, 1],[1, 0, 0, 0, 1, 0],[3, 1, 1, 1, 0, 1]],
[[0, 1, 0, 1, 0, 0],[1, 0, 0, 1, 0, 0],[3, 1, 0, 0, 0, 0]],
[[1, 0, 0, 1, 0, 1],[1, 0, 0, 1, 1, 0],[2, 0, 0, 0, 1, 1]],
[[1, 0, 1, 0, 0, 0],[1, 0, 1, 0, 0, 1],[2, 0, 0, 0, 0, 1]],
[[1, 0, 1, 0, 1, 1],[1, 0, 1, 1, 0, 0],[2, 0, 0, 1, 1, 1]],
[[1, 0, 1, 1, 1, 0],[1, 1, 0, 0, 0, 0],[2, 1, 1, 1, 1, 0]],
[[1, 1, 0, 0, 0, 1],[1, 1, 0, 0, 1, 1],[2, 0, 0, 0, 1, 0]],
[[1, 1, 0, 1, 0, 0],[1, 1, 1, 0, 0, 0],[2, 0, 1, 1, 0, 0]],
[[1, 1, 0, 1, 1, 1],[1, 1, 1, 0, 0, 1],[2, 0, 1, 1, 1, 0]],
[[1, 1, 1, 0, 1, 0],[1, 1, 1, 1, 0, 0],[2, 0, 0, 1, 1, 0]],
[[1, 0, 1, 0, 1, 0],[1, 1, 1, 1, 1, 0],[2, 1, 0, 1, 0, 0]],
[[2, 0, 0, 0, 0, 0],[2, 0, 0, 1, 0, 0],[2, 0, 1, 0, 0, 1],[2, 0, 1, 1, 0, 1]],
[[1, 1, 1, 1, 0, 1],[2, 0, 0, 1, 0, 1],[3, 0, 0, 0, 0, 0],[3, 0, 0, 0, 0, 1],[3, 1, 1, 0, 0, 1]],
[[1, 1, 0, 0, 1, 0],[2, 0, 1, 0, 1, 0],[3, 0, 0, 0, 1, 0],[3, 0, 0, 1, 0, 0],[3, 1, 1, 1, 1, 0]],
[[1, 1, 0, 1, 0, 1],[2, 0, 1, 1, 1, 1],[3, 0, 0, 0, 1, 1],[3, 0, 0, 1, 0, 1],[3, 1, 1, 1, 0, 0]],
[[2, 1, 0, 0, 1, 1],[2, 1, 0, 1, 0, 1],[2, 1, 0, 1, 1, 0],[3, 0, 0, 1, 1, 0],[3, 1, 0, 1, 1, 0]],
[[2, 1, 1, 0, 0, 0],[2, 1, 1, 0, 0, 1],[2, 1, 1, 0, 1, 0],[3, 0, 1, 0, 0, 1],[3, 1, 0, 0, 1, 0]],
[[1, 0, 1, 1, 1, 1],[2, 1, 1, 1, 0, 1],[3, 0, 0, 1, 1, 1],[3, 0, 1, 0, 1, 0],[3, 1, 1, 1, 1, 1]],
[[2, 0, 1, 0, 1, 1],[2, 1, 0, 0, 0, 0],[2, 1, 0, 0, 0, 1],[3, 1, 0, 0, 0, 1],[3, 1, 1, 0, 1, 1]],
[[2, 1, 1, 0, 1, 1],[2, 1, 1, 1, 0, 0],[2, 1, 1, 1, 1, 1],[3, 0, 1, 1, 0, 0],[3, 1, 0, 1, 0, 0]],
[[1, 1, 1, 0, 1, 1],[2, 1, 0, 1, 1, 1],[3, 0, 1, 1, 1, 0],[3, 1, 1, 0, 0, 0],[3, 1, 1, 0, 1, 0]],
[[0, 1, 0, 1, 0, 1],[1, 0, 0, 1, 1, 1],[2, 0, 1, 0, 0, 0],[2, 1, 0, 0, 1, 0],[3, 0, 1, 0, 0, 0]],
[[0, 0, 0, 1, 1, 0],[1, 0, 0, 0, 1, 1],[1, 1, 1, 1, 1, 1],[3, 0, 1, 1, 1, 1],[3, 1, 0, 1, 0, 1]],
[[0, 1, 1, 1, 0, 1],[1, 0, 1, 1, 0, 1],[1, 1, 0, 1, 1, 0],[3, 0, 1, 0, 1, 1],[3, 0, 1, 1, 0, 1]]
]
A partition for sets of sizes: 21*3  1*4 12*5
\end{verbatim}
\newpage

\begin{verbatim} 
[4, 2, 2, 2, 2, 2]
[18, 2, 13]
[
[[0, 0, 0, 0, 0, 1],[0, 0, 0, 0, 1, 0],[0, 0, 0, 0, 1, 1]],
[[0, 0, 0, 1, 0, 0],[0, 0, 1, 0, 0, 0],[0, 0, 1, 1, 0, 0]],
[[0, 0, 0, 1, 1, 1],[0, 0, 1, 0, 0, 1],[0, 0, 1, 1, 1, 0]],
[[0, 0, 1, 0, 1, 0],[0, 1, 0, 0, 0, 0],[0, 1, 1, 0, 1, 0]],
[[0, 0, 1, 1, 0, 1],[0, 1, 0, 0, 0, 1],[0, 1, 1, 1, 0, 0]],
[[0, 0, 1, 0, 1, 1],[0, 1, 0, 0, 1, 0],[0, 1, 1, 0, 0, 1]],
[[0, 1, 0, 0, 1, 1],[1, 0, 0, 0, 0, 0],[3, 1, 0, 0, 1, 1]],
[[0, 1, 0, 1, 1, 0],[1, 0, 0, 0, 0, 1],[3, 1, 0, 1, 1, 1]],
[[0, 0, 0, 1, 0, 1],[0, 1, 1, 0, 1, 1],[0, 1, 1, 1, 1, 0]],
[[0, 0, 1, 1, 1, 1],[0, 1, 0, 1, 1, 1],[0, 1, 1, 0, 0, 0]],
[[0, 1, 1, 1, 1, 1],[1, 0, 0, 0, 1, 0],[3, 1, 1, 1, 0, 1]],
[[0, 1, 0, 1, 0, 0],[1, 0, 0, 1, 0, 0],[3, 1, 0, 0, 0, 0]],
[[1, 0, 0, 1, 0, 1],[1, 0, 0, 1, 1, 0],[2, 0, 0, 0, 1, 1]],
[[1, 0, 1, 0, 0, 0],[1, 0, 1, 0, 0, 1],[2, 0, 0, 0, 0, 1]],
[[1, 0, 1, 0, 1, 1],[1, 0, 1, 1, 0, 0],[2, 0, 0, 1, 1, 1]],
[[1, 0, 1, 1, 1, 0],[1, 1, 0, 0, 0, 0],[2, 1, 1, 1, 1, 0]],
[[1, 1, 0, 0, 0, 1],[1, 1, 0, 0, 1, 1],[2, 0, 0, 0, 1, 0]],
[[1, 1, 0, 1, 0, 0],[1, 1, 1, 0, 0, 0],[2, 0, 1, 1, 0, 0]],
[[1, 1, 0, 1, 1, 1],[1, 1, 1, 0, 0, 1],[3, 0, 0, 0, 0, 0],[3, 0, 1, 1, 1, 0]],
[[1, 1, 1, 0, 1, 1],[1, 1, 1, 1, 0, 0],[3, 0, 0, 0, 0, 1],[3, 0, 0, 1, 1, 0]],
[[1, 1, 1, 1, 1, 1],[2, 0, 0, 0, 0, 0],[3, 0, 0, 0, 1, 0],[3, 0, 0, 0, 1, 1],[3, 1, 1, 1, 1, 0]],
[[2, 0, 0, 1, 0, 0],[2, 0, 0, 1, 0, 1],[2, 0, 0, 1, 1, 0],[3, 0, 1, 0, 0, 0],[3, 0, 1, 1, 1, 1]],
[[2, 0, 1, 0, 0, 1],[2, 0, 1, 0, 1, 0],[2, 0, 1, 0, 1, 1],[3, 0, 0, 1, 0, 0],[3, 0, 1, 1, 0, 0]],
[[2, 0, 1, 1, 1, 0],[2, 0, 1, 1, 1, 1],[2, 1, 0, 0, 0, 0],[3, 0, 0, 1, 0, 1],[3, 1, 0, 1, 0, 0]],
[[2, 1, 0, 0, 1, 1],[2, 1, 0, 1, 0, 0],[2, 1, 0, 1, 0, 1],[3, 0, 0, 1, 1, 1],[3, 1, 0, 1, 0, 1]],
[[2, 1, 1, 0, 0, 0],[2, 1, 1, 0, 0, 1],[2, 1, 1, 0, 1, 0],[3, 0, 1, 0, 0, 1],[3, 1, 0, 0, 1, 0]],
[[1, 0, 1, 1, 1, 1],[2, 1, 1, 1, 0, 1],[3, 1, 0, 0, 0, 1],[3, 1, 1, 0, 0, 0],[3, 1, 1, 0, 1, 1]],
[[1, 0, 1, 0, 1, 0],[1, 1, 0, 0, 1, 0],[1, 1, 0, 1, 1, 0],[2, 1, 0, 0, 0, 1],[3, 1, 1, 1, 1, 1]],
[[2, 1, 0, 1, 1, 0],[2, 1, 1, 0, 1, 1],[2, 1, 1, 1, 1, 1],[3, 0, 1, 0, 1, 1],[3, 1, 1, 0, 0, 1]],
[[2, 0, 1, 1, 0, 1],[2, 1, 0, 1, 1, 1],[2, 1, 1, 1, 0, 0],[3, 1, 1, 0, 1, 0],[3, 1, 1, 1, 0, 0]],
[[1, 0, 1, 1, 0, 1],[1, 1, 0, 1, 0, 1],[1, 1, 1, 0, 1, 0],[2, 0, 1, 0, 0, 0],[3, 0, 1, 0, 1, 0]],
[[0, 0, 0, 1, 1, 0],[0, 1, 0, 1, 0, 1],[0, 1, 1, 1, 0, 1],[1, 0, 0, 0, 1, 1],[3, 0, 1, 1, 0, 1]],
[[1, 0, 0, 1, 1, 1],[1, 1, 1, 1, 0, 1],[1, 1, 1, 1, 1, 0],[2, 1, 0, 0, 1, 0],[3, 1, 0, 1, 1, 0]]
]
A partition for sets of sizes: 18*3  2*4 13*5
\end{verbatim}
\newpage

\begin{verbatim} 
[4, 2, 2, 2, 2, 2]
[19, 0, 14]
[
[[0, 0, 0, 0, 0, 1],[0, 0, 0, 0, 1, 0],[0, 0, 0, 0, 1, 1]],
[[0, 0, 0, 1, 0, 0],[0, 0, 1, 0, 0, 0],[0, 0, 1, 1, 0, 0]],
[[0, 0, 0, 1, 1, 1],[0, 0, 1, 0, 0, 1],[0, 0, 1, 1, 1, 0]],
[[0, 0, 1, 0, 1, 0],[0, 1, 0, 0, 0, 0],[0, 1, 1, 0, 1, 0]],
[[0, 0, 1, 1, 0, 1],[0, 1, 0, 0, 0, 1],[0, 1, 1, 1, 0, 0]],
[[0, 0, 1, 0, 1, 1],[0, 1, 0, 0, 1, 0],[0, 1, 1, 0, 0, 1]],
[[0, 1, 0, 0, 1, 1],[1, 0, 0, 0, 0, 0],[3, 1, 0, 0, 1, 1]],
[[0, 1, 0, 1, 1, 0],[1, 0, 0, 0, 0, 1],[3, 1, 0, 1, 1, 1]],
[[0, 0, 0, 1, 0, 1],[0, 1, 1, 0, 1, 1],[0, 1, 1, 1, 1, 0]],
[[0, 0, 1, 1, 1, 1],[0, 1, 0, 1, 1, 1],[0, 1, 1, 0, 0, 0]],
[[0, 1, 1, 1, 1, 1],[1, 0, 0, 0, 1, 0],[3, 1, 1, 1, 0, 1]],
[[0, 1, 0, 1, 0, 0],[1, 0, 0, 1, 0, 0],[3, 1, 0, 0, 0, 0]],
[[1, 0, 0, 1, 0, 1],[1, 0, 0, 1, 1, 0],[2, 0, 0, 0, 1, 1]],
[[1, 0, 1, 0, 0, 0],[1, 0, 1, 0, 0, 1],[2, 0, 0, 0, 0, 1]],
[[1, 0, 1, 0, 1, 1],[1, 0, 1, 1, 0, 0],[2, 0, 0, 1, 1, 1]],
[[1, 0, 1, 1, 1, 0],[1, 1, 0, 0, 0, 0],[2, 1, 1, 1, 1, 0]],
[[1, 1, 0, 0, 0, 1],[1, 1, 0, 0, 1, 1],[2, 0, 0, 0, 1, 0]],
[[1, 1, 0, 1, 0, 0],[1, 1, 1, 0, 0, 0],[2, 0, 1, 1, 0, 0]],
[[1, 1, 0, 1, 1, 1],[1, 1, 1, 0, 0, 1],[2, 0, 1, 1, 1, 0]],
[[1, 1, 1, 0, 1, 0],[1, 1, 1, 0, 1, 1],[1, 1, 1, 1, 0, 0],[2, 0, 0, 1, 0, 0],[3, 1, 1, 0, 0, 1]],
[[1, 1, 1, 1, 1, 1],[2, 0, 0, 0, 0, 0],[3, 0, 0, 0, 0, 0],[3, 0, 0, 0, 0, 1],[3, 1, 1, 1, 1, 0]],
[[1, 1, 1, 1, 0, 1],[2, 0, 0, 1, 0, 1],[3, 0, 0, 0, 1, 0],[3, 0, 0, 1, 0, 1],[3, 1, 1, 1, 1, 1]],
[[2, 0, 1, 0, 0, 1],[2, 0, 1, 0, 1, 0],[2, 0, 1, 0, 1, 1],[3, 0, 0, 0, 1, 1],[3, 0, 1, 0, 1, 1]],
[[1, 1, 0, 1, 0, 1],[2, 0, 1, 1, 1, 1],[3, 0, 0, 1, 0, 0],[3, 0, 0, 1, 1, 0],[3, 1, 1, 0, 0, 0]],
[[2, 1, 0, 0, 1, 1],[2, 1, 0, 1, 0, 0],[2, 1, 0, 1, 0, 1],[3, 0, 0, 1, 1, 1],[3, 1, 0, 1, 0, 1]],
[[2, 1, 1, 0, 0, 0],[2, 1, 1, 0, 0, 1],[2, 1, 1, 0, 1, 0],[3, 0, 1, 0, 0, 1],[3, 1, 0, 0, 1, 0]],
[[1, 0, 1, 1, 1, 1],[2, 1, 1, 1, 0, 1],[3, 0, 1, 0, 0, 0],[3, 0, 1, 1, 0, 0],[3, 1, 0, 1, 1, 0]],
[[2, 0, 0, 1, 1, 0],[2, 1, 0, 0, 0, 0],[2, 1, 0, 0, 0, 1],[3, 0, 1, 0, 1, 0],[3, 0, 1, 1, 0, 1]],
[[2, 1, 0, 1, 1, 0],[2, 1, 1, 0, 1, 1],[2, 1, 1, 1, 1, 1],[3, 0, 1, 1, 1, 0],[3, 1, 1, 1, 0, 0]],
[[1, 0, 1, 0, 1, 0],[1, 1, 0, 1, 1, 0],[1, 1, 1, 1, 1, 0],[2, 0, 1, 1, 0, 1],[3, 0, 1, 1, 1, 1]],
[[0, 1, 0, 1, 0, 1],[1, 0, 1, 1, 0, 1],[1, 1, 0, 0, 1, 0],[3, 1, 0, 0, 0, 1],[3, 1, 1, 0, 1, 1]],
[[0, 0, 0, 1, 1, 0],[1, 0, 0, 0, 1, 1],[2, 0, 1, 0, 0, 0],[2, 1, 0, 1, 1, 1],[3, 1, 1, 0, 1, 0]],
[[0, 1, 1, 1, 0, 1],[1, 0, 0, 1, 1, 1],[2, 1, 0, 0, 1, 0],[2, 1, 1, 1, 0, 0],[3, 1, 0, 1, 0, 0]]
]
A partition for sets of sizes: 19*3  0*4 14*5
\end{verbatim}
\newpage

\begin{verbatim} 
[4, 2, 2, 2, 2, 2]
[15, 3, 14]
[
[[0, 0, 0, 0, 0, 1],[0, 0, 0, 0, 1, 0],[0, 0, 0, 0, 1, 1]],
[[0, 0, 0, 1, 0, 0],[0, 0, 1, 0, 0, 0],[0, 0, 1, 1, 0, 0]],
[[0, 0, 0, 1, 1, 1],[0, 0, 1, 0, 0, 1],[0, 0, 1, 1, 1, 0]],
[[0, 0, 1, 0, 1, 0],[0, 1, 0, 0, 0, 0],[0, 1, 1, 0, 1, 0]],
[[0, 0, 1, 1, 0, 1],[0, 1, 0, 0, 0, 1],[0, 1, 1, 1, 0, 0]],
[[0, 0, 1, 0, 1, 1],[0, 1, 0, 0, 1, 0],[0, 1, 1, 0, 0, 1]],
[[0, 1, 0, 0, 1, 1],[1, 0, 0, 0, 0, 0],[3, 1, 0, 0, 1, 1]],
[[0, 1, 0, 1, 1, 0],[1, 0, 0, 0, 0, 1],[3, 1, 0, 1, 1, 1]],
[[0, 0, 0, 1, 0, 1],[0, 1, 1, 0, 1, 1],[0, 1, 1, 1, 1, 0]],
[[0, 0, 1, 1, 1, 1],[0, 1, 0, 1, 1, 1],[0, 1, 1, 0, 0, 0]],
[[0, 1, 1, 1, 1, 1],[1, 0, 0, 0, 1, 0],[3, 1, 1, 1, 0, 1]],
[[0, 1, 0, 1, 0, 0],[1, 0, 0, 1, 0, 0],[3, 1, 0, 0, 0, 0]],
[[1, 0, 0, 1, 0, 1],[1, 0, 0, 1, 1, 0],[2, 0, 0, 0, 1, 1]],
[[1, 0, 1, 0, 0, 0],[1, 0, 1, 0, 0, 1],[2, 0, 0, 0, 0, 1]],
[[1, 0, 1, 0, 1, 1],[1, 0, 1, 1, 0, 0],[2, 0, 0, 1, 1, 1]],
[[1, 0, 1, 1, 1, 0],[1, 0, 1, 1, 1, 1],[1, 1, 0, 0, 0, 0],[1, 1, 0, 0, 0, 1]],
[[1, 1, 0, 0, 1, 0],[1, 1, 0, 0, 1, 1],[1, 1, 0, 1, 0, 0],[1, 1, 0, 1, 0, 1]],
[[1, 1, 0, 1, 1, 0],[1, 1, 0, 1, 1, 1],[1, 1, 1, 0, 0, 0],[1, 1, 1, 0, 0, 1]],
[[1, 1, 1, 0, 1, 0],[1, 1, 1, 0, 1, 1],[1, 1, 1, 1, 0, 0],[2, 0, 0, 0, 1, 0],[3, 1, 1, 1, 1, 1]],
[[1, 1, 1, 1, 1, 1],[2, 0, 0, 0, 0, 0],[3, 0, 0, 0, 0, 0],[3, 0, 0, 0, 0, 1],[3, 1, 1, 1, 1, 0]],
[[2, 0, 0, 1, 0, 0],[2, 0, 0, 1, 0, 1],[2, 0, 0, 1, 1, 0],[3, 0, 0, 0, 1, 0],[3, 0, 0, 1, 0, 1]],
[[2, 0, 1, 0, 0, 1],[2, 0, 1, 0, 1, 0],[2, 0, 1, 0, 1, 1],[3, 0, 0, 0, 1, 1],[3, 0, 1, 0, 1, 1]],
[[2, 0, 1, 1, 1, 0],[2, 0, 1, 1, 1, 1],[2, 1, 0, 0, 0, 0],[3, 0, 0, 1, 0, 0],[3, 1, 0, 1, 0, 1]],
[[2, 1, 0, 0, 1, 1],[2, 1, 0, 1, 0, 0],[2, 1, 0, 1, 0, 1],[3, 0, 0, 1, 1, 0],[3, 1, 0, 1, 0, 0]],
[[2, 1, 1, 0, 0, 0],[2, 1, 1, 0, 0, 1],[2, 1, 1, 0, 1, 0],[3, 0, 0, 1, 1, 1],[3, 1, 1, 1, 0, 0]],
[[2, 1, 1, 1, 0, 1],[2, 1, 1, 1, 1, 0],[2, 1, 1, 1, 1, 1],[3, 0, 1, 0, 1, 0],[3, 1, 0, 1, 1, 0]],
[[1, 0, 1, 1, 0, 1],[2, 0, 1, 1, 0, 0],[3, 0, 1, 0, 0, 0],[3, 1, 0, 0, 0, 1],[3, 1, 1, 0, 0, 0]],
[[2, 1, 0, 0, 0, 1],[2, 1, 1, 0, 1, 1],[2, 1, 1, 1, 0, 0],[3, 0, 1, 1, 0, 0],[3, 1, 1, 0, 1, 0]],
[[1, 0, 1, 0, 1, 0],[2, 0, 1, 0, 0, 0],[3, 0, 1, 0, 0, 1],[3, 1, 0, 0, 1, 0],[3, 1, 1, 0, 0, 1]],
[[0, 1, 0, 1, 0, 1],[0, 1, 1, 1, 0, 1],[1, 0, 0, 0, 1, 1],[1, 1, 1, 1, 0, 1],[2, 1, 0, 1, 1, 0]],
[[0, 0, 0, 1, 1, 0],[1, 0, 0, 1, 1, 1],[2, 0, 1, 1, 0, 1],[2, 1, 0, 1, 1, 1],[3, 1, 1, 0, 1, 1]],
[[1, 1, 1, 1, 1, 0],[2, 1, 0, 0, 1, 0],[3, 0, 1, 1, 0, 1],[3, 0, 1, 1, 1, 0],[3, 0, 1, 1, 1, 1]]
]
A partition for sets of sizes: 15*3  3*4 14*5
\end{verbatim}
\newpage

\begin{verbatim} 
[4, 2, 2, 2, 2, 2]
[16, 1, 15]
[
[[0, 0, 0, 0, 0, 1],[0, 0, 0, 0, 1, 0],[0, 0, 0, 0, 1, 1]],
[[0, 0, 0, 1, 0, 0],[0, 0, 1, 0, 0, 0],[0, 0, 1, 1, 0, 0]],
[[0, 0, 0, 1, 1, 1],[0, 0, 1, 0, 0, 1],[0, 0, 1, 1, 1, 0]],
[[0, 0, 1, 0, 1, 0],[0, 1, 0, 0, 0, 0],[0, 1, 1, 0, 1, 0]],
[[0, 0, 1, 1, 0, 1],[0, 1, 0, 0, 0, 1],[0, 1, 1, 1, 0, 0]],
[[0, 0, 1, 0, 1, 1],[0, 1, 0, 0, 1, 0],[0, 1, 1, 0, 0, 1]],
[[0, 1, 0, 0, 1, 1],[1, 0, 0, 0, 0, 0],[3, 1, 0, 0, 1, 1]],
[[0, 1, 0, 1, 1, 0],[1, 0, 0, 0, 0, 1],[3, 1, 0, 1, 1, 1]],
[[0, 0, 0, 1, 0, 1],[0, 1, 1, 0, 1, 1],[0, 1, 1, 1, 1, 0]],
[[0, 0, 1, 1, 1, 1],[0, 1, 0, 1, 1, 1],[0, 1, 1, 0, 0, 0]],
[[0, 1, 1, 1, 1, 1],[1, 0, 0, 0, 1, 0],[3, 1, 1, 1, 0, 1]],
[[0, 1, 0, 1, 0, 0],[1, 0, 0, 1, 0, 0],[3, 1, 0, 0, 0, 0]],
[[1, 0, 0, 1, 0, 1],[1, 0, 0, 1, 1, 0],[2, 0, 0, 0, 1, 1]],
[[1, 0, 1, 0, 0, 0],[1, 0, 1, 0, 0, 1],[2, 0, 0, 0, 0, 1]],
[[1, 0, 1, 0, 1, 1],[1, 0, 1, 1, 0, 0],[2, 0, 0, 1, 1, 1]],
[[1, 0, 1, 1, 1, 0],[1, 1, 0, 0, 0, 0],[2, 1, 1, 1, 1, 0]],
[[1, 1, 0, 0, 0, 1],[1, 1, 0, 0, 1, 0],[1, 1, 0, 1, 0, 0],[1, 1, 0, 1, 1, 1]],
[[1, 1, 0, 1, 0, 1],[1, 1, 0, 1, 1, 0],[1, 1, 1, 0, 0, 0],[2, 0, 0, 0, 0, 0],[3, 1, 1, 0, 1, 1]],
[[1, 1, 1, 0, 1, 0],[1, 1, 1, 0, 1, 1],[1, 1, 1, 1, 0, 0],[2, 0, 0, 0, 1, 0],[3, 1, 1, 1, 1, 1]],
[[1, 1, 1, 1, 1, 1],[2, 0, 0, 1, 0, 0],[3, 0, 0, 0, 0, 0],[3, 0, 0, 0, 0, 1],[3, 1, 1, 0, 1, 0]],
[[1, 1, 0, 0, 1, 1],[2, 0, 0, 1, 0, 1],[3, 0, 0, 0, 1, 0],[3, 0, 0, 1, 0, 1],[3, 1, 0, 0, 0, 1]],
[[2, 0, 1, 0, 0, 1],[2, 0, 1, 0, 1, 0],[2, 0, 1, 0, 1, 1],[3, 0, 0, 0, 1, 1],[3, 0, 1, 0, 1, 1]],
[[2, 0, 1, 1, 1, 0],[2, 0, 1, 1, 1, 1],[2, 1, 0, 0, 0, 0],[3, 0, 0, 1, 0, 0],[3, 1, 0, 1, 0, 1]],
[[2, 1, 0, 0, 1, 1],[2, 1, 0, 1, 0, 0],[2, 1, 0, 1, 0, 1],[3, 0, 0, 1, 1, 0],[3, 1, 0, 1, 0, 0]],
[[2, 1, 1, 0, 0, 0],[2, 1, 1, 0, 0, 1],[2, 1, 1, 0, 1, 0],[3, 0, 0, 1, 1, 1],[3, 1, 1, 1, 0, 0]],
[[1, 0, 1, 1, 1, 1],[2, 1, 1, 1, 0, 1],[3, 0, 1, 0, 0, 0],[3, 0, 1, 1, 0, 0],[3, 1, 0, 1, 1, 0]],
[[2, 0, 0, 1, 1, 0],[2, 0, 1, 1, 0, 0],[2, 1, 0, 0, 0, 1],[3, 0, 1, 0, 0, 1],[3, 1, 0, 0, 1, 0]],
[[1, 0, 1, 1, 0, 1],[2, 1, 0, 1, 1, 0],[3, 0, 1, 0, 1, 0],[3, 0, 1, 1, 1, 1],[3, 1, 1, 1, 1, 0]],
[[1, 0, 1, 0, 1, 0],[1, 1, 1, 1, 0, 1],[2, 0, 1, 0, 0, 0],[2, 0, 1, 1, 0, 1],[2, 1, 0, 0, 1, 0]],
[[0, 1, 0, 1, 0, 1],[0, 1, 1, 1, 0, 1],[2, 1, 1, 1, 1, 1],[3, 0, 1, 1, 1, 0],[3, 1, 1, 0, 0, 1]],
[[0, 0, 0, 1, 1, 0],[1, 0, 0, 1, 1, 1],[2, 1, 0, 1, 1, 1],[2, 1, 1, 0, 1, 1],[3, 0, 1, 1, 0, 1]],
[[1, 0, 0, 0, 1, 1],[1, 1, 1, 0, 0, 1],[1, 1, 1, 1, 1, 0],[2, 1, 1, 1, 0, 0],[3, 1, 1, 0, 0, 0]]
]
A partition for sets of sizes: 16*3  1*4 15*5
\end{verbatim}
\newpage

\begin{verbatim} 
[4, 2, 2, 2, 2, 2]
[13, 2, 16]
[
[[0, 0, 0, 0, 0, 1],[0, 0, 0, 0, 1, 0],[0, 0, 0, 0, 1, 1]],
[[0, 0, 0, 1, 0, 0],[0, 0, 1, 0, 0, 0],[0, 0, 1, 1, 0, 0]],
[[0, 0, 0, 1, 1, 1],[0, 0, 1, 0, 0, 1],[0, 0, 1, 1, 1, 0]],
[[0, 0, 1, 0, 1, 0],[0, 1, 0, 0, 0, 0],[0, 1, 1, 0, 1, 0]],
[[0, 0, 1, 1, 0, 1],[0, 1, 0, 0, 0, 1],[0, 1, 1, 1, 0, 0]],
[[0, 0, 1, 0, 1, 1],[0, 1, 0, 0, 1, 0],[0, 1, 1, 0, 0, 1]],
[[0, 1, 0, 0, 1, 1],[1, 0, 0, 0, 0, 0],[3, 1, 0, 0, 1, 1]],
[[0, 1, 0, 1, 1, 0],[1, 0, 0, 0, 0, 1],[3, 1, 0, 1, 1, 1]],
[[0, 0, 0, 1, 0, 1],[0, 1, 1, 0, 1, 1],[0, 1, 1, 1, 1, 0]],
[[0, 0, 1, 1, 1, 1],[0, 1, 0, 1, 1, 1],[0, 1, 1, 0, 0, 0]],
[[0, 1, 1, 1, 1, 1],[1, 0, 0, 0, 1, 0],[3, 1, 1, 1, 0, 1]],
[[0, 1, 0, 1, 0, 0],[1, 0, 0, 1, 0, 0],[3, 1, 0, 0, 0, 0]],
[[1, 0, 0, 1, 0, 1],[1, 0, 0, 1, 1, 0],[2, 0, 0, 0, 1, 1]],
[[1, 0, 1, 0, 0, 0],[1, 0, 1, 0, 0, 1],[1, 0, 1, 0, 1, 0],[1, 0, 1, 0, 1, 1]],
[[1, 0, 1, 1, 0, 0],[1, 0, 1, 1, 0, 1],[1, 0, 1, 1, 1, 0],[1, 0, 1, 1, 1, 1]],
[[1, 1, 0, 0, 0, 0],[1, 1, 0, 0, 0, 1],[1, 1, 0, 0, 1, 0],[2, 0, 0, 0, 0, 1],[3, 1, 0, 0, 1, 0]],
[[1, 1, 0, 1, 0, 1],[1, 1, 0, 1, 1, 0],[1, 1, 0, 1, 1, 1],[2, 0, 0, 0, 0, 0],[3, 1, 0, 1, 0, 0]],
[[1, 1, 1, 0, 1, 0],[1, 1, 1, 0, 1, 1],[1, 1, 1, 1, 0, 0],[2, 0, 0, 0, 1, 0],[3, 1, 1, 1, 1, 1]],
[[1, 1, 1, 1, 1, 1],[2, 0, 0, 1, 0, 0],[3, 0, 0, 0, 0, 0],[3, 0, 0, 0, 0, 1],[3, 1, 1, 0, 1, 0]],
[[1, 1, 1, 0, 0, 0],[2, 0, 0, 1, 0, 1],[3, 0, 0, 0, 1, 0],[3, 0, 0, 0, 1, 1],[3, 1, 1, 1, 0, 0]],
[[2, 0, 1, 0, 0, 1],[2, 0, 1, 0, 1, 0],[2, 0, 1, 0, 1, 1],[3, 0, 0, 1, 0, 0],[3, 0, 1, 1, 0, 0]],
[[2, 0, 1, 1, 1, 0],[2, 0, 1, 1, 1, 1],[2, 1, 0, 0, 0, 0],[3, 0, 0, 1, 1, 1],[3, 1, 0, 1, 1, 0]],
[[2, 1, 0, 0, 1, 1],[2, 1, 0, 1, 0, 0],[2, 1, 0, 1, 0, 1],[3, 0, 1, 0, 0, 1],[3, 1, 1, 0, 1, 1]],
[[2, 1, 1, 0, 0, 0],[2, 1, 1, 0, 0, 1],[2, 1, 1, 0, 1, 0],[3, 0, 0, 1, 0, 1],[3, 1, 1, 1, 1, 0]],
[[2, 1, 1, 1, 0, 1],[2, 1, 1, 1, 1, 0],[2, 1, 1, 1, 1, 1],[3, 0, 1, 1, 0, 1],[3, 1, 0, 0, 0, 1]],
[[2, 0, 0, 1, 1, 0],[2, 0, 0, 1, 1, 1],[2, 0, 1, 1, 0, 0],[3, 0, 0, 1, 1, 0],[3, 0, 1, 0, 1, 1]],
[[2, 1, 0, 0, 0, 1],[2, 1, 0, 1, 1, 0],[2, 1, 0, 1, 1, 1],[3, 0, 1, 0, 0, 0],[3, 1, 1, 0, 0, 0]],
[[2, 0, 1, 1, 0, 1],[2, 1, 0, 0, 1, 0],[2, 1, 1, 0, 1, 1],[3, 0, 1, 0, 1, 0],[3, 0, 1, 1, 1, 0]],
[[0, 1, 0, 1, 0, 1],[1, 1, 0, 1, 0, 0],[2, 0, 1, 0, 0, 0],[2, 1, 1, 1, 0, 0],[3, 1, 0, 1, 0, 1]],
[[0, 0, 0, 1, 1, 0],[1, 0, 0, 0, 1, 1],[1, 1, 0, 0, 1, 1],[3, 0, 1, 1, 1, 1],[3, 1, 1, 0, 0, 1]],
[[0, 1, 1, 1, 0, 1],[1, 0, 0, 1, 1, 1],[1, 1, 1, 0, 0, 1],[1, 1, 1, 1, 0, 1],[1, 1, 1, 1, 1, 0]]
]
A partition for sets of sizes: 13*3  2*4 16*5
\end{verbatim}
\newpage

\begin{verbatim} 
[4, 2, 2, 2, 2, 2]
[14, 0, 17]
[
[[0, 0, 0, 0, 0, 1],[0, 0, 0, 0, 1, 0],[0, 0, 0, 0, 1, 1]],
[[0, 0, 0, 1, 0, 0],[0, 0, 1, 0, 0, 0],[0, 0, 1, 1, 0, 0]],
[[0, 0, 0, 1, 1, 1],[0, 0, 1, 0, 0, 1],[0, 0, 1, 1, 1, 0]],
[[0, 0, 1, 0, 1, 0],[0, 1, 0, 0, 0, 0],[0, 1, 1, 0, 1, 0]],
[[0, 0, 1, 1, 0, 1],[0, 1, 0, 0, 0, 1],[0, 1, 1, 1, 0, 0]],
[[0, 0, 1, 0, 1, 1],[0, 1, 0, 0, 1, 0],[0, 1, 1, 0, 0, 1]],
[[0, 1, 0, 0, 1, 1],[1, 0, 0, 0, 0, 0],[3, 1, 0, 0, 1, 1]],
[[0, 1, 0, 1, 1, 0],[1, 0, 0, 0, 0, 1],[3, 1, 0, 1, 1, 1]],
[[0, 0, 0, 1, 0, 1],[0, 1, 1, 0, 1, 1],[0, 1, 1, 1, 1, 0]],
[[0, 0, 1, 1, 1, 1],[0, 1, 0, 1, 1, 1],[0, 1, 1, 0, 0, 0]],
[[0, 1, 1, 1, 1, 1],[1, 0, 0, 0, 1, 0],[3, 1, 1, 1, 0, 1]],
[[0, 1, 0, 1, 0, 0],[1, 0, 0, 1, 0, 0],[3, 1, 0, 0, 0, 0]],
[[1, 0, 0, 1, 0, 1],[1, 0, 0, 1, 1, 0],[2, 0, 0, 0, 1, 1]],
[[1, 0, 1, 0, 0, 0],[1, 0, 1, 0, 0, 1],[2, 0, 0, 0, 0, 1]],
[[1, 0, 1, 0, 1, 1],[1, 0, 1, 1, 0, 0],[1, 0, 1, 1, 0, 1],[2, 0, 0, 0, 0, 0],[3, 0, 1, 0, 1, 0]],
[[1, 1, 0, 0, 0, 0],[1, 1, 0, 0, 0, 1],[1, 1, 0, 0, 1, 0],[2, 0, 0, 0, 1, 0],[3, 1, 0, 0, 0, 1]],
[[1, 1, 0, 1, 0, 1],[1, 1, 0, 1, 1, 0],[1, 1, 0, 1, 1, 1],[2, 0, 0, 1, 1, 0],[3, 1, 0, 0, 1, 0]],
[[1, 1, 1, 0, 1, 0],[1, 1, 1, 0, 1, 1],[1, 1, 1, 1, 0, 0],[2, 0, 0, 1, 0, 0],[3, 1, 1, 0, 0, 1]],
[[1, 1, 1, 1, 1, 1],[2, 0, 0, 1, 0, 1],[3, 0, 0, 0, 0, 0],[3, 0, 0, 0, 0, 1],[3, 1, 1, 0, 1, 1]],
[[1, 1, 1, 1, 0, 1],[2, 0, 0, 1, 1, 1],[3, 0, 0, 0, 1, 0],[3, 0, 0, 1, 0, 0],[3, 1, 1, 1, 0, 0]],
[[2, 0, 1, 0, 0, 1],[2, 0, 1, 0, 1, 0],[2, 0, 1, 0, 1, 1],[3, 0, 0, 0, 1, 1],[3, 0, 1, 0, 1, 1]],
[[2, 0, 1, 1, 1, 0],[2, 0, 1, 1, 1, 1],[2, 1, 0, 0, 0, 0],[3, 0, 0, 1, 0, 1],[3, 1, 0, 1, 0, 0]],
[[2, 1, 0, 0, 1, 1],[2, 1, 0, 1, 0, 0],[2, 1, 0, 1, 0, 1],[3, 0, 0, 1, 1, 1],[3, 1, 0, 1, 0, 1]],
[[2, 1, 1, 0, 0, 0],[2, 1, 1, 0, 0, 1],[2, 1, 1, 0, 1, 0],[3, 0, 1, 1, 0, 1],[3, 1, 0, 1, 1, 0]],
[[2, 1, 1, 1, 0, 1],[2, 1, 1, 1, 1, 0],[2, 1, 1, 1, 1, 1],[3, 0, 0, 1, 1, 0],[3, 1, 1, 0, 1, 0]],
[[1, 1, 1, 0, 0, 0],[2, 1, 0, 0, 0, 1],[3, 0, 1, 0, 0, 0],[3, 0, 1, 1, 1, 0],[3, 0, 1, 1, 1, 1]],
[[2, 1, 0, 1, 1, 0],[2, 1, 1, 0, 1, 1],[2, 1, 1, 1, 0, 0],[3, 0, 1, 0, 0, 1],[3, 1, 1, 0, 0, 0]],
[[1, 0, 1, 1, 1, 0],[1, 1, 0, 1, 0, 0],[1, 1, 1, 1, 1, 0],[2, 0, 1, 0, 0, 0],[3, 0, 1, 1, 0, 0]],
[[1, 0, 0, 0, 1, 1],[1, 1, 0, 0, 1, 1],[1, 1, 1, 0, 0, 1],[2, 1, 0, 1, 1, 1],[3, 1, 1, 1, 1, 0]],
[[0, 0, 0, 1, 1, 0],[1, 0, 0, 1, 1, 1],[2, 0, 1, 1, 0, 0],[2, 1, 0, 0, 1, 0],[3, 1, 1, 1, 1, 1]],
[[0, 1, 0, 1, 0, 1],[0, 1, 1, 1, 0, 1],[1, 0, 1, 0, 1, 0],[1, 0, 1, 1, 1, 1],[2, 0, 1, 1, 0, 1]]
]
A partition for sets of sizes: 14*3  0*4 17*5
\end{verbatim}
\newpage

\begin{verbatim} 
[4, 2, 2, 2, 2, 2]
[10, 3, 17]
[
[[0, 0, 0, 0, 0, 1],[0, 0, 0, 0, 1, 0],[0, 0, 0, 0, 1, 1]],
[[0, 0, 0, 1, 0, 0],[0, 0, 1, 0, 0, 0],[0, 0, 1, 1, 0, 0]],
[[0, 0, 0, 1, 1, 1],[0, 0, 1, 0, 0, 1],[0, 0, 1, 1, 1, 0]],
[[0, 0, 1, 0, 1, 0],[0, 1, 0, 0, 0, 0],[0, 1, 1, 0, 1, 0]],
[[0, 0, 1, 1, 0, 1],[0, 1, 0, 0, 0, 1],[0, 1, 1, 1, 0, 0]],
[[0, 0, 1, 0, 1, 1],[0, 1, 0, 0, 1, 0],[0, 1, 1, 0, 0, 1]],
[[0, 1, 0, 0, 1, 1],[1, 0, 0, 0, 0, 0],[3, 1, 0, 0, 1, 1]],
[[0, 1, 0, 1, 1, 0],[1, 0, 0, 0, 0, 1],[3, 1, 0, 1, 1, 1]],
[[0, 0, 0, 1, 0, 1],[0, 1, 1, 0, 1, 1],[0, 1, 1, 1, 1, 0]],
[[0, 0, 1, 1, 1, 1],[0, 1, 0, 1, 1, 1],[0, 1, 1, 0, 0, 0]],
[[0, 1, 1, 1, 1, 1],[1, 0, 0, 0, 1, 0],[1, 0, 0, 0, 1, 1],[2, 1, 1, 1, 1, 0]],
[[0, 1, 1, 1, 0, 1],[1, 0, 0, 1, 0, 0],[1, 0, 0, 1, 0, 1],[2, 1, 1, 1, 0, 0]],
[[1, 0, 0, 1, 1, 1],[1, 0, 1, 0, 0, 0],[1, 1, 0, 0, 0, 0],[1, 1, 1, 1, 1, 1]],
[[1, 0, 1, 0, 1, 1],[1, 0, 1, 1, 0, 0],[1, 0, 1, 1, 0, 1],[2, 0, 0, 0, 0, 0],[3, 0, 1, 0, 1, 0]],
[[1, 0, 1, 0, 0, 1],[1, 1, 0, 0, 0, 1],[1, 1, 0, 0, 1, 0],[2, 0, 0, 0, 0, 1],[3, 0, 1, 0, 1, 1]],
[[1, 1, 0, 1, 0, 1],[1, 1, 0, 1, 1, 0],[1, 1, 0, 1, 1, 1],[2, 0, 0, 0, 1, 0],[3, 1, 0, 1, 1, 0]],
[[1, 1, 1, 0, 1, 0],[1, 1, 1, 0, 1, 1],[1, 1, 1, 1, 0, 0],[2, 0, 0, 0, 1, 1],[3, 1, 1, 1, 1, 0]],
[[1, 0, 1, 0, 1, 0],[1, 0, 1, 1, 1, 0],[1, 1, 0, 0, 1, 1],[2, 0, 0, 1, 0, 1],[3, 1, 0, 0, 1, 0]],
[[2, 0, 0, 1, 0, 0],[2, 0, 0, 1, 1, 0],[2, 0, 0, 1, 1, 1],[3, 0, 0, 0, 0, 0],[3, 0, 0, 1, 0, 1]],
[[2, 0, 1, 0, 0, 1],[2, 0, 1, 0, 1, 0],[2, 0, 1, 0, 1, 1],[3, 0, 0, 0, 0, 1],[3, 0, 1, 0, 0, 1]],
[[2, 0, 1, 1, 1, 0],[2, 0, 1, 1, 1, 1],[2, 1, 0, 0, 0, 0],[3, 0, 0, 1, 0, 0],[3, 1, 0, 1, 0, 1]],
[[2, 1, 0, 0, 1, 1],[2, 1, 0, 1, 0, 0],[2, 1, 0, 1, 0, 1],[3, 0, 0, 0, 1, 0],[3, 1, 0, 0, 0, 0]],
[[2, 1, 1, 0, 0, 0],[2, 1, 1, 0, 0, 1],[2, 1, 1, 0, 1, 0],[3, 0, 0, 0, 1, 1],[3, 1, 1, 0, 0, 0]],
[[0, 1, 0, 1, 0, 0],[1, 1, 1, 0, 0, 0],[2, 1, 1, 1, 0, 1],[2, 1, 1, 1, 1, 1],[3, 0, 1, 1, 1, 0]],
[[2, 1, 0, 1, 1, 0],[2, 1, 0, 1, 1, 1],[2, 1, 1, 0, 1, 1],[3, 0, 0, 1, 1, 0],[3, 1, 1, 1, 0, 0]],
[[2, 0, 1, 1, 0, 1],[2, 1, 0, 0, 0, 1],[2, 1, 0, 0, 1, 0],[3, 1, 0, 0, 0, 1],[3, 1, 1, 1, 1, 1]],
[[1, 0, 1, 1, 1, 1],[2, 0, 1, 1, 0, 0],[3, 0, 1, 1, 0, 1],[3, 1, 0, 1, 0, 0],[3, 1, 1, 0, 1, 0]],
[[1, 1, 1, 1, 0, 1],[2, 0, 1, 0, 0, 0],[3, 0, 0, 1, 1, 1],[3, 0, 1, 1, 1, 1],[3, 1, 1, 1, 0, 1]],
[[0, 0, 0, 1, 1, 0],[3, 0, 1, 0, 0, 0],[3, 0, 1, 1, 0, 0],[3, 1, 1, 0, 0, 1],[3, 1, 1, 0, 1, 1]],
[[0, 1, 0, 1, 0, 1],[1, 0, 0, 1, 1, 0],[1, 1, 0, 1, 0, 0],[1, 1, 1, 0, 0, 1],[1, 1, 1, 1, 1, 0]]
]
A partition for sets of sizes: 10*3  3*4 17*5
\end{verbatim}
\newpage

\begin{verbatim} 
[4, 2, 2, 2, 2, 2]
[11, 1, 18]
[
[[0, 0, 0, 0, 0, 1],[0, 0, 0, 0, 1, 0],[0, 0, 0, 0, 1, 1]],
[[0, 0, 0, 1, 0, 0],[0, 0, 1, 0, 0, 0],[0, 0, 1, 1, 0, 0]],
[[0, 0, 0, 1, 1, 1],[0, 0, 1, 0, 0, 1],[0, 0, 1, 1, 1, 0]],
[[0, 0, 1, 0, 1, 0],[0, 1, 0, 0, 0, 0],[0, 1, 1, 0, 1, 0]],
[[0, 0, 1, 1, 0, 1],[0, 1, 0, 0, 0, 1],[0, 1, 1, 1, 0, 0]],
[[0, 0, 1, 0, 1, 1],[0, 1, 0, 0, 1, 0],[0, 1, 1, 0, 0, 1]],
[[0, 1, 0, 0, 1, 1],[1, 0, 0, 0, 0, 0],[3, 1, 0, 0, 1, 1]],
[[0, 1, 0, 1, 1, 0],[1, 0, 0, 0, 0, 1],[3, 1, 0, 1, 1, 1]],
[[0, 0, 0, 1, 0, 1],[0, 1, 1, 0, 1, 1],[0, 1, 1, 1, 1, 0]],
[[0, 0, 1, 1, 1, 1],[0, 1, 0, 1, 1, 1],[0, 1, 1, 0, 0, 0]],
[[0, 1, 1, 1, 1, 1],[1, 0, 0, 0, 1, 0],[3, 1, 1, 1, 0, 1]],
[[0, 1, 0, 1, 0, 0],[1, 0, 0, 0, 1, 1],[1, 0, 0, 1, 0, 0],[2, 1, 0, 0, 1, 1]],
[[1, 0, 0, 1, 1, 0],[1, 0, 0, 1, 1, 1],[1, 0, 1, 0, 0, 0],[2, 0, 0, 0, 0, 0],[3, 0, 1, 0, 0, 1]],
[[1, 0, 1, 0, 1, 1],[1, 0, 1, 1, 0, 0],[1, 0, 1, 1, 0, 1],[2, 0, 0, 0, 0, 1],[3, 0, 1, 0, 1, 1]],
[[1, 1, 0, 0, 0, 0],[1, 1, 0, 0, 0, 1],[1, 1, 0, 0, 1, 0],[2, 0, 0, 0, 1, 0],[3, 1, 0, 0, 0, 1]],
[[1, 1, 0, 1, 0, 1],[1, 1, 0, 1, 1, 0],[1, 1, 0, 1, 1, 1],[2, 0, 0, 1, 0, 0],[3, 1, 0, 0, 0, 0]],
[[1, 1, 1, 0, 1, 0],[1, 1, 1, 0, 1, 1],[1, 1, 1, 1, 0, 0],[2, 0, 0, 0, 1, 1],[3, 1, 1, 1, 1, 0]],
[[1, 1, 1, 1, 1, 1],[2, 0, 0, 1, 0, 1],[3, 0, 0, 0, 0, 0],[3, 0, 0, 0, 0, 1],[3, 1, 1, 0, 1, 1]],
[[1, 1, 1, 0, 0, 0],[2, 0, 0, 1, 1, 0],[3, 0, 0, 0, 1, 0],[3, 0, 0, 0, 1, 1],[3, 1, 1, 1, 1, 1]],
[[2, 0, 1, 0, 0, 1],[2, 0, 1, 0, 1, 0],[2, 0, 1, 0, 1, 1],[3, 0, 0, 1, 0, 0],[3, 0, 1, 1, 0, 0]],
[[2, 0, 1, 1, 1, 0],[2, 0, 1, 1, 1, 1],[2, 1, 0, 0, 0, 0],[3, 0, 0, 1, 0, 1],[3, 1, 0, 1, 0, 0]],
[[1, 0, 0, 1, 0, 1],[2, 1, 0, 1, 0, 0],[3, 0, 0, 1, 1, 0],[3, 0, 1, 1, 0, 1],[3, 1, 1, 0, 1, 0]],
[[2, 1, 1, 0, 0, 0],[2, 1, 1, 0, 0, 1],[2, 1, 1, 0, 1, 0],[3, 0, 0, 1, 1, 1],[3, 1, 1, 1, 0, 0]],
[[2, 1, 1, 1, 0, 1],[2, 1, 1, 1, 1, 0],[2, 1, 1, 1, 1, 1],[3, 0, 1, 0, 1, 0],[3, 1, 0, 1, 1, 0]],
[[1, 0, 1, 0, 0, 1],[2, 0, 0, 1, 1, 1],[3, 0, 1, 1, 1, 1],[3, 1, 1, 0, 0, 0],[3, 1, 1, 0, 0, 1]],
[[1, 0, 1, 0, 1, 0],[1, 0, 1, 1, 1, 0],[2, 0, 1, 1, 0, 1],[2, 1, 0, 0, 1, 0],[2, 1, 1, 0, 1, 1]],
[[1, 0, 1, 1, 1, 1],[1, 1, 0, 1, 0, 0],[1, 1, 1, 0, 0, 1],[2, 0, 1, 1, 0, 0],[3, 0, 1, 1, 1, 0]],
[[2, 1, 0, 1, 1, 0],[2, 1, 0, 1, 1, 1],[2, 1, 1, 1, 0, 0],[3, 0, 1, 0, 0, 0],[3, 1, 0, 1, 0, 1]],
[[0, 0, 0, 1, 1, 0],[0, 1, 1, 1, 0, 1],[1, 1, 0, 0, 1, 1],[1, 1, 1, 1, 0, 1],[2, 1, 0, 1, 0, 1]],
[[0, 1, 0, 1, 0, 1],[1, 1, 1, 1, 1, 0],[2, 0, 1, 0, 0, 0],[2, 1, 0, 0, 0, 1],[3, 1, 0, 0, 1, 0]]
]
A partition for sets of sizes: 11*3  1*4 18*5
\end{verbatim}
\newpage

\begin{verbatim} 
[4, 2, 2, 2, 2, 2]
[8, 2, 19]
[
[[0, 0, 0, 0, 0, 1],[0, 0, 0, 0, 1, 0],[0, 0, 0, 0, 1, 1]],
[[0, 0, 0, 1, 0, 0],[0, 0, 1, 0, 0, 0],[0, 0, 1, 1, 0, 0]],
[[0, 0, 0, 1, 1, 1],[0, 0, 1, 0, 0, 1],[0, 0, 1, 1, 1, 0]],
[[0, 0, 1, 0, 1, 0],[0, 1, 0, 0, 0, 0],[0, 1, 1, 0, 1, 0]],
[[0, 0, 1, 1, 0, 1],[0, 1, 0, 0, 0, 1],[0, 1, 1, 1, 0, 0]],
[[0, 0, 1, 0, 1, 1],[0, 1, 0, 0, 1, 0],[0, 1, 1, 0, 0, 1]],
[[0, 1, 0, 0, 1, 1],[1, 0, 0, 0, 0, 0],[3, 1, 0, 0, 1, 1]],
[[0, 1, 0, 1, 1, 0],[1, 0, 0, 0, 0, 1],[3, 1, 0, 1, 1, 1]],
[[0, 0, 0, 1, 0, 1],[0, 0, 0, 1, 1, 0],[0, 1, 1, 1, 0, 1],[0, 1, 1, 1, 1, 0]],
[[0, 1, 1, 0, 1, 1],[0, 1, 1, 1, 1, 1],[1, 0, 0, 0, 1, 0],[3, 0, 0, 1, 1, 0]],
[[0, 1, 0, 1, 1, 1],[1, 0, 0, 0, 1, 1],[1, 0, 0, 1, 0, 0],[1, 0, 0, 1, 0, 1],[1, 1, 0, 1, 0, 1]],
[[1, 0, 0, 1, 1, 0],[1, 0, 0, 1, 1, 1],[1, 0, 1, 0, 0, 0],[2, 0, 0, 0, 0, 0],[3, 0, 1, 0, 0, 1]],
[[1, 0, 1, 0, 1, 1],[1, 0, 1, 1, 0, 0],[1, 0, 1, 1, 0, 1],[2, 0, 0, 0, 0, 1],[3, 0, 1, 0, 1, 1]],
[[1, 1, 0, 0, 0, 0],[1, 1, 0, 0, 0, 1],[1, 1, 0, 0, 1, 0],[2, 0, 0, 0, 1, 0],[3, 1, 0, 0, 0, 1]],
[[0, 0, 1, 1, 1, 1],[1, 1, 0, 1, 1, 0],[1, 1, 0, 1, 1, 1],[3, 0, 0, 0, 0, 0],[3, 0, 1, 1, 1, 0]],
[[1, 1, 1, 0, 1, 0],[1, 1, 1, 0, 1, 1],[1, 1, 1, 1, 0, 0],[2, 0, 0, 0, 1, 1],[3, 1, 1, 1, 1, 0]],
[[1, 1, 1, 1, 1, 1],[2, 0, 0, 1, 0, 0],[3, 0, 0, 0, 0, 1],[3, 0, 0, 0, 1, 0],[3, 1, 1, 0, 0, 0]],
[[1, 0, 1, 0, 0, 1],[2, 0, 0, 1, 0, 1],[3, 0, 0, 0, 1, 1],[3, 0, 0, 1, 0, 1],[3, 0, 1, 0, 1, 0]],
[[2, 0, 1, 0, 0, 1],[2, 0, 1, 0, 1, 0],[2, 0, 1, 0, 1, 1],[3, 0, 0, 1, 0, 0],[3, 0, 1, 1, 0, 0]],
[[2, 0, 1, 1, 1, 0],[2, 0, 1, 1, 1, 1],[2, 1, 0, 0, 0, 0],[3, 0, 0, 1, 1, 1],[3, 1, 0, 1, 1, 0]],
[[2, 1, 0, 0, 1, 1],[2, 1, 0, 1, 0, 0],[2, 1, 0, 1, 0, 1],[3, 0, 1, 0, 0, 0],[3, 1, 1, 0, 1, 0]],
[[2, 1, 1, 0, 0, 0],[2, 1, 1, 0, 0, 1],[2, 1, 1, 0, 1, 0],[3, 0, 1, 1, 1, 1],[3, 1, 0, 1, 0, 0]],
[[1, 1, 1, 0, 0, 0],[2, 1, 1, 1, 1, 0],[3, 0, 1, 1, 0, 1],[3, 1, 0, 0, 0, 0],[3, 1, 1, 0, 1, 1]],
[[1, 1, 0, 0, 1, 1],[2, 0, 0, 1, 1, 1],[3, 1, 0, 0, 1, 0],[3, 1, 1, 0, 0, 1],[3, 1, 1, 1, 1, 1]],
[[1, 0, 1, 0, 1, 0],[1, 0, 1, 1, 1, 1],[2, 0, 1, 0, 0, 0],[2, 1, 0, 0, 0, 1],[2, 1, 1, 1, 0, 0]],
[[2, 0, 1, 1, 0, 1],[2, 1, 1, 0, 1, 1],[2, 1, 1, 1, 1, 1],[3, 1, 0, 1, 0, 1],[3, 1, 1, 1, 0, 0]],
[[1, 1, 0, 1, 0, 0],[1, 1, 1, 0, 0, 1],[2, 0, 1, 1, 0, 0],[2, 1, 0, 1, 1, 0],[2, 1, 0, 1, 1, 1]],
[[0, 1, 1, 0, 0, 0],[1, 1, 1, 1, 1, 0],[2, 0, 0, 1, 1, 0],[2, 1, 1, 1, 0, 1],[3, 1, 1, 1, 0, 1]],
[[0, 1, 0, 1, 0, 0],[0, 1, 0, 1, 0, 1],[1, 0, 1, 1, 1, 0],[1, 1, 1, 1, 0, 1],[2, 1, 0, 0, 1, 0]]
]
A partition for sets of sizes:  8*3  2*4 19*5
\end{verbatim}
\newpage

\begin{verbatim} 
[4, 2, 2, 2, 2, 2]
[9, 0, 20]
[
[[0, 0, 0, 0, 0, 1],[0, 0, 0, 0, 1, 0],[0, 0, 0, 0, 1, 1]],
[[0, 0, 0, 1, 0, 0],[0, 0, 1, 0, 0, 0],[0, 0, 1, 1, 0, 0]],
[[0, 0, 0, 1, 1, 1],[0, 0, 1, 0, 0, 1],[0, 0, 1, 1, 1, 0]],
[[0, 0, 1, 0, 1, 0],[0, 1, 0, 0, 0, 0],[0, 1, 1, 0, 1, 0]],
[[0, 0, 1, 1, 0, 1],[0, 1, 0, 0, 0, 1],[0, 1, 1, 1, 0, 0]],
[[0, 0, 1, 0, 1, 1],[0, 1, 0, 0, 1, 0],[0, 1, 1, 0, 0, 1]],
[[0, 1, 0, 0, 1, 1],[1, 0, 0, 0, 0, 0],[3, 1, 0, 0, 1, 1]],
[[0, 1, 0, 1, 1, 0],[1, 0, 0, 0, 0, 1],[3, 1, 0, 1, 1, 1]],
[[0, 0, 0, 1, 0, 1],[0, 1, 1, 0, 1, 1],[0, 1, 1, 1, 1, 0]],
[[0, 0, 1, 1, 1, 1],[0, 1, 1, 1, 0, 1],[0, 1, 1, 1, 1, 1],[1, 0, 0, 0, 1, 0],[3, 0, 1, 1, 1, 1]],
[[0, 1, 0, 1, 1, 1],[1, 0, 0, 0, 1, 1],[1, 0, 0, 1, 0, 0],[1, 0, 0, 1, 0, 1],[1, 1, 0, 1, 0, 1]],
[[1, 0, 0, 1, 1, 0],[1, 0, 0, 1, 1, 1],[1, 0, 1, 0, 0, 0],[2, 0, 0, 0, 0, 0],[3, 0, 1, 0, 0, 1]],
[[1, 0, 1, 0, 1, 1],[1, 0, 1, 1, 0, 0],[1, 0, 1, 1, 0, 1],[2, 0, 0, 0, 0, 1],[3, 0, 1, 0, 1, 1]],
[[1, 1, 0, 0, 0, 0],[1, 1, 0, 0, 0, 1],[1, 1, 0, 0, 1, 0],[2, 0, 0, 0, 1, 0],[3, 1, 0, 0, 0, 1]],
[[0, 0, 0, 1, 1, 0],[1, 1, 0, 1, 1, 0],[1, 1, 0, 1, 1, 1],[1, 1, 1, 0, 0, 0],[1, 1, 1, 1, 1, 1]],
[[1, 1, 1, 0, 1, 0],[1, 1, 1, 0, 1, 1],[1, 1, 1, 1, 0, 0],[2, 0, 0, 0, 1, 1],[3, 1, 1, 1, 1, 0]],
[[1, 1, 1, 0, 0, 1],[1, 1, 1, 1, 0, 1],[2, 0, 0, 1, 0, 1],[2, 0, 0, 1, 1, 0],[2, 0, 0, 1, 1, 1]],
[[2, 0, 0, 1, 0, 0],[2, 0, 1, 0, 0, 0],[2, 0, 1, 0, 0, 1],[3, 0, 0, 0, 0, 0],[3, 0, 0, 1, 0, 1]],
[[1, 1, 0, 0, 1, 1],[2, 0, 1, 0, 1, 0],[3, 0, 0, 0, 0, 1],[3, 0, 0, 0, 1, 0],[3, 1, 1, 0, 1, 0]],
[[2, 0, 1, 1, 1, 0],[2, 0, 1, 1, 1, 1],[2, 1, 0, 0, 0, 0],[3, 0, 0, 0, 1, 1],[3, 1, 0, 0, 1, 0]],
[[2, 1, 0, 0, 1, 1],[2, 1, 0, 1, 0, 0],[2, 1, 0, 1, 0, 1],[3, 0, 0, 1, 0, 0],[3, 1, 0, 1, 1, 0]],
[[2, 1, 1, 0, 0, 0],[2, 1, 1, 0, 0, 1],[2, 1, 1, 0, 1, 0],[3, 0, 0, 1, 1, 0],[3, 1, 1, 1, 0, 1]],
[[2, 1, 1, 1, 0, 1],[2, 1, 1, 1, 1, 0],[2, 1, 1, 1, 1, 1],[3, 0, 0, 1, 1, 1],[3, 1, 1, 0, 1, 1]],
[[2, 0, 1, 1, 0, 0],[2, 1, 0, 0, 0, 1],[2, 1, 0, 1, 1, 0],[3, 1, 0, 1, 0, 0],[3, 1, 1, 1, 1, 1]],
[[1, 0, 1, 0, 0, 1],[1, 0, 1, 0, 1, 0],[1, 0, 1, 1, 1, 1],[2, 1, 1, 1, 0, 0],[3, 1, 0, 0, 0, 0]],
[[0, 1, 0, 1, 0, 0],[3, 0, 1, 0, 0, 0],[3, 0, 1, 0, 1, 0],[3, 0, 1, 1, 1, 0],[3, 1, 1, 0, 0, 0]],
[[0, 1, 0, 1, 0, 1],[1, 1, 0, 1, 0, 0],[2, 1, 0, 1, 1, 1],[2, 1, 1, 0, 1, 1],[3, 0, 1, 1, 0, 1]],
[[0, 1, 1, 0, 0, 0],[1, 1, 1, 1, 1, 0],[2, 0, 1, 1, 0, 1],[2, 1, 0, 0, 1, 0],[3, 1, 1, 0, 0, 1]],
[[1, 0, 1, 1, 1, 0],[2, 0, 1, 0, 1, 1],[3, 0, 1, 1, 0, 0],[3, 1, 0, 1, 0, 1],[3, 1, 1, 1, 0, 0]]
]
A partition for sets of sizes:  9*3  0*4 20*5
\end{verbatim}

\newpage

\begin{verbatim} 
[4, 2, 2, 2, 2, 2]
[5, 3, 20]
[
[[0, 0, 0, 0, 0, 1],[0, 0, 0, 0, 1, 0],[0, 0, 0, 0, 1, 1]],
[[0, 0, 0, 1, 0, 0],[0, 0, 1, 0, 0, 0],[0, 0, 1, 1, 0, 0]],
[[0, 0, 0, 1, 1, 1],[0, 0, 1, 0, 0, 1],[0, 0, 1, 1, 1, 0]],
[[0, 0, 1, 0, 1, 0],[0, 1, 0, 0, 0, 0],[0, 1, 1, 0, 1, 0]],
[[0, 0, 1, 1, 0, 1],[0, 1, 0, 0, 0, 1],[0, 1, 1, 1, 0, 0]],
[[0, 0, 1, 0, 1, 1],[0, 1, 0, 0, 1, 0],[1, 0, 0, 0, 0, 0],[3, 1, 1, 0, 0, 1]],
[[0, 1, 0, 1, 0, 0],[0, 1, 0, 1, 0, 1],[0, 1, 0, 1, 1, 0],[0, 1, 0, 1, 1, 1]],
[[0, 1, 1, 0, 0, 0],[0, 1, 1, 0, 0, 1],[0, 1, 1, 1, 1, 0],[0, 1, 1, 1, 1, 1]],
[[0, 0, 1, 1, 1, 1],[0, 1, 1, 1, 0, 1],[1, 0, 0, 0, 0, 1],[1, 0, 0, 0, 1, 0],[2, 1, 0, 0, 0, 1]],
[[0, 0, 0, 1, 1, 0],[0, 1, 1, 0, 1, 1],[1, 0, 0, 0, 1, 1],[1, 0, 0, 1, 0, 0],[2, 1, 1, 0, 1, 0]],
[[1, 0, 0, 1, 1, 0],[1, 0, 0, 1, 1, 1],[1, 0, 1, 0, 0, 0],[2, 0, 0, 0, 0, 0],[3, 0, 1, 0, 0, 1]],
[[1, 0, 1, 0, 1, 1],[1, 0, 1, 1, 0, 0],[1, 0, 1, 1, 0, 1],[2, 0, 0, 0, 0, 1],[3, 0, 1, 0, 1, 1]],
[[1, 1, 0, 0, 0, 0],[1, 1, 0, 0, 0, 1],[1, 1, 0, 0, 1, 0],[2, 0, 0, 0, 1, 0],[3, 1, 0, 0, 0, 1]],
[[1, 1, 0, 1, 0, 1],[1, 1, 0, 1, 1, 0],[1, 1, 0, 1, 1, 1],[2, 0, 0, 0, 1, 1],[3, 1, 0, 1, 1, 1]],
[[1, 1, 1, 0, 1, 0],[1, 1, 1, 0, 1, 1],[1, 1, 1, 1, 0, 0],[2, 0, 0, 1, 0, 1],[3, 1, 1, 0, 0, 0]],
[[1, 1, 1, 1, 1, 1],[2, 0, 0, 1, 0, 0],[3, 0, 0, 0, 0, 0],[3, 0, 0, 0, 0, 1],[3, 1, 1, 0, 1, 0]],
[[1, 0, 1, 0, 0, 1],[1, 1, 1, 1, 0, 1],[2, 0, 0, 1, 1, 0],[2, 0, 0, 1, 1, 1],[2, 1, 0, 1, 0, 1]],
[[2, 0, 1, 0, 0, 1],[2, 0, 1, 0, 1, 0],[2, 0, 1, 0, 1, 1],[3, 0, 0, 0, 1, 0],[3, 0, 1, 0, 1, 0]],
[[2, 0, 1, 1, 1, 0],[2, 0, 1, 1, 1, 1],[2, 1, 0, 0, 0, 0],[3, 0, 0, 0, 1, 1],[3, 1, 0, 0, 1, 0]],
[[2, 1, 0, 0, 1, 1],[2, 1, 0, 1, 0, 0],[2, 1, 0, 1, 1, 0],[3, 0, 0, 1, 0, 0],[3, 1, 0, 1, 0, 1]],
[[2, 1, 1, 0, 0, 0],[2, 1, 1, 0, 0, 1],[2, 1, 1, 0, 1, 1],[3, 0, 0, 1, 0, 1],[3, 1, 1, 1, 1, 1]],
[[2, 1, 1, 1, 0, 1],[2, 1, 1, 1, 1, 0],[2, 1, 1, 1, 1, 1],[3, 0, 0, 1, 1, 1],[3, 1, 1, 0, 1, 1]],
[[2, 0, 1, 1, 0, 0],[2, 0, 1, 1, 0, 1],[2, 1, 0, 0, 1, 0],[3, 0, 1, 1, 0, 1],[3, 1, 1, 1, 1, 0]],
[[1, 0, 1, 1, 1, 0],[1, 0, 1, 1, 1, 1],[1, 1, 0, 1, 0, 0],[2, 0, 1, 0, 0, 0],[3, 1, 1, 1, 0, 1]],
[[1, 0, 0, 1, 0, 1],[1, 0, 1, 0, 1, 0],[1, 1, 1, 1, 1, 0],[2, 1, 0, 1, 1, 1],[3, 0, 0, 1, 1, 0]],
[[0, 0, 0, 1, 0, 1],[3, 0, 1, 1, 0, 0],[3, 0, 1, 1, 1, 1],[3, 1, 0, 0, 0, 0],[3, 1, 0, 1, 1, 0]],
[[0, 1, 0, 0, 1, 1],[1, 1, 1, 0, 0, 0],[1, 1, 1, 0, 0, 1],[3, 0, 1, 1, 1, 0],[3, 1, 1, 1, 0, 0]],
[[1, 1, 0, 0, 1, 1],[2, 1, 1, 1, 0, 0],[3, 0, 1, 0, 0, 0],[3, 1, 0, 0, 1, 1],[3, 1, 0, 1, 0, 0]]
]
A partition for sets of sizes:  5*3  3*4 20*5
\end{verbatim}
\newpage

\begin{verbatim} 
[4, 2, 2, 2, 2, 2]
[6, 1, 21]
[
[[0, 0, 0, 0, 0, 1],[0, 0, 0, 0, 1, 0],[0, 0, 0, 0, 1, 1]],
[[0, 0, 0, 1, 0, 0],[0, 0, 1, 0, 0, 0],[0, 0, 1, 1, 0, 0]],
[[0, 0, 0, 1, 1, 1],[0, 0, 1, 0, 0, 1],[0, 0, 1, 1, 1, 0]],
[[0, 0, 1, 0, 1, 0],[0, 1, 0, 0, 0, 0],[0, 1, 1, 0, 1, 0]],
[[0, 0, 1, 1, 0, 1],[0, 1, 0, 0, 0, 1],[0, 1, 1, 1, 0, 0]],
[[0, 0, 1, 0, 1, 1],[0, 1, 0, 0, 1, 0],[0, 1, 1, 0, 0, 1]],
[[0, 1, 0, 0, 1, 1],[0, 1, 0, 1, 0, 0],[0, 1, 1, 0, 0, 0],[0, 1, 1, 1, 1, 1]],
[[0, 1, 0, 1, 1, 1],[0, 1, 1, 0, 1, 1],[0, 1, 1, 1, 0, 1],[1, 0, 0, 0, 0, 0],[3, 1, 0, 0, 0, 1]],
[[0, 0, 1, 1, 1, 1],[0, 1, 1, 1, 1, 0],[1, 0, 0, 0, 0, 1],[1, 0, 0, 0, 1, 0],[2, 1, 0, 0, 1, 0]],
[[0, 0, 0, 1, 0, 1],[0, 1, 0, 1, 1, 0],[1, 0, 0, 0, 1, 1],[1, 0, 0, 1, 0, 0],[2, 1, 0, 1, 0, 0]],
[[1, 0, 0, 1, 1, 0],[1, 0, 0, 1, 1, 1],[1, 0, 1, 0, 0, 0],[2, 0, 0, 0, 0, 0],[3, 0, 1, 0, 0, 1]],
[[1, 0, 1, 0, 1, 1],[1, 0, 1, 1, 0, 0],[1, 0, 1, 1, 0, 1],[2, 0, 0, 0, 0, 1],[3, 0, 1, 0, 1, 1]],
[[1, 1, 0, 0, 0, 0],[1, 1, 0, 0, 0, 1],[1, 1, 0, 0, 1, 0],[2, 0, 0, 0, 1, 1],[3, 1, 0, 0, 0, 0]],
[[1, 1, 0, 1, 0, 1],[1, 1, 0, 1, 1, 0],[1, 1, 0, 1, 1, 1],[2, 0, 0, 0, 1, 0],[3, 1, 0, 1, 1, 0]],
[[1, 1, 1, 0, 1, 0],[1, 1, 1, 0, 1, 1],[1, 1, 1, 1, 0, 0],[2, 0, 0, 1, 0, 0],[3, 1, 1, 0, 0, 1]],
[[1, 1, 1, 1, 1, 1],[2, 0, 0, 1, 0, 1],[3, 0, 0, 0, 0, 0],[3, 0, 0, 0, 0, 1],[3, 1, 1, 0, 1, 1]],
[[1, 1, 1, 1, 0, 1],[2, 0, 0, 1, 1, 0],[3, 0, 0, 0, 1, 0],[3, 0, 0, 0, 1, 1],[3, 1, 1, 0, 1, 0]],
[[2, 0, 1, 0, 0, 1],[2, 0, 1, 0, 1, 0],[2, 0, 1, 0, 1, 1],[3, 0, 0, 1, 0, 0],[3, 0, 1, 1, 0, 0]],
[[2, 0, 1, 1, 1, 0],[2, 0, 1, 1, 1, 1],[2, 1, 0, 0, 0, 0],[3, 0, 0, 1, 0, 1],[3, 1, 0, 1, 0, 0]],
[[2, 1, 0, 0, 1, 1],[2, 1, 0, 1, 0, 1],[2, 1, 0, 1, 1, 0],[3, 0, 0, 1, 1, 1],[3, 1, 0, 1, 1, 1]],
[[2, 1, 1, 0, 0, 0],[2, 1, 1, 0, 0, 1],[2, 1, 1, 0, 1, 0],[3, 0, 0, 1, 1, 0],[3, 1, 1, 1, 0, 1]],
[[2, 1, 1, 1, 0, 1],[2, 1, 1, 1, 1, 0],[2, 1, 1, 1, 1, 1],[3, 0, 1, 1, 1, 0],[3, 1, 0, 0, 1, 0]],
[[1, 0, 1, 0, 0, 1],[2, 0, 0, 1, 1, 1],[3, 0, 1, 0, 0, 0],[3, 1, 0, 0, 1, 1],[3, 1, 0, 1, 0, 1]],
[[1, 0, 0, 1, 0, 1],[2, 0, 1, 1, 0, 0],[3, 0, 1, 0, 1, 0],[3, 1, 1, 1, 0, 0],[3, 1, 1, 1, 1, 1]],
[[2, 0, 1, 1, 0, 1],[2, 1, 0, 1, 1, 1],[2, 1, 1, 1, 0, 0],[3, 1, 1, 0, 0, 0],[3, 1, 1, 1, 1, 0]],
[[0, 1, 0, 1, 0, 1],[1, 0, 1, 0, 1, 0],[1, 1, 0, 0, 1, 1],[1, 1, 0, 1, 0, 0],[1, 1, 1, 0, 0, 0]],
[[2, 0, 1, 0, 0, 0],[2, 1, 0, 0, 0, 1],[2, 1, 1, 0, 1, 1],[3, 0, 1, 1, 0, 1],[3, 0, 1, 1, 1, 1]],
[[0, 0, 0, 1, 1, 0],[1, 0, 1, 1, 1, 0],[1, 0, 1, 1, 1, 1],[1, 1, 1, 0, 0, 1],[1, 1, 1, 1, 1, 0]]
]
A partition for sets of sizes:  6*3  1*4 21*5
\end{verbatim}
\newpage

\begin{verbatim} 
[4, 2, 2, 2, 2, 2]
[3, 2, 22]
[
[[0, 0, 0, 0, 0, 1],[0, 0, 0, 0, 1, 0],[0, 0, 0, 0, 1, 1]],
[[0, 0, 0, 1, 0, 0],[0, 0, 1, 0, 0, 0],[0, 0, 1, 1, 0, 0]],
[[0, 0, 0, 1, 1, 1],[0, 0, 1, 0, 0, 1],[0, 0, 1, 1, 1, 0]],
[[0, 0, 1, 0, 1, 0],[0, 0, 1, 0, 1, 1],[0, 1, 0, 0, 0, 0],[0, 1, 0, 0, 0, 1]],
[[0, 0, 0, 1, 0, 1],[0, 0, 1, 1, 1, 1],[0, 1, 0, 0, 1, 0],[0, 1, 1, 0, 0, 0]],
[[0, 0, 0, 1, 1, 0],[0, 1, 0, 0, 1, 1],[0, 1, 0, 1, 0, 0],[0, 1, 0, 1, 1, 0],[0, 1, 0, 1, 1, 1]],
[[0, 1, 0, 1, 0, 1],[0, 1, 1, 0, 0, 1],[0, 1, 1, 0, 1, 0],[1, 0, 0, 0, 0, 0],[3, 1, 0, 1, 1, 0]],
[[0, 1, 1, 1, 0, 0],[0, 1, 1, 1, 0, 1],[0, 1, 1, 1, 1, 0],[1, 0, 0, 0, 0, 1],[3, 1, 1, 1, 1, 0]],
[[0, 1, 1, 1, 1, 1],[1, 0, 0, 0, 1, 0],[1, 0, 0, 0, 1, 1],[1, 0, 0, 1, 0, 0],[1, 1, 1, 0, 1, 0]],
[[1, 0, 0, 1, 1, 0],[1, 0, 0, 1, 1, 1],[1, 0, 1, 0, 0, 0],[2, 0, 0, 0, 0, 0],[3, 0, 1, 0, 0, 1]],
[[1, 0, 1, 0, 1, 1],[1, 0, 1, 1, 0, 0],[1, 0, 1, 1, 0, 1],[2, 0, 0, 0, 0, 1],[3, 0, 1, 0, 1, 1]],
[[1, 1, 0, 0, 0, 0],[1, 1, 0, 0, 0, 1],[1, 1, 0, 0, 1, 0],[2, 0, 0, 0, 1, 0],[3, 1, 0, 0, 0, 1]],
[[1, 1, 0, 1, 0, 1],[1, 1, 0, 1, 1, 0],[1, 1, 0, 1, 1, 1],[2, 0, 0, 0, 1, 1],[3, 1, 0, 1, 1, 1]],
[[1, 0, 0, 1, 0, 1],[1, 1, 1, 0, 1, 1],[1, 1, 1, 1, 0, 0],[2, 0, 0, 1, 0, 0],[3, 0, 0, 1, 1, 0]],
[[1, 1, 1, 1, 1, 1],[2, 0, 0, 1, 0, 1],[3, 0, 0, 0, 0, 0],[3, 0, 0, 0, 0, 1],[3, 1, 1, 0, 1, 1]],
[[1, 1, 1, 1, 0, 1],[2, 0, 0, 1, 1, 0],[3, 0, 0, 0, 1, 0],[3, 0, 0, 0, 1, 1],[3, 1, 1, 0, 1, 0]],
[[2, 0, 1, 0, 0, 1],[2, 0, 1, 0, 1, 0],[2, 0, 1, 0, 1, 1],[3, 0, 0, 1, 0, 0],[3, 0, 1, 1, 0, 0]],
[[2, 0, 1, 1, 1, 0],[2, 0, 1, 1, 1, 1],[2, 1, 0, 0, 0, 0],[3, 0, 0, 1, 0, 1],[3, 1, 0, 1, 0, 0]],
[[2, 1, 0, 0, 1, 1],[2, 1, 0, 1, 0, 0],[2, 1, 0, 1, 0, 1],[3, 0, 0, 1, 1, 1],[3, 1, 0, 1, 0, 1]],
[[2, 1, 1, 0, 0, 0],[2, 1, 1, 0, 0, 1],[2, 1, 1, 0, 1, 0],[3, 0, 1, 0, 0, 0],[3, 1, 0, 0, 1, 1]],
[[2, 1, 1, 1, 0, 1],[2, 1, 1, 1, 1, 0],[2, 1, 1, 1, 1, 1],[3, 0, 1, 1, 1, 0],[3, 1, 0, 0, 1, 0]],
[[1, 0, 1, 0, 0, 1],[2, 0, 0, 1, 1, 1],[3, 0, 1, 0, 1, 0],[3, 1, 1, 0, 0, 0],[3, 1, 1, 1, 0, 0]],
[[2, 1, 0, 1, 1, 0],[2, 1, 0, 1, 1, 1],[2, 1, 1, 1, 0, 0],[3, 0, 1, 1, 0, 1],[3, 1, 0, 0, 0, 0]],
[[2, 0, 1, 1, 0, 1],[2, 1, 0, 0, 1, 0],[2, 1, 1, 0, 1, 1],[3, 1, 1, 0, 0, 1],[3, 1, 1, 1, 0, 1]],
[[1, 0, 1, 0, 1, 0],[1, 1, 0, 0, 1, 1],[1, 1, 1, 1, 1, 0],[2, 0, 1, 0, 0, 0],[3, 0, 1, 1, 1, 1]],
[[0, 1, 1, 0, 1, 1],[1, 1, 1, 0, 0, 1],[2, 0, 1, 1, 0, 0],[2, 1, 0, 0, 0, 1],[3, 1, 1, 1, 1, 1]],
[[0, 0, 1, 1, 0, 1],[1, 0, 1, 1, 1, 0],[1, 0, 1, 1, 1, 1],[1, 1, 0, 1, 0, 0],[1, 1, 1, 0, 0, 0]]
]
A partition for sets of sizes:  3*3  2*4 22*5

[4, 2, 2, 2, 2, 2]
[4, 0, 23]
[
[[0, 0, 0, 0, 0, 1],[0, 0, 0, 0, 1, 0],[0, 0, 0, 0, 1, 1]],
[[0, 0, 0, 1, 0, 0],[0, 0, 1, 0, 0, 0],[0, 0, 1, 1, 0, 0]],
[[0, 0, 0, 1, 1, 1],[0, 0, 1, 0, 0, 1],[0, 0, 1, 1, 1, 0]],
[[0, 0, 1, 0, 1, 0],[0, 1, 0, 0, 0, 0],[0, 1, 1, 0, 1, 0]],
[[0, 0, 1, 1, 0, 1],[0, 0, 1, 1, 1, 1],[0, 1, 0, 0, 0, 1],[1, 0, 0, 0, 0, 0],[3, 1, 0, 0, 1, 1]],
[[0, 1, 0, 0, 1, 0],[0, 1, 0, 0, 1, 1],[0, 1, 0, 1, 0, 0],[1, 0, 0, 0, 0, 1],[3, 1, 0, 1, 0, 0]],
[[0, 1, 0, 1, 1, 1],[0, 1, 1, 0, 0, 0],[0, 1, 1, 0, 0, 1],[1, 0, 0, 0, 1, 1],[3, 1, 0, 1, 0, 1]],
[[0, 1, 1, 1, 0, 0],[0, 1, 1, 1, 0, 1],[0, 1, 1, 1, 1, 0],[1, 0, 0, 0, 1, 0],[3, 1, 1, 1, 0, 1]],
[[0, 1, 0, 1, 0, 1],[0, 1, 1, 1, 1, 1],[1, 0, 0, 1, 0, 0],[1, 0, 0, 1, 0, 1],[2, 0, 1, 0, 1, 1]],
[[1, 0, 0, 1, 1, 0],[1, 0, 0, 1, 1, 1],[1, 0, 1, 0, 0, 0],[2, 0, 0, 0, 0, 0],[3, 0, 1, 0, 0, 1]],
[[1, 0, 1, 0, 1, 1],[1, 0, 1, 1, 0, 0],[1, 0, 1, 1, 0, 1],[2, 0, 0, 0, 0, 1],[3, 0, 1, 0, 1, 1]],
[[1, 1, 0, 0, 0, 0],[1, 1, 0, 0, 0, 1],[1, 1, 0, 0, 1, 0],[2, 0, 0, 0, 1, 0],[3, 1, 0, 0, 0, 1]],
[[1, 1, 0, 1, 0, 1],[1, 1, 0, 1, 1, 0],[1, 1, 0, 1, 1, 1],[2, 0, 0, 0, 1, 1],[3, 1, 0, 1, 1, 1]],
[[1, 1, 1, 0, 1, 0],[1, 1, 1, 0, 1, 1],[1, 1, 1, 1, 0, 0],[2, 0, 0, 1, 0, 0],[3, 1, 1, 0, 0, 1]],
[[1, 1, 1, 1, 1, 1],[2, 0, 0, 1, 0, 1],[3, 0, 0, 0, 0, 0],[3, 0, 0, 0, 0, 1],[3, 1, 1, 0, 1, 1]],
[[1, 1, 1, 1, 0, 1],[2, 0, 0, 1, 1, 0],[3, 0, 0, 0, 1, 0],[3, 0, 0, 0, 1, 1],[3, 1, 1, 0, 1, 0]],
[[2, 0, 1, 0, 0, 1],[2, 0, 1, 0, 1, 0],[2, 0, 1, 1, 0, 0],[3, 0, 0, 1, 0, 1],[3, 0, 1, 0, 1, 0]],
[[2, 0, 1, 1, 1, 0],[2, 0, 1, 1, 1, 1],[2, 1, 0, 0, 0, 0],[3, 0, 0, 1, 1, 1],[3, 1, 0, 1, 1, 0]],
[[2, 1, 0, 0, 1, 1],[2, 1, 0, 1, 0, 0],[2, 1, 0, 1, 0, 1],[3, 0, 1, 1, 0, 0],[3, 1, 1, 1, 1, 0]],
[[2, 1, 1, 0, 0, 0],[2, 1, 1, 0, 0, 1],[2, 1, 1, 0, 1, 0],[3, 0, 0, 1, 0, 0],[3, 1, 1, 1, 1, 1]],
[[2, 1, 1, 1, 0, 1],[2, 1, 1, 1, 1, 0],[2, 1, 1, 1, 1, 1],[3, 0, 1, 1, 1, 0],[3, 1, 0, 0, 1, 0]],
[[1, 0, 1, 0, 0, 1],[2, 0, 0, 1, 1, 1],[3, 0, 0, 1, 1, 0],[3, 1, 0, 0, 0, 0],[3, 1, 1, 0, 0, 0]],
[[2, 1, 0, 0, 0, 1],[2, 1, 0, 0, 1, 0],[2, 1, 0, 1, 1, 1],[3, 0, 1, 0, 0, 0],[3, 1, 1, 1, 0, 0]],
[[1, 0, 1, 1, 1, 0],[1, 1, 0, 0, 1, 1],[1, 1, 1, 0, 0, 0],[2, 0, 1, 0, 0, 0],[3, 0, 1, 1, 0, 1]],
[[1, 1, 0, 1, 0, 0],[1, 1, 1, 0, 0, 1],[1, 1, 1, 1, 1, 0],[2, 1, 1, 1, 0, 0],[3, 0, 1, 1, 1, 1]],
[[0, 0, 0, 1, 1, 0],[0, 0, 1, 0, 1, 1],[0, 1, 1, 0, 1, 1],[2, 0, 1, 1, 0, 1],[2, 1, 1, 0, 1, 1]],
[[0, 0, 0, 1, 0, 1],[0, 1, 0, 1, 1, 0],[1, 0, 1, 0, 1, 0],[1, 0, 1, 1, 1, 1],[2, 1, 0, 1, 1, 0]]
]
A partition for sets of sizes:  4*3  0*4 23*5
\end{verbatim}
\newpage

\begin{verbatim} 
[4, 2, 2, 2, 2, 2]
[0, 3, 23]
[
[[0, 0, 0, 0, 0, 1],[0, 0, 0, 0, 1, 0],[0, 0, 0, 1, 0, 0],[0, 0, 0, 1, 1, 1]],
[[0, 0, 0, 1, 0, 1],[0, 0, 0, 1, 1, 0],[0, 0, 1, 0, 0, 0],[0, 0, 1, 0, 1, 1]],
[[0, 0, 1, 0, 0, 1],[0, 0, 1, 0, 1, 0],[0, 0, 1, 1, 0, 0],[0, 0, 1, 1, 1, 1]],
[[0, 0, 1, 1, 0, 1],[0, 0, 1, 1, 1, 0],[0, 1, 0, 0, 0, 0],[1, 0, 0, 0, 0, 0],[3, 1, 0, 0, 1, 1]],
[[0, 1, 0, 0, 1, 0],[0, 1, 0, 0, 1, 1],[0, 1, 0, 1, 0, 0],[1, 0, 0, 0, 0, 1],[3, 1, 0, 1, 0, 0]],
[[0, 1, 0, 1, 1, 1],[0, 1, 1, 0, 0, 0],[0, 1, 1, 0, 0, 1],[1, 0, 0, 0, 1, 1],[3, 1, 0, 1, 0, 1]],
[[0, 1, 1, 1, 0, 0],[0, 1, 1, 1, 0, 1],[0, 1, 1, 1, 1, 0],[1, 0, 0, 0, 1, 0],[3, 1, 1, 1, 0, 1]],
[[0, 1, 0, 1, 0, 1],[0, 1, 1, 1, 1, 1],[1, 0, 0, 1, 0, 0],[1, 0, 0, 1, 0, 1],[2, 0, 1, 0, 1, 1]],
[[1, 0, 0, 1, 1, 0],[1, 0, 0, 1, 1, 1],[1, 0, 1, 0, 0, 0],[2, 0, 0, 0, 0, 0],[3, 0, 1, 0, 0, 1]],
[[1, 0, 1, 0, 1, 1],[1, 0, 1, 1, 0, 0],[1, 0, 1, 1, 0, 1],[2, 0, 0, 0, 0, 1],[3, 0, 1, 0, 1, 1]],
[[1, 1, 0, 0, 0, 0],[1, 1, 0, 0, 0, 1],[1, 1, 0, 0, 1, 0],[2, 0, 0, 0, 1, 0],[3, 1, 0, 0, 0, 1]],
[[1, 1, 0, 1, 0, 1],[1, 1, 0, 1, 1, 0],[1, 1, 0, 1, 1, 1],[2, 0, 0, 0, 1, 1],[3, 1, 0, 1, 1, 1]],
[[1, 1, 1, 0, 1, 0],[1, 1, 1, 0, 1, 1],[1, 1, 1, 1, 0, 0],[2, 0, 0, 1, 0, 0],[3, 1, 1, 0, 0, 1]],
[[1, 1, 1, 1, 1, 1],[2, 0, 0, 1, 0, 1],[3, 0, 0, 0, 0, 0],[3, 0, 0, 0, 0, 1],[3, 1, 1, 0, 1, 1]],
[[1, 1, 1, 1, 0, 1],[2, 0, 0, 1, 1, 0],[3, 0, 0, 0, 1, 0],[3, 0, 0, 0, 1, 1],[3, 1, 1, 0, 1, 0]],
[[2, 0, 1, 0, 0, 1],[2, 0, 1, 0, 1, 0],[2, 0, 1, 1, 0, 0],[3, 0, 0, 1, 0, 1],[3, 0, 1, 0, 1, 0]],
[[2, 0, 1, 1, 1, 0],[2, 0, 1, 1, 1, 1],[2, 1, 0, 0, 0, 0],[3, 0, 0, 1, 1, 1],[3, 1, 0, 1, 1, 0]],
[[2, 1, 0, 0, 1, 1],[2, 1, 0, 1, 0, 0],[2, 1, 0, 1, 0, 1],[3, 0, 1, 1, 0, 0],[3, 1, 1, 1, 1, 0]],
[[2, 1, 1, 0, 0, 0],[2, 1, 1, 0, 0, 1],[2, 1, 1, 0, 1, 0],[3, 0, 0, 1, 0, 0],[3, 1, 1, 1, 1, 1]],
[[2, 1, 1, 1, 0, 1],[2, 1, 1, 1, 1, 0],[2, 1, 1, 1, 1, 1],[3, 0, 1, 1, 1, 0],[3, 1, 0, 0, 1, 0]],
[[1, 0, 1, 0, 0, 1],[2, 0, 0, 1, 1, 1],[3, 0, 0, 1, 1, 0],[3, 1, 0, 0, 0, 0],[3, 1, 1, 0, 0, 0]],
[[2, 1, 0, 0, 0, 1],[2, 1, 0, 0, 1, 0],[2, 1, 0, 1, 1, 1],[3, 0, 1, 0, 0, 0],[3, 1, 1, 1, 0, 0]],
[[1, 0, 1, 1, 1, 0],[1, 1, 0, 0, 1, 1],[1, 1, 1, 0, 0, 0],[2, 0, 1, 0, 0, 0],[3, 0, 1, 1, 0, 1]],
[[0, 1, 0, 0, 0, 1],[1, 1, 0, 1, 0, 0],[2, 1, 0, 1, 1, 0],[2, 1, 1, 1, 0, 0],[3, 0, 1, 1, 1, 1]],
[[0, 0, 0, 0, 1, 1],[0, 1, 1, 0, 1, 0],[1, 0, 1, 0, 1, 0],[1, 1, 1, 1, 1, 0],[2, 0, 1, 1, 0, 1]],
[[0, 1, 0, 1, 1, 0],[0, 1, 1, 0, 1, 1],[1, 0, 1, 1, 1, 1],[1, 1, 1, 0, 0, 1],[2, 1, 1, 0, 1, 1]]
]
A partition for sets of sizes:  0*3  3*4 23*5

[4, 2, 2, 2, 2, 2]
[1, 1, 24]
[
[[0, 0, 0, 0, 0, 1],[0, 0, 0, 0, 1, 0],[0, 0, 0, 0, 1, 1]],
[[0, 0, 0, 1, 0, 0],[0, 0, 0, 1, 0, 1],[0, 0, 0, 1, 1, 0],[0, 0, 0, 1, 1, 1]],
[[0, 0, 1, 0, 0, 0],[0, 0, 1, 0, 0, 1],[0, 0, 1, 0, 1, 0],[0, 1, 0, 0, 0, 0],[0, 1, 1, 0, 1, 1]],
[[0, 0, 1, 1, 0, 1],[0, 0, 1, 1, 1, 0],[0, 0, 1, 1, 1, 1],[0, 1, 0, 0, 0, 1],[0, 1, 1, 1, 0, 1]],
[[0, 1, 0, 0, 1, 0],[0, 1, 0, 0, 1, 1],[0, 1, 0, 1, 0, 0],[1, 0, 0, 0, 0, 0],[3, 1, 0, 1, 0, 1]],
[[0, 1, 0, 1, 1, 1],[0, 1, 1, 0, 0, 0],[0, 1, 1, 0, 0, 1],[1, 0, 0, 0, 0, 1],[3, 1, 0, 1, 1, 1]],
[[0, 1, 1, 1, 0, 0],[0, 1, 1, 1, 1, 0],[0, 1, 1, 1, 1, 1],[1, 0, 0, 0, 1, 0],[3, 1, 1, 1, 1, 1]],
[[0, 1, 1, 0, 1, 0],[1, 0, 0, 0, 1, 1],[1, 0, 0, 1, 0, 0],[1, 0, 0, 1, 0, 1],[1, 1, 1, 0, 0, 0]],
[[1, 0, 0, 1, 1, 0],[1, 0, 0, 1, 1, 1],[1, 0, 1, 0, 0, 0],[2, 0, 0, 0, 0, 0],[3, 0, 1, 0, 0, 1]],
[[1, 0, 1, 0, 1, 1],[1, 0, 1, 1, 0, 0],[1, 0, 1, 1, 0, 1],[2, 0, 0, 0, 0, 1],[3, 0, 1, 0, 1, 1]],
[[1, 1, 0, 0, 0, 0],[1, 1, 0, 0, 0, 1],[1, 1, 0, 0, 1, 0],[2, 0, 0, 0, 1, 0],[3, 1, 0, 0, 0, 1]],
[[1, 1, 0, 1, 0, 1],[1, 1, 0, 1, 1, 0],[1, 1, 0, 1, 1, 1],[2, 0, 0, 1, 0, 0],[3, 1, 0, 0, 0, 0]],
[[1, 1, 1, 0, 1, 0],[1, 1, 1, 0, 1, 1],[1, 1, 1, 1, 0, 0],[2, 0, 0, 0, 1, 1],[3, 1, 1, 1, 1, 0]],
[[1, 1, 1, 1, 1, 1],[2, 0, 0, 1, 0, 1],[3, 0, 0, 0, 0, 0],[3, 0, 0, 0, 0, 1],[3, 1, 1, 0, 1, 1]],
[[0, 0, 1, 0, 1, 1],[1, 0, 1, 0, 0, 1],[2, 0, 0, 1, 1, 0],[2, 0, 0, 1, 1, 1],[3, 0, 0, 0, 1, 1]],
[[2, 0, 1, 0, 0, 1],[2, 0, 1, 0, 1, 0],[2, 0, 1, 0, 1, 1],[3, 0, 0, 0, 1, 0],[3, 0, 1, 0, 1, 0]],
[[2, 0, 1, 1, 1, 0],[2, 0, 1, 1, 1, 1],[2, 1, 0, 0, 0, 0],[3, 0, 0, 1, 0, 1],[3, 1, 0, 1, 0, 0]],
[[2, 1, 0, 0, 1, 1],[2, 1, 0, 1, 0, 0],[2, 1, 0, 1, 0, 1],[3, 0, 0, 1, 0, 0],[3, 1, 0, 1, 1, 0]],
[[2, 1, 1, 0, 0, 0],[2, 1, 1, 0, 0, 1],[2, 1, 1, 0, 1, 0],[3, 0, 0, 1, 1, 0],[3, 1, 1, 1, 0, 1]],
[[2, 1, 1, 1, 0, 1],[2, 1, 1, 1, 1, 0],[2, 1, 1, 1, 1, 1],[3, 0, 1, 1, 1, 0],[3, 1, 0, 0, 1, 0]],
[[2, 0, 1, 1, 0, 0],[2, 1, 0, 1, 1, 0],[2, 1, 1, 0, 1, 1],[3, 0, 1, 1, 0, 0],[3, 0, 1, 1, 0, 1]],
[[1, 0, 1, 0, 1, 0],[2, 0, 1, 1, 0, 1],[3, 0, 1, 0, 0, 0],[3, 1, 0, 0, 1, 1],[3, 1, 1, 1, 0, 0]],
[[2, 1, 0, 0, 0, 1],[2, 1, 0, 0, 1, 0],[2, 1, 1, 1, 0, 0],[3, 0, 0, 1, 1, 1],[3, 1, 1, 0, 0, 0]],
[[0, 1, 0, 1, 1, 0],[1, 1, 0, 0, 1, 1],[2, 0, 1, 0, 0, 0],[2, 1, 0, 1, 1, 1],[3, 1, 1, 0, 1, 0]],
[[0, 0, 1, 1, 0, 0],[1, 0, 1, 1, 1, 0],[1, 1, 0, 1, 0, 0],[3, 0, 1, 1, 1, 1],[3, 1, 1, 0, 0, 1]],
[[0, 1, 0, 1, 0, 1],[1, 0, 1, 1, 1, 1],[1, 1, 1, 0, 0, 1],[1, 1, 1, 1, 0, 1],[1, 1, 1, 1, 1, 0]]
]
A partition for sets of sizes:  1*3  1*4 24*5
\end{verbatim}
\newpage

\section{\texorpdfstring{Zero sum partitions of $\left((\zet_4)^{2}\times\zet_2\right)^*$ (with $b\leq 3$)}{Zero sum partitions of ((Z4)**2 x Z2)* (with b <= 3)}}\label{a:442}
\begin{verbatim}
[4, 4, 2]
[7, 0, 2]
[
[[0, 0, 1],[0, 1, 0],[0, 3, 1]],
[[0, 2, 0],[1, 0, 0],[3, 2, 0]],
[[0, 1, 1],[1, 0, 1],[3, 3, 0]],
[[1, 1, 0],[1, 1, 1],[2, 2, 1]],
[[1, 2, 1],[1, 3, 0],[2, 3, 1]],
[[2, 0, 0],[3, 1, 1],[3, 3, 1]],
[[2, 1, 1],[3, 1, 0],[3, 2, 1]],
[[2, 0, 1],[2, 1, 0],[2, 3, 0],[3, 0, 0],[3, 0, 1]],
[[0, 2, 1],[0, 3, 0],[1, 2, 0],[1, 3, 1],[2, 2, 0]]
]
A partition for sets of sizes:  7*3  0*4  2*5

[4, 4, 2]
[2, 0, 5]
[
[[0, 0, 1],[0, 1, 0],[0, 3, 1]],
[[0, 2, 0],[1, 0, 0],[3, 2, 0]],
[[0, 1, 1],[0, 2, 1],[1, 0, 1],[1, 1, 0],[2, 0, 1]],
[[1, 2, 0],[1, 2, 1],[1, 3, 0],[2, 0, 0],[3, 1, 1]],
[[1, 1, 1],[2, 1, 0],[3, 0, 0],[3, 3, 0],[3, 3, 1]],
[[2, 1, 1],[2, 2, 0],[2, 3, 1],[3, 0, 1],[3, 2, 1]],
[[0, 3, 0],[1, 3, 1],[2, 2, 1],[2, 3, 0],[3, 1, 0]]
]
A partition for sets of sizes:  2*3  0*4  5*5
\end{verbatim}

\newpage
\section{\texorpdfstring{Zero sum partitions of $\left((\zet_4)^{3}\right)^*$ (with $b\leq 3$)}{Zero sum partitions of ((Z4)**3)* (with b <= 3)}}\label{a:444}
\begin{verbatim}
[4, 4, 4]
[21, 0, 0]
[
[[0, 0, 1],[0, 1, 0],[0, 3, 3]],
[[0, 0, 2],[0, 1, 1],[0, 3, 1]],
[[0, 1, 3],[1, 0, 0],[3, 3, 1]],
[[0, 2, 2],[0, 3, 0],[0, 3, 2]],
[[0, 1, 2],[1, 0, 2],[3, 3, 0]],
[[0, 2, 0],[1, 0, 1],[3, 2, 3]],
[[1, 0, 3],[1, 1, 0],[2, 3, 1]],
[[1, 1, 2],[1, 1, 3],[2, 2, 3]],
[[1, 2, 1],[1, 2, 2],[2, 0, 1]],
[[1, 3, 0],[1, 3, 2],[2, 2, 2]],
[[2, 0, 0],[3, 0, 1],[3, 0, 3]],
[[2, 0, 2],[3, 0, 0],[3, 0, 2]],
[[2, 1, 1],[3, 1, 1],[3, 2, 2]],
[[2, 2, 1],[3, 1, 0],[3, 1, 3]],
[[1, 2, 0],[1, 3, 1],[2, 3, 3]],
[[2, 3, 2],[3, 2, 0],[3, 3, 2]],
[[2, 3, 0],[3, 2, 1],[3, 3, 3]],
[[1, 1, 1],[1, 2, 3],[2, 1, 0]],
[[0, 2, 3],[2, 1, 2],[2, 1, 3]],
[[0, 0, 3],[1, 3, 3],[3, 1, 2]],
[[0, 2, 1],[2, 0, 3],[2, 2, 0]]
]
A partition for sets of sizes: 21*3  0*4  0*5

[4, 4, 4]
[16, 0, 3]
[
[[0, 0, 1],[0, 1, 0],[0, 3, 3]],
[[0, 0, 2],[0, 1, 1],[0, 3, 1]],
[[0, 1, 3],[1, 0, 0],[3, 3, 1]],
[[0, 2, 2],[0, 3, 0],[0, 3, 2]],
[[0, 1, 2],[1, 0, 2],[3, 3, 0]],
[[0, 2, 0],[1, 0, 1],[3, 2, 3]],
[[1, 0, 3],[1, 1, 0],[2, 3, 1]],
[[1, 1, 2],[1, 1, 3],[2, 2, 3]],
[[1, 2, 1],[1, 2, 2],[2, 0, 1]],
[[1, 3, 0],[1, 3, 2],[2, 2, 2]],
[[2, 0, 0],[3, 0, 1],[3, 0, 3]],
[[2, 0, 2],[3, 0, 0],[3, 0, 2]],
[[2, 1, 1],[3, 1, 1],[3, 2, 2]],
[[2, 2, 1],[3, 1, 0],[3, 1, 3]],
[[1, 2, 0],[1, 3, 1],[2, 3, 3]],
[[2, 3, 2],[3, 2, 0],[3, 3, 2]],
[[2, 1, 0],[2, 1, 3],[2, 2, 0],[3, 1, 2],[3, 3, 3]],
[[0, 2, 3],[1, 1, 1],[2, 0, 3],[2, 3, 0],[3, 2, 1]],
[[0, 0, 3],[0, 2, 1],[1, 2, 3],[1, 3, 3],[2, 1, 2]]
]
A partition for sets of sizes: 16*3  0*4  3*5
\end{verbatim}

\newpage

\begin{verbatim} 
[4, 4, 4]
[11, 0, 6]
[
[[0, 0, 1],[0, 1, 0],[0, 3, 3]],
[[0, 0, 2],[0, 1, 1],[0, 3, 1]],
[[0, 1, 3],[1, 0, 0],[3, 3, 1]],
[[0, 2, 2],[0, 3, 0],[0, 3, 2]],
[[0, 1, 2],[1, 0, 2],[3, 3, 0]],
[[0, 2, 0],[1, 0, 1],[3, 2, 3]],
[[1, 0, 3],[1, 1, 0],[2, 3, 1]],
[[1, 1, 2],[1, 1, 3],[2, 2, 3]],
[[1, 2, 1],[1, 2, 2],[2, 0, 1]],
[[1, 3, 0],[1, 3, 2],[2, 2, 2]],
[[2, 0, 0],[3, 0, 1],[3, 0, 3]],
[[2, 0, 2],[2, 0, 3],[2, 1, 0],[3, 0, 0],[3, 3, 3]],
[[2, 1, 3],[2, 2, 0],[2, 2, 1],[3, 0, 2],[3, 3, 2]],
[[2, 3, 0],[2, 3, 2],[2, 3, 3],[3, 1, 1],[3, 2, 2]],
[[1, 3, 3],[2, 1, 1],[3, 1, 0],[3, 1, 3],[3, 2, 1]],
[[0, 0, 3],[1, 1, 1],[1, 2, 0],[1, 2, 3],[1, 3, 1]],
[[0, 2, 1],[0, 2, 3],[2, 1, 2],[3, 1, 2],[3, 2, 0]]
]
A partition for sets of sizes: 11*3  0*4  6*5

[4, 4, 4]
[6, 0, 9]
[
[[0, 0, 1],[0, 1, 0],[0, 3, 3]],
[[0, 0, 2],[0, 1, 1],[0, 3, 1]],
[[0, 1, 3],[1, 0, 0],[3, 3, 1]],
[[0, 2, 2],[0, 3, 0],[0, 3, 2]],
[[0, 1, 2],[1, 0, 2],[3, 3, 0]],
[[0, 2, 0],[1, 0, 1],[3, 2, 3]],
[[1, 0, 3],[1, 1, 0],[1, 1, 1],[2, 0, 0],[3, 2, 0]],
[[1, 2, 0],[1, 2, 1],[1, 2, 2],[2, 0, 3],[3, 2, 2]],
[[1, 3, 1],[1, 3, 2],[1, 3, 3],[2, 1, 1],[3, 2, 1]],
[[2, 0, 2],[2, 1, 0],[2, 1, 2],[3, 1, 1],[3, 1, 3]],
[[2, 1, 3],[2, 2, 0],[2, 2, 1],[3, 0, 1],[3, 3, 3]],
[[2, 3, 0],[2, 3, 1],[2, 3, 2],[3, 0, 3],[3, 3, 2]],
[[2, 2, 2],[2, 2, 3],[2, 3, 3],[3, 0, 0],[3, 1, 0]],
[[0, 0, 3],[0, 2, 3],[1, 1, 2],[1, 1, 3],[2, 0, 1]],
[[0, 2, 1],[1, 2, 3],[1, 3, 0],[3, 0, 2],[3, 1, 2]]
]
A partition for sets of sizes:  6*3  0*4  9*5

[4, 4, 4]
[1, 0, 12]
[
[[0, 0, 1],[0, 1, 0],[0, 3, 3]],
[[0, 0, 2],[0, 1, 1],[0, 1, 2],[0, 3, 1],[0, 3, 2]],
[[0, 2, 1],[0, 2, 2],[0, 2, 3],[1, 0, 0],[3, 2, 2]],
[[0, 2, 0],[0, 3, 0],[1, 0, 1],[1, 0, 2],[2, 3, 1]],
[[1, 0, 3],[1, 1, 0],[1, 1, 1],[2, 0, 0],[3, 2, 0]],
[[1, 2, 0],[1, 2, 1],[1, 2, 2],[2, 0, 2],[3, 2, 3]],
[[1, 3, 1],[1, 3, 2],[1, 3, 3],[2, 0, 1],[3, 3, 1]],
[[1, 2, 3],[2, 0, 3],[3, 0, 0],[3, 0, 1],[3, 2, 1]],
[[2, 1, 3],[2, 2, 0],[2, 2, 1],[3, 0, 2],[3, 3, 2]],
[[2, 3, 0],[2, 3, 2],[2, 3, 3],[3, 0, 3],[3, 3, 0]],
[[2, 1, 1],[2, 1, 2],[2, 2, 2],[3, 1, 0],[3, 3, 3]],
[[0, 0, 3],[1, 1, 2],[1, 1, 3],[3, 1, 1],[3, 1, 3]],
[[0, 1, 3],[1, 3, 0],[2, 1, 0],[2, 2, 3],[3, 1, 2]]
]
A partition for sets of sizes:  1*3  0*4 12*5
\end{verbatim}

\newpage
\section{\texorpdfstring{Zero sum partitions of $\left(\zet_8\times(\zet_2)^{2}\right)^*$ (with $b \leq 11$)}{Zero sum partitions of (Z8 x (Z2)**2)* (with b <= 11)}}\label{a:822}
\begin{verbatim}
[8, 2, 2]
[9, 1, 0]
[
[[0, 0, 1],[0, 1, 0],[0, 1, 1]],
[[1, 0, 0],[1, 0, 1],[6, 0, 1]],
[[1, 1, 1],[2, 0, 0],[5, 1, 1]],
[[2, 1, 0],[2, 1, 1],[4, 0, 1]],
[[3, 0, 1],[6, 0, 0],[7, 0, 1]],
[[4, 0, 0],[5, 0, 0],[7, 0, 0]],
[[4, 1, 1],[5, 0, 1],[7, 1, 0]],
[[2, 0, 1],[3, 1, 0],[3, 1, 1]],
[[1, 1, 0],[3, 0, 0],[4, 1, 0]],
[[5, 1, 0],[6, 1, 0],[6, 1, 1],[7, 1, 1]]
]
A partition for sets of sizes:  9*3  1*4  0*5

[8, 2, 2]
[6, 2, 1]
[
[[0, 0, 1],[0, 1, 0],[0, 1, 1]],
[[1, 0, 0],[1, 0, 1],[6, 0, 1]],
[[1, 1, 1],[2, 0, 0],[5, 1, 1]],
[[2, 1, 0],[2, 1, 1],[4, 0, 1]],
[[3, 0, 1],[6, 0, 0],[7, 0, 1]],
[[4, 0, 0],[5, 0, 0],[7, 0, 0]],
[[4, 1, 1],[6, 1, 0],[7, 1, 0],[7, 1, 1]],
[[2, 0, 1],[3, 0, 0],[5, 1, 0],[6, 1, 1]],
[[1, 1, 0],[3, 1, 0],[3, 1, 1],[4, 1, 0],[5, 0, 1]]
]
A partition for sets of sizes:  6*3  2*4  1*5

[8, 2, 2]
[7, 0, 2]
[
[[0, 0, 1],[0, 1, 0],[0, 1, 1]],
[[1, 0, 0],[1, 0, 1],[6, 0, 1]],
[[1, 1, 1],[2, 0, 0],[5, 1, 1]],
[[2, 1, 0],[2, 1, 1],[4, 0, 1]],
[[3, 0, 1],[6, 0, 0],[7, 0, 1]],
[[4, 0, 0],[5, 0, 0],[7, 0, 0]],
[[4, 1, 1],[5, 0, 1],[7, 1, 0]],
[[2, 0, 1],[4, 1, 0],[5, 1, 0],[6, 1, 0],[7, 1, 1]],
[[1, 1, 0],[3, 0, 0],[3, 1, 0],[3, 1, 1],[6, 1, 1]]
]
A partition for sets of sizes:  7*3  0*4  2*5

[8, 2, 2]
[4, 1, 3]
[
[[0, 0, 1],[0, 1, 0],[0, 1, 1]],
[[1, 0, 0],[1, 0, 1],[6, 0, 1]],
[[1, 1, 1],[2, 0, 0],[5, 1, 1]],
[[2, 1, 0],[2, 1, 1],[4, 0, 1]],
[[3, 0, 1],[3, 1, 0],[3, 1, 1],[7, 0, 0]],
[[3, 0, 0],[4, 1, 0],[4, 1, 1],[6, 0, 0],[7, 0, 1]],
[[2, 0, 1],[5, 0, 0],[5, 0, 1],[5, 1, 0],[7, 1, 0]],
[[1, 1, 0],[4, 0, 0],[6, 1, 0],[6, 1, 1],[7, 1, 1]]
]
A partition for sets of sizes:  4*3  1*4  3*5
\end{verbatim}

\newpage

\begin{verbatim} 
[8, 2, 2]
[1, 2, 4]
[
[[0, 0, 1],[0, 1, 0],[0, 1, 1]],
[[1, 0, 0],[1, 0, 1],[1, 1, 0],[5, 1, 1]],
[[2, 0, 0],[2, 0, 1],[2, 1, 0],[2, 1, 1]],
[[3, 0, 0],[3, 0, 1],[3, 1, 0],[3, 1, 1],[4, 0, 0]],
[[4, 0, 1],[4, 1, 0],[4, 1, 1],[5, 0, 0],[7, 0, 0]],
[[5, 1, 0],[6, 1, 0],[7, 0, 1],[7, 1, 0],[7, 1, 1]],
[[1, 1, 1],[5, 0, 1],[6, 0, 0],[6, 0, 1],[6, 1, 1]]
]
A partition for sets of sizes:  1*3  2*4  4*5

[8, 2, 2]
[2, 0, 5]
[
[[0, 0, 1],[0, 1, 0],[0, 1, 1]],
[[1, 0, 0],[1, 0, 1],[6, 0, 1]],
[[1, 1, 1],[2, 0, 0],[2, 0, 1],[4, 0, 0],[7, 1, 0]],
[[3, 0, 0],[3, 0, 1],[4, 1, 0],[7, 0, 0],[7, 1, 1]],
[[4, 0, 1],[4, 1, 1],[5, 0, 0],[5, 0, 1],[6, 1, 1]],
[[1, 1, 0],[5, 1, 0],[5, 1, 1],[6, 1, 0],[7, 0, 1]],
[[2, 1, 0],[2, 1, 1],[3, 1, 0],[3, 1, 1],[6, 0, 0]]
]
A partition for sets of sizes:  2*3  0*4  5*5
\end{verbatim}

\newpage
\section{\texorpdfstring{Zero sum partitions of $\left(\zet_4\times(\zet_2)^{2}\right)^*+\left((\zet_2)^{2}\right)^*$}{Zero sum partitions of (Z4 x (Z2)**2)* + ((Z2)**2)*}}\label{a:422+22}
\begin{verbatim}
[4, 2, 2, 2, 2]
[15, 0, 0]
[
[[0, 0, 1, 0, 1],[0, 0, 0, 1, 0],[0, 0, 1, 1, 1]],
[[0, 0, 0, 0, 1],[1, 0, 0, 1, 0],[3, 0, 0, 1, 1]],
[[0, 0, 1, 0, 1],[1, 0, 0, 1, 1],[3, 0, 1, 1, 0]],
[[1, 0, 0, 0, 1],[1, 0, 1, 1, 0],[2, 0, 1, 1, 1]],
[[1, 0, 1, 0, 1],[1, 0, 1, 1, 1],[2, 0, 0, 1, 0]],
[[1, 0, 0, 0, 1],[1, 0, 0, 1, 0],[2, 0, 0, 1, 1]],
[[1, 0, 1, 1, 0],[1, 0, 1, 1, 1],[2, 0, 0, 0, 1]],
[[0, 0, 1, 1, 0],[1, 0, 0, 1, 1],[3, 0, 1, 0, 1]],
[[2, 0, 1, 0, 1],[3, 0, 0, 1, 0],[3, 0, 1, 1, 1]],
[[2, 0, 0, 1, 0],[3, 0, 1, 0, 1],[3, 0, 1, 1, 1]],
[[2, 0, 1, 1, 1],[3, 0, 1, 1, 0],[3, 0, 0, 0, 1]],
[[0, 0, 1, 1, 0],[1, 0, 1, 0, 1],[3, 0, 0, 1, 1]],
[[0, 0, 0, 1, 1],[2, 0, 1, 1, 0],[2, 0, 1, 0, 1]],
[[0, 0, 1, 1, 1],[2, 0, 0, 0, 1],[2, 0, 1, 1, 0]],
[[2, 0, 0, 1, 1],[3, 0, 0, 0, 1],[3, 0, 0, 1, 0]]
]
A partition for sets of sizes: 15*3  0*4  0*5

[4, 2, 2, 2, 2]
[11, 3, 0]
[
[[0, 0, 1, 0, 1],[0, 0, 0, 1, 0],[0, 0, 1, 1, 1]],
[[0, 0, 0, 0, 1],[1, 0, 0, 1, 0],[3, 0, 0, 1, 1]],
[[0, 0, 1, 0, 1],[1, 0, 0, 1, 1],[3, 0, 1, 1, 0]],
[[1, 0, 0, 0, 1],[1, 0, 1, 1, 0],[2, 0, 1, 1, 1]],
[[1, 0, 1, 0, 1],[1, 0, 1, 1, 1],[2, 0, 0, 1, 0]],
[[1, 0, 0, 0, 1],[1, 0, 0, 1, 0],[2, 0, 0, 1, 1]],
[[1, 0, 1, 1, 0],[1, 0, 1, 1, 1],[2, 0, 0, 0, 1]],
[[0, 0, 1, 1, 0],[1, 0, 0, 1, 1],[3, 0, 1, 0, 1]],
[[2, 0, 1, 0, 1],[3, 0, 0, 1, 0],[3, 0, 1, 1, 1]],
[[2, 0, 0, 0, 1],[3, 0, 0, 1, 1],[3, 0, 0, 1, 0]],
[[2, 0, 1, 1, 0],[3, 0, 0, 0, 1],[3, 0, 1, 1, 1]],
[[0, 0, 1, 1, 0],[2, 0, 1, 0, 1],[3, 0, 1, 0, 1],[3, 0, 1, 1, 0]],
[[1, 0, 1, 0, 1],[2, 0, 0, 1, 1],[2, 0, 1, 1, 1],[3, 0, 0, 0, 1]],
[[0, 0, 1, 1, 1],[0, 0, 0, 1, 1],[2, 0, 1, 1, 0],[2, 0, 0, 1, 0]]
]
A partition for sets of sizes: 11*3  3*4  0*5

[4, 2, 2, 2, 2]
[7, 6, 0]
[
[[0, 0, 1, 0, 1],[0, 0, 0, 1, 0],[0, 0, 1, 1, 1]],
[[0, 0, 0, 0, 1],[1, 0, 0, 1, 0],[3, 0, 0, 1, 1]],
[[0, 0, 1, 0, 1],[1, 0, 0, 1, 1],[3, 0, 1, 1, 0]],
[[1, 0, 0, 0, 1],[1, 0, 1, 1, 0],[2, 0, 1, 1, 1]],
[[1, 0, 1, 0, 1],[1, 0, 1, 1, 1],[2, 0, 0, 1, 0]],
[[1, 0, 0, 0, 1],[1, 0, 0, 1, 0],[2, 0, 0, 1, 1]],
[[1, 0, 1, 1, 0],[1, 0, 1, 1, 1],[2, 0, 0, 0, 1]],
[[1, 0, 1, 0, 1],[2, 0, 1, 0, 1],[2, 0, 0, 0, 1],[3, 0, 0, 0, 1]],
[[0, 0, 1, 1, 0],[2, 0, 1, 1, 0],[3, 0, 0, 1, 0],[3, 0, 0, 1, 0]],
[[0, 0, 1, 1, 0],[2, 0, 1, 0, 1],[3, 0, 1, 1, 0],[3, 0, 1, 0, 1]],
[[1, 0, 0, 1, 1],[2, 0, 0, 1, 0],[2, 0, 1, 1, 0],[3, 0, 1, 1, 1]],
[[0, 0, 0, 1, 1],[2, 0, 1, 1, 1],[3, 0, 0, 1, 1],[3, 0, 1, 1, 1]],
[[0, 0, 1, 1, 1],[2, 0, 0, 1, 1],[3, 0, 1, 0, 1],[3, 0, 0, 0, 1]]
]
A partition for sets of sizes:  7*3  6*4  0*5
\end{verbatim}

\newpage

\begin{verbatim} 
[4, 2, 2, 2, 2]
[3, 9, 0]
[
[[0, 0, 1, 0, 1],[0, 0, 0, 1, 0],[0, 0, 1, 1, 1]],
[[0, 0, 0, 0, 1],[1, 0, 0, 1, 0],[3, 0, 0, 1, 1]],
[[0, 0, 1, 0, 1],[1, 0, 0, 1, 1],[3, 0, 1, 1, 0]],
[[1, 0, 0, 0, 1],[1, 0, 1, 0, 1],[1, 0, 0, 0, 1],[1, 0, 1, 0, 1]],
[[1, 0, 1, 1, 0],[1, 0, 1, 1, 1],[1, 0, 0, 1, 0],[1, 0, 0, 1, 1]],
[[0, 0, 1, 1, 0],[1, 0, 1, 1, 0],[1, 0, 1, 1, 1],[2, 0, 1, 1, 1]],
[[2, 0, 0, 0, 1],[2, 0, 0, 1, 0],[2, 0, 1, 0, 1],[2, 0, 1, 1, 0]],
[[2, 0, 0, 1, 1],[2, 0, 1, 1, 1],[2, 0, 0, 0, 1],[2, 0, 1, 0, 1]],
[[0, 0, 1, 1, 0],[2, 0, 0, 1, 0],[3, 0, 0, 0, 1],[3, 0, 1, 0, 1]],
[[3, 0, 0, 1, 0],[3, 0, 0, 1, 1],[3, 0, 1, 1, 0],[3, 0, 1, 1, 1]],
[[0, 0, 0, 1, 1],[2, 0, 1, 1, 0],[3, 0, 0, 1, 0],[3, 0, 1, 1, 1]],
[[0, 0, 1, 1, 1],[2, 0, 0, 1, 1],[3, 0, 1, 0, 1],[3, 0, 0, 0, 1]]
]
A partition for sets of sizes:  3*3  9*4  0*5

[4, 2, 2, 2, 2]
[12, 1, 1]
[
[[0, 0, 1, 0, 1],[0, 0, 0, 1, 0],[0, 0, 1, 1, 1]],
[[0, 0, 0, 0, 1],[1, 0, 0, 1, 0],[3, 0, 0, 1, 1]],
[[0, 0, 1, 0, 1],[1, 0, 0, 1, 1],[3, 0, 1, 1, 0]],
[[1, 0, 0, 0, 1],[1, 0, 1, 1, 0],[2, 0, 1, 1, 1]],
[[1, 0, 1, 0, 1],[1, 0, 1, 1, 1],[2, 0, 0, 1, 0]],
[[1, 0, 0, 0, 1],[1, 0, 0, 1, 0],[2, 0, 0, 1, 1]],
[[1, 0, 1, 1, 0],[1, 0, 1, 1, 1],[2, 0, 0, 0, 1]],
[[0, 0, 1, 1, 0],[1, 0, 0, 1, 1],[3, 0, 1, 0, 1]],
[[2, 0, 1, 0, 1],[3, 0, 0, 1, 0],[3, 0, 1, 1, 1]],
[[2, 0, 0, 0, 1],[3, 0, 0, 1, 1],[3, 0, 0, 1, 0]],
[[2, 0, 1, 1, 0],[3, 0, 0, 0, 1],[3, 0, 1, 1, 1]],
[[0, 0, 0, 1, 1],[2, 0, 1, 1, 0],[2, 0, 1, 0, 1]],
[[1, 0, 1, 0, 1],[2, 0, 0, 1, 1],[2, 0, 1, 1, 1],[3, 0, 0, 0, 1]],
[[0, 0, 1, 1, 1],[0, 0, 1, 1, 0],[2, 0, 0, 1, 0],[3, 0, 1, 0, 1],[3, 0, 1, 1, 0]]
]
A partition for sets of sizes: 12*3  1*4  1*5

[4, 2, 2, 2, 2]
[8, 4, 1]
[
[[0, 0, 1, 0, 1],[0, 0, 0, 1, 0],[0, 0, 1, 1, 1]],
[[0, 0, 0, 0, 1],[1, 0, 0, 1, 0],[3, 0, 0, 1, 1]],
[[0, 0, 1, 0, 1],[1, 0, 0, 1, 1],[3, 0, 1, 1, 0]],
[[1, 0, 0, 0, 1],[1, 0, 1, 1, 0],[2, 0, 1, 1, 1]],
[[1, 0, 1, 0, 1],[1, 0, 1, 1, 1],[2, 0, 0, 1, 0]],
[[1, 0, 0, 0, 1],[1, 0, 0, 1, 0],[2, 0, 0, 1, 1]],
[[1, 0, 1, 1, 0],[1, 0, 1, 1, 1],[2, 0, 0, 0, 1]],
[[0, 0, 1, 1, 0],[1, 0, 0, 1, 1],[3, 0, 1, 0, 1]],
[[2, 0, 1, 0, 1],[2, 0, 1, 1, 0],[2, 0, 0, 0, 1],[2, 0, 0, 1, 0]],
[[0, 0, 1, 1, 0],[2, 0, 0, 1, 1],[3, 0, 0, 1, 0],[3, 0, 1, 1, 1]],
[[3, 0, 0, 1, 1],[3, 0, 1, 0, 1],[3, 0, 0, 0, 1],[3, 0, 1, 1, 1]],
[[0, 0, 1, 1, 1],[0, 0, 0, 1, 1],[1, 0, 1, 0, 1],[3, 0, 0, 0, 1]],
[[2, 0, 1, 0, 1],[2, 0, 1, 1, 0],[2, 0, 1, 1, 1],[3, 0, 1, 1, 0],[3, 0, 0, 1, 0]]
]
A partition for sets of sizes:  8*3  4*4  1*5
\end{verbatim}

\newpage

\begin{verbatim} 
[4, 2, 2, 2, 2]
[4, 7, 1]
[
[[0, 0, 1, 0, 1],[0, 0, 0, 1, 0],[0, 0, 1, 1, 1]],
[[0, 0, 0, 0, 1],[1, 0, 0, 1, 0],[3, 0, 0, 1, 1]],
[[0, 0, 1, 0, 1],[1, 0, 0, 1, 1],[3, 0, 1, 1, 0]],
[[1, 0, 0, 0, 1],[1, 0, 1, 1, 0],[2, 0, 1, 1, 1]],
[[1, 0, 1, 0, 1],[1, 0, 1, 1, 1],[1, 0, 0, 0, 1],[1, 0, 0, 1, 1]],
[[1, 0, 0, 1, 0],[1, 0, 1, 0, 1],[3, 0, 0, 0, 1],[3, 0, 1, 1, 0]],
[[1, 0, 1, 1, 1],[2, 0, 0, 0, 1],[2, 0, 0, 1, 1],[3, 0, 1, 0, 1]],
[[2, 0, 1, 0, 1],[2, 0, 1, 1, 0],[2, 0, 0, 0, 1],[2, 0, 0, 1, 0]],
[[0, 0, 1, 1, 0],[2, 0, 0, 1, 1],[3, 0, 0, 1, 0],[3, 0, 1, 1, 1]],
[[0, 0, 1, 1, 0],[2, 0, 1, 1, 1],[3, 0, 0, 1, 1],[3, 0, 0, 1, 0]],
[[0, 0, 0, 1, 1],[2, 0, 1, 0, 1],[3, 0, 0, 0, 1],[3, 0, 1, 1, 1]],
[[0, 0, 1, 1, 1],[1, 0, 1, 1, 0],[2, 0, 0, 1, 0],[2, 0, 1, 1, 0],[3, 0, 1, 0, 1]]
]
A partition for sets of sizes:  4*3  7*4  1*5

[4, 2, 2, 2, 2]
[0, 10, 1]
[
[[0, 0, 1, 0, 1],[0, 0, 1, 1, 0],[0, 0, 0, 0, 1],[0, 0, 0, 1, 0]],
[[0, 0, 1, 1, 1],[0, 0, 0, 1, 1],[1, 0, 0, 0, 1],[3, 0, 1, 0, 1]],
[[0, 0, 1, 1, 1],[1, 0, 0, 1, 0],[1, 0, 0, 1, 1],[2, 0, 1, 1, 0]],
[[1, 0, 1, 0, 1],[1, 0, 1, 1, 0],[1, 0, 0, 0, 1],[1, 0, 0, 1, 0]],
[[1, 0, 1, 1, 1],[1, 0, 0, 1, 1],[3, 0, 0, 1, 0],[3, 0, 1, 1, 0]],
[[1, 0, 1, 1, 1],[2, 0, 0, 0, 1],[2, 0, 1, 0, 1],[3, 0, 0, 1, 1]],
[[2, 0, 0, 1, 0],[2, 0, 1, 1, 0],[2, 0, 0, 0, 1],[2, 0, 1, 0, 1]],
[[0, 0, 1, 0, 1],[2, 0, 0, 1, 1],[3, 0, 0, 0, 1],[3, 0, 1, 1, 1]],
[[1, 0, 1, 0, 1],[2, 0, 1, 1, 1],[2, 0, 1, 1, 1],[3, 0, 1, 0, 1]],
[[3, 0, 0, 1, 1],[3, 0, 1, 1, 0],[3, 0, 0, 1, 0],[3, 0, 1, 1, 1]],
[[0, 0, 1, 1, 0],[1, 0, 1, 1, 0],[2, 0, 0, 1, 1],[2, 0, 0, 1, 0],[3, 0, 0, 0, 1]]
]
A partition for sets of sizes:  0*3 10*4  1*5


[4, 2, 2, 2, 2]
[9, 2, 2]
[
[[0, 0, 1, 0, 1],[0, 0, 0, 1, 0],[0, 0, 1, 1, 1]],
[[0, 0, 0, 0, 1],[1, 0, 0, 1, 0],[3, 0, 0, 1, 1]],
[[0, 0, 1, 0, 1],[1, 0, 0, 1, 1],[3, 0, 1, 1, 0]],
[[1, 0, 0, 0, 1],[1, 0, 1, 1, 0],[2, 0, 1, 1, 1]],
[[1, 0, 1, 0, 1],[1, 0, 1, 1, 1],[2, 0, 0, 1, 0]],
[[1, 0, 0, 0, 1],[1, 0, 0, 1, 0],[2, 0, 0, 1, 1]],
[[1, 0, 1, 1, 0],[1, 0, 1, 1, 1],[2, 0, 0, 0, 1]],
[[0, 0, 1, 1, 0],[1, 0, 0, 1, 1],[3, 0, 1, 0, 1]],
[[2, 0, 1, 0, 1],[3, 0, 0, 1, 0],[3, 0, 1, 1, 1]],
[[2, 0, 0, 0, 1],[2, 0, 0, 1, 0],[2, 0, 1, 0, 1],[2, 0, 1, 1, 0]],
[[3, 0, 0, 1, 1],[3, 0, 1, 0, 1],[3, 0, 0, 0, 1],[3, 0, 1, 1, 1]],
[[0, 0, 1, 1, 0],[1, 0, 1, 0, 1],[2, 0, 1, 1, 0],[2, 0, 1, 1, 1],[3, 0, 0, 1, 0]],
[[0, 0, 1, 1, 1],[0, 0, 0, 1, 1],[2, 0, 0, 1, 1],[3, 0, 0, 0, 1],[3, 0, 1, 1, 0]]
]
A partition for sets of sizes:  9*3  2*4  2*5
\end{verbatim}

\newpage

\begin{verbatim} 
[4, 2, 2, 2, 2]
[5, 5, 2]
[
[[0, 0, 1, 0, 1],[0, 0, 0, 1, 0],[0, 0, 1, 1, 1]],
[[0, 0, 0, 0, 1],[1, 0, 0, 1, 0],[3, 0, 0, 1, 1]],
[[0, 0, 1, 0, 1],[1, 0, 0, 1, 1],[3, 0, 1, 1, 0]],
[[1, 0, 0, 0, 1],[1, 0, 1, 1, 0],[2, 0, 1, 1, 1]],
[[1, 0, 1, 0, 1],[1, 0, 1, 1, 1],[2, 0, 0, 1, 0]],
[[1, 0, 0, 0, 1],[1, 0, 0, 1, 0],[1, 0, 1, 0, 1],[1, 0, 1, 1, 0]],
[[1, 0, 0, 1, 1],[1, 0, 1, 1, 1],[3, 0, 0, 0, 1],[3, 0, 1, 0, 1]],
[[2, 0, 0, 1, 1],[2, 0, 1, 0, 1],[2, 0, 0, 0, 1],[2, 0, 1, 1, 1]],
[[0, 0, 1, 1, 0],[2, 0, 0, 1, 0],[3, 0, 0, 1, 0],[3, 0, 1, 1, 0]],
[[3, 0, 0, 1, 1],[3, 0, 1, 1, 1],[3, 0, 0, 0, 1],[3, 0, 1, 0, 1]],
[[0, 0, 1, 1, 0],[0, 0, 1, 1, 1],[0, 0, 0, 1, 1],[2, 0, 0, 0, 1],[2, 0, 0, 1, 1]],
[[2, 0, 1, 1, 0],[2, 0, 1, 0, 1],[2, 0, 1, 1, 0],[3, 0, 0, 1, 0],[3, 0, 1, 1, 1]]
]
A partition for sets of sizes:  5*3  5*4  2*5

[4, 2, 2, 2, 2]
[1, 8, 2]
[
[[0, 0, 1, 0, 1],[0, 0, 0, 1, 0],[0, 0, 1, 1, 1]],
[[0, 0, 0, 0, 1],[0, 0, 0, 1, 1],[1, 0, 0, 0, 1],[3, 0, 0, 1, 1]],
[[0, 0, 1, 1, 0],[1, 0, 0, 1, 0],[1, 0, 0, 1, 1],[2, 0, 1, 1, 1]],
[[0, 0, 1, 1, 0],[1, 0, 1, 0, 1],[1, 0, 1, 1, 0],[2, 0, 1, 0, 1]],
[[1, 0, 0, 0, 1],[1, 0, 0, 1, 0],[1, 0, 1, 0, 1],[1, 0, 1, 1, 0]],
[[1, 0, 0, 1, 1],[1, 0, 1, 1, 1],[3, 0, 0, 0, 1],[3, 0, 1, 0, 1]],
[[2, 0, 0, 1, 1],[2, 0, 1, 1, 0],[2, 0, 0, 1, 1],[2, 0, 1, 1, 0]],
[[1, 0, 1, 1, 1],[2, 0, 1, 1, 1],[2, 0, 0, 0, 1],[3, 0, 0, 0, 1]],
[[0, 0, 1, 1, 1],[2, 0, 0, 0, 1],[3, 0, 0, 1, 1],[3, 0, 1, 0, 1]],
[[0, 0, 1, 0, 1],[3, 0, 1, 1, 0],[3, 0, 0, 1, 0],[3, 0, 1, 1, 0],[3, 0, 1, 1, 1]],
[[2, 0, 0, 1, 0],[2, 0, 0, 1, 0],[2, 0, 1, 0, 1],[3, 0, 0, 1, 0],[3, 0, 1, 1, 1]]
]
A partition for sets of sizes:  1*3  8*4  2*5

[4, 2, 2, 2, 2]
[10, 0, 3]
[
[[0, 0, 1, 0, 1],[0, 0, 0, 1, 0],[0, 0, 1, 1, 1]],
[[0, 0, 0, 0, 1],[1, 0, 0, 1, 0],[3, 0, 0, 1, 1]],
[[0, 0, 1, 0, 1],[1, 0, 0, 1, 1],[3, 0, 1, 1, 0]],
[[1, 0, 0, 0, 1],[1, 0, 1, 1, 0],[2, 0, 1, 1, 1]],
[[1, 0, 1, 0, 1],[1, 0, 1, 1, 1],[2, 0, 0, 1, 0]],
[[1, 0, 0, 0, 1],[1, 0, 0, 1, 0],[2, 0, 0, 1, 1]],
[[1, 0, 1, 1, 0],[1, 0, 1, 1, 1],[2, 0, 0, 0, 1]],
[[0, 0, 1, 1, 0],[1, 0, 0, 1, 1],[3, 0, 1, 0, 1]],
[[2, 0, 1, 0, 1],[3, 0, 0, 1, 0],[3, 0, 1, 1, 1]],
[[2, 0, 0, 0, 1],[3, 0, 0, 1, 1],[3, 0, 0, 1, 0]],
[[2, 0, 1, 0, 1],[2, 0, 1, 1, 0],[2, 0, 1, 1, 1],[3, 0, 1, 0, 1],[3, 0, 0, 0, 1]],
[[0, 0, 1, 1, 0],[1, 0, 1, 0, 1],[2, 0, 1, 1, 0],[2, 0, 0, 1, 0],[3, 0, 1, 1, 1]],
[[0, 0, 1, 1, 1],[0, 0, 0, 1, 1],[2, 0, 0, 1, 1],[3, 0, 0, 0, 1],[3, 0, 1, 1, 0]]
]
A partition for sets of sizes: 10*3  0*4  3*5
\end{verbatim}

\newpage

\begin{verbatim} 
[4, 2, 2, 2, 2]
[6, 3, 3]
[
[[0, 0, 1, 0, 1],[0, 0, 0, 1, 0],[0, 0, 1, 1, 1]],
[[0, 0, 0, 0, 1],[1, 0, 0, 1, 0],[3, 0, 0, 1, 1]],
[[0, 0, 1, 0, 1],[1, 0, 0, 1, 1],[3, 0, 1, 1, 0]],
[[1, 0, 0, 0, 1],[1, 0, 1, 1, 0],[2, 0, 1, 1, 1]],
[[1, 0, 1, 0, 1],[1, 0, 1, 1, 1],[2, 0, 0, 1, 0]],
[[1, 0, 0, 0, 1],[1, 0, 0, 1, 0],[2, 0, 0, 1, 1]],
[[1, 0, 1, 0, 1],[1, 0, 1, 1, 0],[3, 0, 0, 0, 1],[3, 0, 0, 1, 0]],
[[0, 0, 1, 1, 0],[2, 0, 1, 0, 1],[3, 0, 1, 0, 1],[3, 0, 1, 1, 0]],
[[0, 0, 1, 1, 0],[2, 0, 0, 0, 1],[3, 0, 0, 1, 0],[3, 0, 1, 0, 1]],
[[2, 0, 1, 0, 1],[2, 0, 1, 1, 0],[2, 0, 1, 1, 1],[3, 0, 0, 1, 1],[3, 0, 1, 1, 1]],
[[0, 0, 0, 1, 1],[1, 0, 0, 1, 1],[2, 0, 0, 0, 1],[2, 0, 1, 1, 0],[3, 0, 1, 1, 1]],
[[0, 0, 1, 1, 1],[1, 0, 1, 1, 1],[2, 0, 0, 1, 0],[2, 0, 0, 1, 1],[3, 0, 0, 0, 1]]
]
A partition for sets of sizes:  6*3  3*4  3*5

[4, 2, 2, 2, 2]
[2, 6, 3]
[
[[0, 0, 1, 0, 1],[0, 0, 0, 1, 0],[0, 0, 1, 1, 1]],
[[0, 0, 0, 0, 1],[1, 0, 0, 1, 0],[3, 0, 0, 1, 1]],
[[0, 0, 1, 0, 1],[0, 0, 1, 1, 0],[1, 0, 0, 0, 1],[3, 0, 0, 1, 0]],
[[0, 0, 1, 1, 0],[1, 0, 0, 1, 1],[1, 0, 1, 1, 0],[2, 0, 0, 1, 1]],
[[1, 0, 1, 1, 1],[1, 0, 0, 0, 1],[3, 0, 0, 0, 1],[3, 0, 1, 1, 1]],
[[1, 0, 1, 0, 1],[1, 0, 1, 1, 0],[3, 0, 1, 0, 1],[3, 0, 1, 1, 0]],
[[2, 0, 0, 1, 0],[2, 0, 1, 0, 1],[2, 0, 0, 0, 1],[2, 0, 1, 1, 0]],
[[1, 0, 1, 0, 1],[2, 0, 1, 1, 1],[2, 0, 0, 1, 1],[3, 0, 0, 0, 1]],
[[2, 0, 1, 1, 0],[2, 0, 1, 0, 1],[2, 0, 1, 1, 1],[3, 0, 0, 1, 1],[3, 0, 1, 1, 1]],
[[0, 0, 1, 1, 1],[1, 0, 0, 1, 1],[1, 0, 1, 1, 1],[3, 0, 1, 0, 1],[3, 0, 1, 1, 0]],
[[0, 0, 0, 1, 1],[1, 0, 0, 1, 0],[2, 0, 0, 0, 1],[2, 0, 0, 1, 0],[3, 0, 0, 1, 0]]
]
A partition for sets of sizes:  2*3  6*4  3*5

[4, 2, 2, 2, 2]
[7, 1, 4]
[
[[0, 0, 1, 0, 1],[0, 0, 0, 1, 0],[0, 0, 1, 1, 1]],
[[0, 0, 0, 0, 1],[1, 0, 0, 1, 0],[3, 0, 0, 1, 1]],
[[0, 0, 1, 0, 1],[1, 0, 0, 1, 1],[3, 0, 1, 1, 0]],
[[1, 0, 0, 0, 1],[1, 0, 1, 1, 0],[2, 0, 1, 1, 1]],
[[1, 0, 1, 0, 1],[1, 0, 1, 1, 1],[2, 0, 0, 1, 0]],
[[1, 0, 0, 0, 1],[1, 0, 0, 1, 0],[2, 0, 0, 1, 1]],
[[1, 0, 1, 1, 0],[1, 0, 1, 1, 1],[2, 0, 0, 0, 1]],
[[1, 0, 1, 0, 1],[2, 0, 1, 0, 1],[2, 0, 0, 0, 1],[3, 0, 0, 0, 1]],
[[2, 0, 1, 1, 0],[2, 0, 0, 1, 0],[2, 0, 0, 1, 1],[3, 0, 0, 1, 0],[3, 0, 1, 0, 1]],
[[2, 0, 1, 0, 1],[2, 0, 1, 1, 0],[2, 0, 1, 1, 1],[3, 0, 0, 1, 1],[3, 0, 1, 1, 1]],
[[0, 0, 1, 1, 0],[3, 0, 1, 1, 0],[3, 0, 1, 1, 1],[3, 0, 0, 1, 0],[3, 0, 1, 0, 1]],
[[0, 0, 1, 1, 1],[0, 0, 0, 1, 1],[0, 0, 1, 1, 0],[1, 0, 0, 1, 1],[3, 0, 0, 0, 1]]
]
A partition for sets of sizes:  7*3  1*4  4*5
\end{verbatim}

\newpage

\begin{verbatim} 
[4, 2, 2, 2, 2]
[3, 4, 4]
[
[[0, 0, 1, 0, 1],[0, 0, 0, 1, 0],[0, 0, 1, 1, 1]],
[[0, 0, 0, 0, 1],[1, 0, 0, 1, 0],[3, 0, 0, 1, 1]],
[[0, 0, 1, 0, 1],[1, 0, 0, 1, 1],[3, 0, 1, 1, 0]],
[[1, 0, 0, 0, 1],[1, 0, 1, 0, 1],[1, 0, 0, 0, 1],[1, 0, 1, 0, 1]],
[[1, 0, 1, 1, 0],[1, 0, 1, 1, 1],[1, 0, 0, 1, 0],[1, 0, 0, 1, 1]],
[[0, 0, 1, 1, 0],[1, 0, 1, 1, 0],[1, 0, 1, 1, 1],[2, 0, 1, 1, 1]],
[[2, 0, 0, 0, 1],[2, 0, 0, 1, 0],[2, 0, 1, 0, 1],[2, 0, 1, 1, 0]],
[[2, 0, 0, 1, 1],[2, 0, 1, 1, 1],[2, 0, 0, 0, 1],[3, 0, 0, 1, 0],[3, 0, 1, 1, 1]],
[[0, 0, 1, 1, 0],[3, 0, 0, 0, 1],[3, 0, 0, 1, 1],[3, 0, 0, 0, 1],[3, 0, 1, 0, 1]],
[[2, 0, 0, 1, 0],[2, 0, 0, 1, 1],[2, 0, 1, 1, 0],[3, 0, 1, 0, 1],[3, 0, 0, 1, 0]],
[[0, 0, 1, 1, 1],[0, 0, 0, 1, 1],[2, 0, 1, 0, 1],[3, 0, 1, 1, 0],[3, 0, 1, 1, 1]]
]
A partition for sets of sizes:  3*3  4*4  4*5

[4, 2, 2, 2, 2]
[4, 2, 5]
[
[[0, 0, 1, 0, 1],[0, 0, 0, 1, 0],[0, 0, 1, 1, 1]],
[[0, 0, 0, 0, 1],[1, 0, 0, 1, 0],[3, 0, 0, 1, 1]],
[[0, 0, 1, 0, 1],[1, 0, 0, 1, 1],[3, 0, 1, 1, 0]],
[[1, 0, 0, 0, 1],[1, 0, 1, 1, 0],[2, 0, 1, 1, 1]],
[[1, 0, 1, 0, 1],[1, 0, 1, 1, 1],[1, 0, 0, 0, 1],[1, 0, 0, 1, 1]],
[[1, 0, 0, 1, 0],[1, 0, 1, 0, 1],[3, 0, 0, 0, 1],[3, 0, 1, 1, 0]],
[[1, 0, 1, 1, 1],[2, 0, 0, 0, 1],[3, 0, 0, 1, 0],[3, 0, 0, 1, 1],[3, 0, 1, 1, 1]],
[[2, 0, 1, 1, 0],[2, 0, 0, 0, 1],[2, 0, 0, 1, 1],[3, 0, 0, 0, 1],[3, 0, 1, 0, 1]],
[[0, 0, 1, 1, 0],[2, 0, 0, 1, 0],[2, 0, 0, 1, 1],[2, 0, 0, 1, 0],[2, 0, 1, 0, 1]],
[[0, 0, 1, 1, 0],[1, 0, 1, 1, 0],[2, 0, 1, 0, 1],[2, 0, 1, 1, 1],[3, 0, 0, 1, 0]],
[[0, 0, 1, 1, 1],[0, 0, 0, 1, 1],[2, 0, 1, 1, 0],[3, 0, 1, 0, 1],[3, 0, 1, 1, 1]]
]
A partition for sets of sizes:  4*3  2*4  5*5

[4, 2, 2, 2, 2]
[0, 5, 5]
[
[[0, 0, 1, 0, 1],[0, 0, 1, 1, 0],[0, 0, 0, 0, 1],[0, 0, 0, 1, 0]],
[[0, 0, 1, 1, 1],[0, 0, 0, 1, 1],[1, 0, 0, 0, 1],[3, 0, 1, 0, 1]],
[[0, 0, 1, 1, 1],[1, 0, 0, 1, 0],[1, 0, 0, 1, 1],[2, 0, 1, 1, 0]],
[[1, 0, 1, 0, 1],[1, 0, 1, 1, 0],[1, 0, 0, 0, 1],[1, 0, 0, 1, 0]],
[[1, 0, 1, 1, 1],[1, 0, 0, 1, 1],[3, 0, 0, 1, 0],[3, 0, 1, 1, 0]],
[[1, 0, 1, 1, 1],[2, 0, 0, 0, 1],[3, 0, 0, 0, 1],[3, 0, 1, 0, 1],[3, 0, 0, 1, 0]],
[[2, 0, 1, 1, 0],[2, 0, 1, 1, 1],[2, 0, 0, 0, 1],[3, 0, 0, 1, 1],[3, 0, 0, 1, 1]],
[[0, 0, 1, 0, 1],[1, 0, 1, 0, 1],[2, 0, 0, 1, 0],[2, 0, 1, 0, 1],[3, 0, 1, 1, 1]],
[[2, 0, 0, 1, 0],[2, 0, 0, 1, 1],[2, 0, 1, 1, 1],[3, 0, 1, 1, 1],[3, 0, 0, 0, 1]],
[[0, 0, 1, 1, 0],[1, 0, 1, 1, 0],[2, 0, 0, 1, 1],[2, 0, 1, 0, 1],[3, 0, 1, 1, 0]]
]
A partition for sets of sizes:  0*3  5*4  5*5
\end{verbatim}

\newpage

\begin{verbatim} 
[4, 2, 2, 2, 2]
[5, 0, 6]
[
[[0, 0, 1, 0, 1],[0, 0, 0, 1, 0],[0, 0, 1, 1, 1]],
[[0, 0, 0, 0, 1],[1, 0, 0, 1, 0],[3, 0, 0, 1, 1]],
[[0, 0, 1, 0, 1],[1, 0, 0, 1, 1],[3, 0, 1, 1, 0]],
[[1, 0, 0, 0, 1],[1, 0, 1, 1, 0],[2, 0, 1, 1, 1]],
[[1, 0, 1, 0, 1],[1, 0, 1, 1, 1],[2, 0, 0, 1, 0]],
[[1, 0, 0, 0, 1],[1, 0, 0, 1, 0],[1, 0, 0, 1, 1],[2, 0, 0, 0, 1],[3, 0, 0, 0, 1]],
[[1, 0, 1, 1, 1],[2, 0, 0, 1, 1],[3, 0, 0, 0, 1],[3, 0, 0, 1, 0],[3, 0, 1, 1, 1]],
[[2, 0, 1, 1, 0],[2, 0, 0, 0, 1],[2, 0, 0, 1, 0],[3, 0, 0, 1, 1],[3, 0, 1, 1, 0]],
[[0, 0, 1, 1, 0],[1, 0, 1, 0, 1],[2, 0, 1, 1, 0],[2, 0, 1, 1, 1],[3, 0, 0, 1, 0]],
[[0, 0, 0, 1, 1],[1, 0, 1, 1, 0],[2, 0, 1, 0, 1],[2, 0, 1, 0, 1],[3, 0, 1, 0, 1]],
[[0, 0, 1, 1, 1],[0, 0, 1, 1, 0],[2, 0, 0, 1, 1],[3, 0, 1, 0, 1],[3, 0, 1, 1, 1]]
]
A partition for sets of sizes:  5*3  0*4  6*5

[4, 2, 2, 2, 2]
[1, 3, 6]
[
[[0, 0, 1, 0, 1],[0, 0, 0, 1, 0],[0, 0, 1, 1, 1]],
[[0, 0, 0, 0, 1],[0, 0, 0, 1, 1],[1, 0, 0, 0, 1],[3, 0, 0, 1, 1]],
[[0, 0, 1, 1, 0],[1, 0, 0, 1, 0],[1, 0, 0, 1, 1],[2, 0, 1, 1, 1]],
[[0, 0, 1, 1, 0],[1, 0, 1, 0, 1],[1, 0, 1, 1, 0],[2, 0, 1, 0, 1]],
[[1, 0, 0, 0, 1],[1, 0, 0, 1, 0],[1, 0, 0, 1, 1],[2, 0, 0, 0, 1],[3, 0, 0, 0, 1]],
[[1, 0, 1, 1, 1],[2, 0, 0, 1, 0],[3, 0, 0, 0, 1],[3, 0, 0, 1, 0],[3, 0, 1, 1, 0]],
[[2, 0, 1, 1, 0],[2, 0, 1, 1, 1],[2, 0, 0, 0, 1],[3, 0, 1, 0, 1],[3, 0, 1, 0, 1]],
[[0, 0, 1, 1, 1],[1, 0, 1, 0, 1],[2, 0, 0, 1, 1],[2, 0, 1, 1, 0],[3, 0, 1, 1, 1]],
[[0, 0, 1, 0, 1],[1, 0, 1, 1, 1],[1, 0, 1, 1, 0],[3, 0, 0, 1, 1],[3, 0, 1, 1, 1]],
[[2, 0, 0, 1, 0],[2, 0, 0, 1, 1],[2, 0, 1, 0, 1],[3, 0, 1, 1, 0],[3, 0, 0, 1, 0]]
]
A partition for sets of sizes:  1*3  3*4  6*5

[4, 2, 2, 2, 2]
[2, 1, 7]
[
[[0, 0, 1, 0, 1],[0, 0, 0, 1, 0],[0, 0, 1, 1, 1]],
[[0, 0, 0, 0, 1],[1, 0, 0, 1, 0],[3, 0, 0, 1, 1]],
[[0, 0, 1, 0, 1],[0, 0, 1, 1, 0],[1, 0, 0, 0, 1],[3, 0, 0, 1, 0]],
[[0, 0, 1, 1, 0],[1, 0, 0, 1, 1],[1, 0, 1, 1, 0],[1, 0, 0, 0, 1],[1, 0, 0, 1, 0]],
[[1, 0, 1, 0, 1],[1, 0, 1, 1, 1],[1, 0, 0, 1, 1],[2, 0, 0, 1, 1],[3, 0, 0, 1, 0]],
[[1, 0, 1, 1, 1],[2, 0, 0, 0, 1],[3, 0, 0, 0, 1],[3, 0, 1, 1, 0],[3, 0, 0, 0, 1]],
[[2, 0, 1, 1, 0],[2, 0, 1, 1, 1],[2, 0, 0, 0, 1],[3, 0, 1, 0, 1],[3, 0, 1, 0, 1]],
[[2, 0, 1, 0, 1],[2, 0, 1, 1, 0],[2, 0, 1, 1, 1],[3, 0, 0, 1, 1],[3, 0, 1, 1, 1]],
[[1, 0, 1, 0, 1],[1, 0, 1, 1, 0],[2, 0, 0, 1, 0],[2, 0, 0, 1, 0],[2, 0, 0, 1, 1]],
[[0, 0, 1, 1, 1],[0, 0, 0, 1, 1],[2, 0, 1, 0, 1],[3, 0, 1, 1, 1],[3, 0, 1, 1, 0]]
]
A partition for sets of sizes:  2*3  1*4  7*5

[4, 2, 2, 2, 2]
[0, 0, 9]
[
[[0, 0, 1, 0, 1],[0, 0, 1, 1, 0],[0, 0, 1, 1, 1],[0, 0, 0, 0, 1],[0, 0, 1, 0, 1]],
[[0, 0, 0, 1, 1],[0, 0, 1, 1, 0],[0, 0, 1, 1, 1],[1, 0, 0, 0, 1],[3, 0, 0, 1, 1]],
[[1, 0, 0, 1, 0],[1, 0, 0, 1, 1],[1, 0, 1, 0, 1],[2, 0, 0, 0, 1],[3, 0, 1, 0, 1]],
[[1, 0, 0, 0, 1],[1, 0, 0, 1, 0],[1, 0, 0, 1, 1],[2, 0, 0, 1, 0],[3, 0, 0, 1, 0]],
[[1, 0, 1, 1, 1],[2, 0, 0, 1, 1],[3, 0, 0, 0, 1],[3, 0, 0, 1, 0],[3, 0, 1, 1, 1]],
[[2, 0, 1, 1, 0],[2, 0, 1, 1, 1],[2, 0, 0, 0, 1],[3, 0, 1, 1, 0],[3, 0, 1, 1, 0]],
[[1, 0, 1, 0, 1],[2, 0, 1, 1, 0],[3, 0, 1, 1, 1],[3, 0, 0, 0, 1],[3, 0, 1, 0, 1]],
[[1, 0, 1, 1, 1],[1, 0, 1, 1, 0],[2, 0, 0, 1, 1],[2, 0, 1, 0, 1],[2, 0, 1, 1, 1]],
[[0, 0, 0, 1, 0],[1, 0, 1, 1, 0],[2, 0, 1, 0, 1],[2, 0, 0, 1, 0],[3, 0, 0, 1, 1]]
]
A partition for sets of sizes:  0*3  0*4  9*5
\end{verbatim}

\section{Program to check a partition}\label{a:program}
\normalsize
We offer a simple program in Python 3 that allows to check easily if a partition is a zero-sum partition. It can be executed in a terminal (or in an online Python environment like https://www.online-python.com/ or https://trinket.io/python), with the three elements of the input: description of the group, sizes of the sets in the partition, and the partition itself, copied-pasted from the annexes.
\scriptsize
\begin{verbatim}
import json

group = json.loads(input())
sizes = json.loads(input())
sets = json.loads(input())
ok = True
for set in sets:
    sums = [0 for pos in range(len(group))]
    for elem in set:
        for pos in range(len(group)):
            sums[pos] = (sums[pos] + elem[pos]) % group[pos]
    if not sums == [0 for i in range(len(group))]:
        ok = False
        break
if ok:
    print("Zero-sum partition")
else:
    print("Not a zero-sum partition")
\end{verbatim}

\end{appendices}

\end{document}